\documentclass[11pt]{article}
\usepackage{fullpage}
\usepackage{amsmath,amssymb,amsthm}
\usepackage{stmaryrd}
\usepackage[english]{babel}
\usepackage[left=3cm,right=3cm,top=3cm,bottom=4cm]{geometry}

\usepackage{tikz,tikz-cd}
\usetikzlibrary{fadings, calc, backgrounds}
\usetikzlibrary{snakes}
\usetikzlibrary{decorations.pathmorphing,shapes}
\usetikzlibrary{positioning}
\usepackage{braket}
\usepackage{mathtools,stmaryrd}

\SetSymbolFont{stmry}{bold}{U}{stmry}{m}{n}
\usepackage{lmodern,bm,bbm,mathrsfs}
\usepackage{manfnt}

\usepackage[shortlabels]{enumitem} 

\usepackage[scr=esstix]{mathalfa}

\usepackage{hyperref}
\usepackage[textsize=tiny]{todonotes}

\numberwithin{equation}{section} 

\definecolor{steinGray}{gray}{0.9}
\definecolor{ballGray}{gray}{0.95}
\definecolor{curveGray}{gray}{0.2}
\definecolor{textGray}{gray}{0.3}

\newtheorem*{theorem}{Theorem}
\newtheorem*{corollary}{Corollary}
\newtheorem*{lemma}{Lemma}
\newtheorem*{proposition}{Proposition}
\newtheorem*{conjecture}{Conjecture}

\theoremstyle{definition}
\newtheorem*{definition}{Definition}
\newtheorem*{remark}{Remark}

\newtheorem*{example}{Example}

\renewcommand{\bar}[1]{\overline{#1}}

\renewcommand{\hat}[1]{\widehat{#1}}

\renewcommand{\tilde}[1]{\widetilde{#1}}

\newcommand{\bC}{\mathbb{C}}

\newcommand{\bD}{\mathbb{D}}

\newcommand{\bF}{\mathbb{F}}

\newcommand{\bP}{\mathbb{P}}
\newcommand{\bQ}{\mathbb{Q}}
\newcommand{\cQ}{\mathcal{Q}}
\newcommand{\bR}{\mathbb{R}}

\newcommand{\bZ}{\mathbb{Z}}
\newcommand{\cA}{\mathcal{A}}
\newcommand{\cC}{\mathcal{C}}
\newcommand{\rd}{\mathrm{d}}

\newcommand{\cF}{\mathcal{F}}

\newcommand{\cH}{\mathcal{H}}

\newcommand{\cI}{\mathcal{I}}
\newcommand{\cJ}{\mathcal{J}}

\newcommand{\cL}{\mathcal{L}}
\newcommand{\cM}{\mathcal{M}}
\newcommand{\cN}{\mathcal{N}}
\newcommand{\bfN}{\mathbf{N}}

\newcommand{\cO}{\mathcal{O}}

\newcommand{\rOb}{\mathrm{Ob}}
\newcommand{\cP}{\mathcal{P}}

\newcommand{\cU}{\mathcal{U}}

\newcommand{\cX}{\mathcal{X}}

\newcommand{\cZ}{\mathcal{Z}}

\newcommand{\bfM}{\mathbf{M}}
\newcommand{\bfm}{\mathbf{m}}
\newcommand{\fm}{\mathfrak{m}}

\newcommand{\pt}{\mathrm{pt}}

\newcommand{\st}{\mathrm{st}}

\newcommand{\spec}{\mathrm{Spec}}

\newcommand{\vir}{\mathrm{vir}}
\newcommand{\sr}{\mathrm{sr}}

\newcommand{\Tr}{\mathrm{Tr}}
\newcommand{\red}{\mathrm{red}}
\newcommand{\Cpx}{\mathtt{Cpx}}
\newcommand{\CY}{\mathtt{CY}}
\newcommand{\BM}{\mathrm{BM}}
\newcommand{\rCY}{\mathrm{CY}}
\newcommand{\PT}{\mathrm{PT}}

\newcommand{\Tot}{\mathrm{Tot}}
\newcommand{\ord}{\mathrm{ord}}

\newcommand{\cok}{\mathrm{cok}}
\newcommand{\con}{\mathrm{con}}
\newcommand{\GW}{\mathrm{GW}}
\newcommand{\sF}{\mathrm{sF}}
\newcommand{\sm}{\mathrm{sm}}

\newcommand{\ev}{\mathrm{ev}}
\newcommand{\eff}{\mathrm{eff}}

\newcommand{\Aff}{\mathrm{Aff}}
\newcommand{\Aut}{\mathrm{Aut}}
\newcommand{\Herm}{\mathrm{Herm}}
\newcommand{\rmell}{\mathrm{ell}}
\newcommand{\divergence}{\mathrm{div}}
\newcommand{\sdiv}{\mathfrak{div}}
\newcommand{\proj}{\mathrm{proj}}

\newcommand{\tilchi}{\widetilde{\chi}}

\DeclareMathOperator{\Ext}{Ext}

\DeclareMathOperator{\ch}{ch}

\DeclareMathOperator{\Hom}{Hom}
\DeclareMathOperator{\id}{id}
\DeclareMathOperator{\im}{im}

\DeclareMathOperator{\pr}{pr}

\newcommand{\abs}[1]{\left| #1 \right|}
\DeclarePairedDelimiterX{\lseries}[1]{(}{)}{\mkern-2mu\delimsize(#1\delimsize)\mkern-2mu}
\DeclarePairedDelimiterX{\pseries}[1]{[}{]}{\mkern-2mu\delimsize[#1\delimsize]\mkern-2mu}

\makeatletter
\newcommand{\colim@}[2]{%
  \vtop{\m@th\ialign{##\cr
    \hfil$#1\operator@font colim$\hfil\cr
    \noalign{\nointerlineskip\kern1.5\ex@}#2\cr
    \noalign{\nointerlineskip\kern-\ex@}\cr}}%
}
\newcommand{\colim}{%
  \mathop{\mathpalette\colim@{\rightarrowfill@\scriptscriptstyle}}\nmlimits@
}
\renewcommand{\varprojlim}{%
  \mathop{\mathpalette\varlim@{\leftarrowfill@\scriptscriptstyle}}\nmlimits@
}
\makeatother

\tikzset{%
  vertex/.style={shape=circle,fill=black,minimum size=6pt,inner sep=0},
  framing/.style={shape=rectangle,fill=black,minimum size=6pt,inner sep=0},
  baseline={([yshift=-0.8ex]current bounding box.center)}
}

\title{Universal Curve Counting in Calabi--Yau $4$-folds}
\author{Ivan Karpov, Felix Thimm}
\date{\today}

\begin{document}

\maketitle

\begin{abstract}
Recently, Pardon introduced a universal curve counting framework, which allows him to show that curve counting on $3$-folds is ``generated'' by local curves (in a precise sense). We implement this universal curve counting program for Calabi--Yau $4$-folds, where a significant technical difficulty comes from the dependence of the virtual dimension on the genus, which leads to more complicated generators than just local curves. 

Our main techniques are combinations of deformation-theoretic and homotopy-theoretic methods with an adaptation of Pardon's generic transversality results, which allow us to establish a full set of (computable) generators for CY$4$ curve counting.
 
As a result, we prove (a slight correction of) the Gromov--Witten/Pandharipande--Thomas correspondence with primary descendents for CY$4$-folds, formulated by Cao, Maulik, and Toda. From this, we additionally obtain integrality of Gopakumar--Vafa style invariants for CY$4$-folds, as conjectured by Klemm--Pandharipande.
\end{abstract}

\setcounter{tocdepth}{2}
\tableofcontents

\section{Introduction}

\subsection{Curve Counting on CY4-folds}

\subsubsection{}
Fix a complex projective algebraic variety $X$. In enumerative geometry, we study certain invariants of moduli spaces of certain objects attached to $X$. In this work, we work with curve counting invariants, which are defined using natural compactifications of the moduli space of smooth embedded curves in $X$.

One of the natural compactifications is the moduli space of possibly disconnected stable maps
\begin{equation*}
    \bar{\cM}_{g,n}(X,\beta)\coloneqq \left\{f:(C,\mathbf{p})\to X\text{ stable}\ |\ \substack{C\text{ nodal of genus }g,\\ f_*[C]=\beta}\right\}/\cong,
\end{equation*}
where $\beta\in H_2(X)$ is a given curve class, see \cite{Kontsevich1995, Behrend-Manin-Mgn}. The analogous space with only connected source curves is denoted by $\bar{\cM}_{g,n}^\con(X,\beta)$. By \cite{Behrend1997GW}, this space carries a perfect obstruction theory in the sense of \cite{bf97_intrinsic_normal_cone}, which then defines a virtual fundamental class $\big[\bar{\cM}^\con_{g,n}(X,\beta)\big]^\vir$ of homological degree equal to the virtual dimension
\begin{equation}\label{eq:gw-vdim}
    -K_X\cdot\beta\,+\, (\dim(X)-3)(1-g)\,+\,n.
\end{equation}
Since we deal with the Calabi--Yau $4$-fold case in this paper, this simplifies to $1-g+n$, and we only study the genus $0$ and $1$ cases. Given integral classes $\gamma_i\in H^{m_i}(X,\bZ)$ for $1\leq i\leq r$ satisfying $\sum_i (m_i-2)=2$, we define the \textit{primary Gromov--Witten (GW) invariants}
\begin{equation}\label{intro:eq:raw-primary-gw}
    \begin{aligned}
        \GW_{0,\beta}(\gamma_\bullet)&\coloneqq \int_{\big[\bar{\cM}^\con_{0,n}(X,\beta)\big]^\vir}\prod_{i=1}^r\ev_i^*(\gamma_i)\in\bQ,\\
        \GW_{1,\beta}&\coloneqq \int_{\big[\bar{\cM}^\con_{1,0}(X,\beta)\big]^\vir}1\in \bQ\\
    \end{aligned}
\end{equation}
by integrating over the virtual fundamental class. 

\subsubsection{}
A problem with GW invariants is the contribution of multiple covers. So, instead of studying the GW invariants as defined above, one instead considers \textit{Gopakumar--Vafa type invariants}, which correct for these multiple cover contributions, and in the CY$4$ case for certain meeting contributions. So, fix an ample divisor \(H\), and complete all
series in effective curve classes by the degree \(H\cdot\beta\).  We write
\(\beta>0\) when \(\beta\) is nonzero and effective. Klemm and Pandharipande \cite{KP} define the rational numbers
\(n_{0,\beta}(\gamma_\bullet)\) and \(n_{1,\beta}\) by
\begin{equation}
 \begin{aligned}
 \sum_{\beta>0}\mathrm{GW}_{0,\beta}(\gamma_\bullet)Q^\beta
 &={}
 \sum_{\beta>0}n_{0,\beta}(\gamma_\bullet)
       \sum_{d\geq1}d^{r-3}Q^{d\beta},\\
 \sum_{\beta>0}\mathrm{GW}_{1,\beta}Q^\beta
 &={}
 \sum_{\beta>0}n_{1,\beta}
       \sum_{d\geq1}\frac{\sigma(d)}dQ^{d\beta}
 +\frac1{24}\sum_{\beta>0}n_{0,\beta}(c_2(TX))
       \log(1-Q^\beta)\\
 &\hspace{2em}
 -\frac1{24}\sum_{\beta_1,\beta_2>0}m_{\beta_1,\beta_2}
       \log(1-Q^{\beta_1+\beta_2}),
 \end{aligned}
 \label{eq:gv-kp-definitions}
\end{equation}
where the last sum in
\eqref{eq:gv-kp-definitions} is taken over ordered pairs, \(\sigma(d)=\sum_{a\mid d}a\), and the meeting numbers
\(m_{\beta_1,\beta_2}\) are determined by the rules recalled in
\eqref{app:eq:kp-meeting-rules}.

\subsubsection{}\label{intro:sec:PT}
Another natural way to view smooth embedded curves in $X$ is as smooth pure $1$-dimensional subvarieties given by ideal sheaves. One compactification of this is the moduli space of stable pairs
\begin{equation}\label{intro:eq:stable-pairs-moduli}
    P_n(X,\beta)\coloneqq \left\{[\cO_X\xrightarrow{s} \cF]\in D^b(X) |\ \substack{[\mathrm{supp}(\cF)]=\beta,\, \chi(\cF)=n\\ \cF\text{ pure $1$-dim.},\, \cok(s)\text{ $0$-dim.}}\right\}
\end{equation}
as considered in \cite{pt07,PandharipandeThomas2010}. If $X$ is a Fano or CY $3$-fold, this moduli space also carries a perfect obstruction theory, which yields a virtual class and hence invariants as above. However, if $X$ is a CY $4$-fold, as we consider in this paper, the moduli space of stable pairs is equipped with a symmetric $3$-term obstruction theory instead, which yields a virtual fundamental class  
\begin{equation*}
    \big[P_n(X,\beta)\big]^\vir\in H_{2n}(P_n(X,\beta);\bZ)
\end{equation*}
of virtual dimension $2n$, by \cite{Borisov_2017, oh2022counting}. Given the same integral classes $\gamma_i\in H^{m_i}(X,\bZ)$ for $1\leq i\leq r$ with $\sum_i (m_i-2)=2$ as above, this allows us to define the \textit{primary Pandharipande--Thomas (PT) invariants}
\begin{equation}\label{intro:eq:primary-pt}
    \begin{aligned}
        \PT_{0,\beta}&\coloneqq \int_{\big[\cP_{0}(X,\beta)\big]^\vir}1 \in\bZ,\\
        \PT_{1,\beta}(\gamma_\bullet)&\coloneqq \int_{\big[\cP_{1}(X,\beta)\big]^\vir}\prod_{i=1}^r(\pi_{\cP})_*\big(\pi_X^*(\gamma_i)\cup \ch_3(\bF)\big)\in \bZ,
    \end{aligned}
\end{equation}
where $\pi_\cP$ and $\pi_X$ are the projections of $X\times P_n(X,\beta)$, and $\cI=[\cO\to \bF]$ is the universal pair. Note that $\ch_3(\bF)$ is the Poincar\'e dual of the fundamental cycle of $\bF$.

\subsubsection{}\label{intro:sec:conjectures}
The GV type invariant $n_{0,\beta}(\gamma_\bullet)$ and $n_{1,\beta}$ are defined using the rational GW invariants. So, \`a priori, they are rational. However, their definition is chosen to correct for multiple cover and meeting contributions, so that Klemm and Pandharipande conjecture the following.
\begin{conjecture}[{\cite[Conj. 0 and 1]{KP}}]
    The invariants $n_{0,\beta}(\gamma_\bullet)$ and $n_{1,\beta}$ are integers.
\end{conjecture}

Since all curve-counting invariants above are defined using compactifications of the same space of smooth embedded curves, they essentially encode the same information. Hence, one naturally expects the invariants to be related in some way, which was made precise in the following PT/GV correspondence, conjectured by Cao, Maulik, and Toda.
\begin{conjecture}[{\cite[Conj. 0.1 and 0.2]{MCT}}]
    Take $\gamma_\bullet=(\tilde{\gamma}_\bullet,\delta_\bullet)$ with $\sum_i (m_i-2)=2$, where $\tilde{\gamma}_i$ have $m_i>2$, and $\delta_j$ have $m_j=2$.\footnote{Note that the order of insertions is not relevant, and any insertion with $m_i<2$ yields a vanishing expression by degree reasons.} For a suitable choice of orientation, we have
    \begin{equation}\label{intro:eq:cmt-conj}
        \begin{aligned}
            \PT_{1,\beta}(\gamma_\bullet)&=\big(\prod_j\textstyle\int_\beta \delta_j\big)\cdot \sum_{\substack{\beta_1+\beta_2=\beta\\ \beta_1,\beta_2\geq 0}}n_{0,\beta_1}(\tilde{\gamma}_\bullet)\PT_{0,\beta_2},\\
            \sum_{\beta\geq 0}\PT_{0,\beta}Q^\beta&=\prod_{\beta>0}M(Q^\beta)^{n_{1,\beta}},
        \end{aligned}
    \end{equation}
    where $n_{0,0}(\gamma_\bullet)\coloneqq 0$, $\PT_{0,0}\coloneqq 1$, and $M(x)=\prod_{k\geq 1}(1-x^k)^{-k}$ is the MacMahon function.
\end{conjecture}

\begin{remark}
    \begin{itemize}
        \item Note here that we introduced a slight correction for the case of divisor insertions compared to \cite[Conj. 0.1 and 0.2]{MCT}. We prove this corrected formula and show in Example~\ref{ex:mct-divisor-counterexample} that it indeed leads to different computations.
        \item Cao, Maulik, and Toda give numerical evidence for their GW/PT conjectures above (see \cite[Conj. 0.1 and 0.2]{MCT}). Note also that the case of local $\bP^2$ has been recently verified by Cao \cite{cao-local-p2}.
    \end{itemize}
\end{remark}

\subsubsection{}
Similar conjectures in the $3$-fold case were recently proven by Ionel--Parker \cite{Ionel2018} and Pardon \cite{pardon2023}. In order to improve Ionel-Parker's theory, Pardon introduced a ``ring of curve enumeration problems''. Deformation invariant curve counting theories pair with the element of this ring, and this pairing contains the information about certain enumerative invariants. 

Our main result is adapting his techniques to the CY $4$-fold setting and finding a complete set of computable generators for this ring. The main application is as follows.

\begin{theorem}\label{intro:thm:main-application}
   All Conjectures in \S\ref{intro:sec:conjectures} above hold for \textit{primary insertions}.
\end{theorem}

\begin{remark}
    While for us the integrality conjectures above are a consequence of the PT/GV correspondence, it is conceivable that one may obtain a more direct proof by combining our techniques of handling  higher-dimensional generators with Ionel--Parker's method. We do not pursue this direction here.
\end{remark}

\subsection{Pardon's Universal Curve Counting Framework}

\subsubsection{}
We summarize here briefly Pardon's universal counting framework for $3$-folds. For a family $X\to B$ of complex threefolds, let $\cZ(X/B)\to B$ denote its relative space of $1$-cycles. Pardon defines (co-)limits of bivariant (co-)homology groups
\begin{equation*}
        H_c^*(\cZ/\Cpx_3)\coloneqq \varinjlim_{X\to B}H_c^*(\cZ(X/B)/B),\quad H_*^\BM(\cZ/\Cpx_3) \coloneqq \varprojlim_{X\to B}H_*^\BM(\cZ(X/B)/B),
\end{equation*}
where the limits are taken, in a suitable sense, over all families of (possibly non-compact) $3$-folds. 

The first group is the Grothendieck group of $1$-cycles in complex threefolds, which we call the \textit{Pardon group of curve enumeration problems}, since its elements may be regarded as such. In particular, any projective $3$-fold $X$ with curve class $\beta$ and suitable classes $\gamma_\bullet$ defines a canonical element $(X,\beta,\gamma_\bullet)$ of this group which corresponds to counting curves on $X$ with insertions $\gamma_{\bullet}$, see \S\ref{subsubsec:primary-insertions} for details. In this paper, we work with  $H_c^*(\cZ/\CY_4)$ and $H_*^\BM(\cZ/\CY_4)$, which are analogously constructed in \S\ref{sec:unem}.

The second is \textit{Pardon's group of curve enumeration theories}. By deformation invariance, GW- and PT-theories define natural elements in this group (using pushforwards of virtual classes along each support map to the space of $1$-cycles), see \S\ref{subsec:normalized-family-theories} for details. The natural pairing between the two groups evaluates an enumeration theory on an enumeration problem, recovering the usual GW- and PT-invariants as illustrated below.
\begin{equation*}
\begin{tikzpicture}[
    grp/.style={inner sep=2pt, text height=1.8ex, text depth=.7ex},
    squig/.style={-{cm to}, line join=round,
        decorate, decoration={zigzag, segment length=4, amplitude=.9,
                              pre=lineto, pre length=2pt,
                              post=lineto, post length=3pt}},
  ]
  \node[grp] (HL) at (-2,0) {$H^\BM_{\mathrm{rvd}}(\cZ/\Cpx_3)$};
  \node[grp] (HR) at ( 2,0) {$H^*_c(\cZ/\Cpx_3)$};
  \node[inner sep=0pt, anchor=base east] at (HL.base west) {$
      \begin{gathered}
        (\rho_\GW)_*\big[\overline{\cM}_g\big]^\vir\\
        (\rho_\PT)_*\big[P_n\big]^\vir
      \end{gathered}\;\in\;$};
  \node[inner sep=0pt, anchor=base west] at (HR.base east)
      {$\;\ni\;(X,\beta;\gamma_\bullet)$};
  \node (P) at (0,-1.6) {$\langle -,-\rangle$};
  \node (I) at (0,-3)   {GW, PT invariants};
  \draw[squig] (HL.south) -- (P);
  \draw[squig] (HR.south) -- (P);
  \draw[squig] (P) -- (I);
\end{tikzpicture}
\end{equation*}
Here, $rvd$ denotes the real virtual dimension. Hence, we can understand any deformation invariant curve counting theory by studying $H_c^*(\cZ/\Cpx_3)$ in cohomological degree equal to the real virtual dimension.

\subsubsection{}
A rough version of Pardon's main result is then.
\begin{theorem}[{{\cite[Theorem 1.2]{pardon2023}}}]
    The Grothendieck group of semi-Fano $1$-cycles in complex threefolds $H_c^*(\cZ_\sF/\Cpx_3;\bQ)$, in cohomological degree equal to the real virtual dimension, is freely generated as a $\bQ$-algebra by elements $x_{g,m,k}$ corresponding to equivariant local curves with genus $g\geq 0$, cycle multiplicity $m\geq 1$, Chern number $k\geq 0$, and $(m-1)k=0$.
\end{theorem}

The main technique Pardon uses for this theorem is a technical regularization procedure called generic transversality (see \S\ref{sec:regdef}), which, roughly, deforms a variety (or a family thereof) so that, after a stratification, the cycles under consideration have a smooth chart $S$ of expected dimension. Since the cycle space itself has ill-behaved deformation theory, the chart is constructed from deformations of associated stable maps, which can be obtained by normalizing irreducible components of the cycle. For the resulting object, the real virtual dimension matches this expected dimension, and the group can be generated by Poincar\'e duals of points in $S$, which Pardon then connects to the equivariant local curves elements $x_{g,m,k}$.

More precisely, Pardon's theorem formally expresses all relevant elements of his group as  polynomials in the elements $x_{g,m,0}$. Both $\GW$ and $\PT$ are multiplicative, and hence their values on all CY or Fano threefold problems are determined by their values on these generators. 

Thus, to prove $\GW/\PT$-correspondence in full generality, it is enough to check it for local curves. This particular case follows from the calculations of Bryan--Pandharipande and Okounkov--Pandharipande \cite{BryanPandharipande2008,OkounkovPandharipande2010}. Thus, Pardon establishes $\GW/\PT$ for all CY or Fano threefolds $X$.

\subsection{Theorems \& Techniques}

\subsubsection{}
In the CY$4$ setting, a major difficulty is that, in contrast to the $3$-fold case, the real virtual dimension and expected dimension of the charts $S$ above depend on the genus, which results in our having to consider much more complicated classes than simply Poincar\'e duals of points in $S$, see \eqref{inc:eq:associated-graded-generators}. That is, we have to work with classes of \textit{higher degrees.}

Our techniques, outlined just below, manage to work out explicit generators for these additional cohomological degrees, resulting in the theorem below. The main application, Theorem~\ref{intro:thm:main-application}, then follows as a consequence, after some numerical checks, which were essentially already explained in the existing literature, see \S\ref{liftschecks}.

\begin{theorem}[{Theorem~\ref{rcl:thm:indecomposable-quotient-generators}}]\label{intro:thm:reduction-to-gens}
The group
\begin{equation*}
 H_c^q(\mathcal Z_{\mathrm{con}}/\CY_4;\mathbb Q)
\end{equation*}
is spanned by the class $p_m$ of a local rational curve in degree $q=2$; vanishes in degree $q=1$; and in degree $q=0$ is spanned by
\begin{equation*}
    r_m,\quad e_m,\quad \ell_{d_1,d_2}^{\mathrm{con}},\qquad m,d_1,d_2\geq 1,
\end{equation*}
where $r_m$ is the class of a ``Bott family'' of local rational curves over $\bP^1$, $e_m$ is the class of a local elliptic curve, and $\ell_{d_1,d_2}^{\mathrm{con}}$ is the class of a nodal rational--rational meeting. See \S\ref{subsec:standard-connected-classes} for the definition of these classes.
\end{theorem}

\subsubsection{}
We now outline our main techniques. First, in \S\ref{sec:regdef}, we provide our perspective on Pardon's generic transversality technique, which made the proof more accessible to us as algebraic geometers. The base version of this perspective was initially developed during work of the second-named author with Jim Bryan on an orbifold version of Pardon's work. Here, we use this to prove CY preservation and establish modifications that help us work with the more complicated generators.

On the connected semi-regular locus, we can then Whitney stratify by multiplicity, topological type, and incidence type. We compute the dimension of semi-charts in the resulting strata to establish the following table of loci
\begin{equation*}
\begin{array}{c|c|c}
 \text{degree}&\text{locus}&\text{representative on a stratum}\\ \hline
 2&\text{smooth embedded rational curves}&\text{point}\\
 1&\text{smooth embedded rational curves}&\text{oriented loop}\\
 0&\text{smooth embedded rational curves}&\text{oriented surface}\\
 0&\text{smooth embedded elliptic curves}&\text{point}\\
 0&\text{two embedded rational curves meeting at one point}&\text{point}
\end{array}
\end{equation*}
that contribute when studying $H_c^j\left(\cZ_{\mathrm{con}}^{\mathrm{sr}\circ}/\CY_4\right)$. This leads us to consider techniques to deal with the two complicated generators: oriented surface classes of smooth embedded rational curves, and ordinary rational--rational meetings.

We study the locus of ordinary rational--rational meetings in \S\ref{inc:sec:proof}. There, we deform each such point to an explicit toric model using an analogue of the \textit{plumbing} procedure from low-dimensional topology which goes back to Mumford, Milnor and Hirzebruch (see, e.g., \cite{Mumford61}).

\subsubsection{}\label{intro:sec:rational-locus-gen}
To find generators for the locus of smooth rational embedded curves, in \S\ref{rat:sec:rational-normal-data}, we establish a surjective map from a bordism group
\begin{equation*}
    \Phi_{m,d}\colon
 \Omega_d^{\mathrm{SO}}(\mathcal N_0)\otimes\mathbb Q
 \twoheadrightarrow
 H_c^{2-d}(R_m^{\mathrm{sm}}/\CY_4;\mathbb Q),
\end{equation*}
where $\cN_0$ is roughly the classifying space of oriented smooth $S^2$-bundles with rank $3$ complex vector bundles subject to a CY-type condition on the total space, see Definition~\ref{def:N0-classifying-space} for details. The map is constructed by analytification followed by a regularizing deformation using generic transversality. Bordism invariance requires a version of generic transversality which allows fixing the geometry over a subcomplex in the base, over which sufficient regularity is already given. We prove surjectivity using deformation to the normal cone and the Atiyah--Hirzebruch spectral sequence.

With this map, a computation of the bordism groups in \S\ref{rat:sec:degree-two-normal-generator} immediately shows that each $H_c^{2-d}(R_m^{\mathrm{sm}}/\CY_4;\mathbb Q)$ has a single generator in degrees $0$ and $2$, and vanishes in degree $1$. This allows us to choose simple computable geometries as our generators $p_m$ and $r_m$.

\begin{remark}
    The techniques used here are reminiscent of the ones used in \cite{JU,karpov-cy4-orientations} to construct orientations. We find it interesting that such similar techniques, which reduce problems in complex geometry to the computation of certain bordism groups, make an appearance in these two very different contexts.
\end{remark}

\subsubsection{}
Proving the desired conjectures using Pardon's argument, requires GW and PT theory to define enumeration theories in Pardon's group. That is, one needs virtual classes on the moduli spaces over arbitrary analytic CY$4$-families $X\to B$ over simplicial bases, compatible with base change and restriction to open subfamilies. For GW theory, these are provided by Pardon's construction, see \S\ref{subsec:normalized-family-theories}. 

However, for PT theory on CY$4$-folds, there is an important caveat: since the obstruction theories involved are symmetric $3$-term obstruction theories, rather than perfect, the virtual is constructed differently, following \cite{Borisov_2017,oh2022counting}, and depends on a choice of orientation. The family version hence requires canonical orientations \textit{in families}. These were constructed by the first author in \cite{karpov-cy4-orientations}, building on \cite{JU}. The resulting family virtual classes are constructed in \cite{karpov-bj}, as an analytic, relative version of \cite{Borisov_2017}. We record this in Proposition~\ref{prop:family-cy4-virtual-package}. The same orientations also fix the signs in the computations of \cite{MCT} that we use in \S\ref{liftschecks}.

\subsubsection{Plan of the Paper}
In \S\ref{background}, we collect some background results in the theory of complex geometry over simplicial bases as used in Pardon's work. In \S\ref{sec:unem}, we introduce the Pardon rings of curve enumeration theories and problems, and define the elements in them which we consider. Additionally, we introduce a GV-type enumeration theory in the Pardon group, by mimicking the construction of \cite{KP}, and we introduce convolution logarithms of the resulting enumeration theories, which later help us to reduce to the locus of connected cycles. In \S\ref{sec:regdef}, we discuss generic transversality. In \S\ref{inc:sec:incidence}, we construct an incidence stratification that helps us identify which loci need to be studied. In \S\ref{rat:sec:rational-normal-data}, we study the smooth rational embedded locus using the techniques outlined in \S\ref{intro:sec:rational-locus-gen} above. In \S\ref{inc:sec:proof}, we study the ordinary rational--rational meeting locus. In \S\ref{sec:connected-reduction}, we put all of this together to prove Theorem~\ref{intro:thm:reduction-to-gens} (Theorem~\ref{rcl:thm:indecomposable-quotient-generators}), reducing curve enumeration to the generators. In \S\ref{liftschecks}, we collect the remaining numerical calculations, finishing the proof of the PT/GV correspondence and GV integrality.

\subsection{Acknowledgments}
We would like to  thank Kai Behrend, Jim Bryan, Dominic Joyce, Adeel Khan, Martijn Kool, Henry Liu, Davesh Maulik, Miguel Moreira, John Pardon, S\'{e}bastien Picard, and Richard Thomas for their helpful input to the present paper.

\subsubsection{AI disclosure}
The first-named author used ChatGPT 5.6 and ChatGPT 6 to find several references: most importantly, \cite{OzbagciPopescuPampu16} and \cite{Guadagni18} which played a very important role in developing the material of Section \ref{inc:sec:proof}. He also used ChatGPT 6 to edit parts of this manuscript and to proofread several proofs.

The second-named author used AI/LLMs in order to search for references, for generation of tikz graphics, and to double check and help fix errors and oversights in a preliminary version.
\section{Background}\label{background}

This section fixes the language concerning families over simplicial complexes and deformations by regluing that we use below. We follow \cite[\S\S 2.3--2.4 and 5.2]{pardon2023}. We also record the elementary facts about subdivisions which will allow us to choose constructible parameter spaces.

\subsection{Complex Geometry over Simplicial Bases}

\subsubsection{}

\begin{definition}\label{bg:def:admissible-family}
    Following \cite[Definition~2.7]{pardon2023}, a \textit{family of complex manifolds over a simplicial complex} $B$ is a separated submersion $X\to B_\bC$, where $B_\bC$ is the complexification of $B$.
\end{definition}

\begin{remark}
    By \cite[Corollary~2.6 and Lemma~2.8]{pardon2023}, the complexification $B_\bC$ is the colimit of the complexifications of the simplices of $B$. Consequently, a family over $B$ is equivalently given by compatible families over these complexified simplices. We will use this description when making constructions by induction over the simplices.
\end{remark}

\subsubsection{}

\begin{definition}
    An open subset $W$ of a finite simplicial complex $B$ is \textit{constructible} if its complement is a subcomplex. For a subdivision $B'$ of $B$ and an open subset $U\subseteq B$, we denote by $U^\Delta(B')$ the largest subcomplex of $B'$ contained in $U$. If $v$ is a vertex of $B'$, we denote its open star by $\operatorname{st}^{\circ}_{B'}(v)$. 
\end{definition}

\begin{definition}
    We say that an open subset $U$ of a finite simplicial complex $B$ is \textit{Stein} if there is a Stein open subset $\Omega\subseteq B_\bC$ such that $\Omega\cap B=U$.
\end{definition}

\subsubsection{}

\begin{lemma}\label{bg:lemma:constructible-stein}
    Every constructible open subset $U$ of a finite simplicial complex is Stein.
\end{lemma}

\begin{proof}
    As evident by the construction in \cite[\S 2.4]{pardon2023}, $B_\bC$ is Stein as it is the closed analytic subspace of the affine space with coordinates indexed by the vertices of $B$ which is cut out by the affine hyperplane equation and the monomials corresponding to non-faces. Thus, $B_\bC\times\bC$ is also Stein. Denote by $z_v$ the barycentric coordinate corresponding to a vertex $v$ on $B_\bC$, and for any non-empty simplex $S$ of $B$ consider the monomials
    \begin{equation*}
        m_S\coloneqq \prod_{v\in S} z_v.
    \end{equation*}

    Now $U$ is constructible with complement a subcomplex, which we call $L$. We can consider the subset
    \begin{equation*}
        \Omega_L\coloneqq \left\{g_L\neq 0\right\}\subseteq B_\bC,
    \end{equation*}
    where $g_L\coloneqq \sum_{S\in B, S\notin L} m_S$. We can identify $\Omega_L$ with its closed analytic subspace cut out by $g_Lw=1$ in $B_\bC\times\bC$. This shows $\Omega_L$ is Stein, and its intersection with $B$ is precisely $U$, finishing the proof.
\end{proof}

\subsubsection{}

\begin{lemma}\label{bg:lemma:star-constructible-stein}
    The open star $\st^\circ_{B'}(v)$ of a vertex $v$ with respect to a subdivision $B'\subseteq B$ is constructible and Stein.
\end{lemma}

\begin{proof}
    The complement of the open star of a vertex $v$ is the subcomplex consisting of the simplices which do not contain $v$, so the open star is constructible, and hence Stein by Lemma~\ref{bg:lemma:constructible-stein} above.
\end{proof}

\subsubsection{}

\begin{proposition}\label{bg:prop:bump-functions}
    Let $A$ be a subcomplex of a finite simplicial complex $B$. After replacing $B$ by its barycentric subdivision, there exists a continuous function $\chi_A:B\to[0,1]$, real analytic on every simplex, which vanishes on the boundary of $A$ and is positive at every point in the interior of $A$.
\end{proposition}

\begin{proof}
    A vertex of the barycentric subdivision is the barycenter $b_\sigma$ of a simplex $\sigma\subseteq B$. Set $\chi_A(b_\sigma)=1$ when $b_\sigma$ lies in the interior of $A$, and set $\chi_A(b_\sigma)=0$ otherwise. Extend $\chi_A$ affinely over every simplex of the subdivision. Every simplex meeting the interior of $A$ away from its boundary has a vertex of the first kind, while every simplex in the boundary has only vertices of the second kind. The resulting function has the required properties.
\end{proof}

\subsubsection{}

\begin{proposition}\label{bg:prop:subdivision-star-cover}
    Let $B$ be a finite simplicial complex, let $K\subseteq B$ be compact, and let $\{U_i\}_{i\in I}$ be a finite open Stein cover of $K$. There exist a subdivision $B'$ of $B$, finitely many vertices $v_a$ of $B'$, and a map $a\mapsto i(a)$ such that the open stars $\operatorname{st}^{\circ}_{B'}(v_a)$ cover $K$, each of them is contained in $U_{i(a)}$, and each of them is constructible and Stein in $B'$.
\end{proposition}

\begin{proof}
    Equip $B$ with the metric induced by a Euclidean realization. If $n=\dim B$, the diameter of every simplex in the barycentric subdivision of a simplex of diameter $r$ is at most $nr/(n+1)$; see \cite[Section~2.1, pp.~120--121]{hatcher2002algebraic}. It follows by iteration that the diameters of the simplices tend uniformly to zero. Since the diameter of a closed star is at most twice the largest diameter of a simplex which it contains, the Lebesgue number lemma gives a subdivision $B'$ such that every closed star meeting $K$ is contained in some $U_i$. Choose finitely many of the corresponding open stars which cover $K$. Open stars are constructible and Stein by Lemma~\ref{bg:lemma:star-constructible-stein}.
\end{proof}

\subsubsection{}
\begin{corollary}\label{bg:prop:subdivision-constructible-cover}
    Let $B$ be a finite simplicial complex, let $K\subseteq B$ be compact, and let $\{U_i\}_{i\in I}$ be a finite open Stein cover of $K$. There exist a subdivision $B'$ of $B$ such that there are constructible open subsets $U_i^\st\subseteq U_i$ covering $K$, with each $U_i^\st$ Stein in $B'$.
\end{corollary}

\begin{proof}
    We can apply Proposition~\ref{bg:prop:subdivision-star-cover} above to obtain a subdivision $B'$, finitely many vertices $v_a$ of $B'$, and a map $a\mapsto i(a)$ such that the open stars $\operatorname{st}^{\circ}_{B'}(v_a)$ cover $K$, and each of them is contained in $U_{i(a)}$. For $i\in I$, we then define the open subset
    \begin{equation*}
        U_i^\st \coloneqq \bigcup_{a:i(a)=i}\operatorname{st}^{\circ}_{B'}(v_a),
    \end{equation*}
    of $U_i$. It covers $K$, because the open stars do and each of them is contained in some $U_i$. Its complement is
    \begin{equation*}
        B'\setminus U_i^\st = \bigcap_{a:i(a)=i}B'\setminus\operatorname{st}^{\circ}_{B'}(v_a),
    \end{equation*}
    which is a subcomplex as an intersection of finitely many subcomplexes. Hence, $U_i^\st$ is constructible and Stein by Lemma~\ref{bg:lemma:constructible-stein}. 
\end{proof}

\subsection{Regluing Deformations}

\subsubsection{}\label{bg:sec:regluing-construction-base-version}

We first record the topological point needed to justify the construction by regluing from \cite[Definition~5.2]{pardon2023}. If $S\subseteq X$ and $\epsilon>0$, we write $\overline{N}_\epsilon(S)$ for its closed $\epsilon$ neighborhood.

\begin{proposition}\label{bg:prop:hausdorff-gluing-deformation}
    Let $(X,d)$ be a connected Hausdorff metric space, covered by open subsets $X_1$ and $X_2$, and put $U=X_1\cap X_2$. Let $\gamma:U_1\xrightarrow{\sim}U_2$ be a homeomorphism between open subsets of $U$. Choose connected open subsets $X_1^-\subseteq X_1$ and $X_2^-\subseteq X_2$, and assume that $\tilde{U}_1\coloneqq X_1^-\cap U_1\cap\gamma^{-1}(X_2^-)$ is nonempty. Suppose that, for some $\epsilon>0$,
    \begin{equation*}
        d(x,\gamma(x))<\epsilon\quad\text{for all }x\in U_1,
        \qquad
        \overline{X_1^-}\cap\overline{N}_\epsilon(\overline{X_2^-})\subseteq U_1,
    \end{equation*}
    where the closures are taken in $X$. Then the space obtained by gluing $X_1^-$ and $X_2^-$ along the restriction of $\gamma$ to $\tilde{U}_1$ is connected and Hausdorff.
\end{proposition}

\begin{proof}
    Put $\tilde{U}_2=\gamma(\tilde{U}_1)$. Since the $\tilde{U}_i$ are open, the quotient map from $X_1^-\sqcup X_2^-$ to the glued space is open. Its equivalence relation is closed if and only if the graph of $\gamma|_{\tilde{U}_1}$ is closed in $X_1^-\times X_2^-$. An open quotient $q:X_1^-\sqcup X_2^-\to (X_1^-\sqcup X_2^-)/R$ by a closed equivalence relation $R$ is Hausdorff: indeed, two points which are not equivalent have a product neighborhood $V_1\times V_2$ disjoint from $R$, and the saturations $q^{-1}(q(V_i))$ of the two factors give disjoint open neighborhoods in the quotient, since $q$ is open.

    Consider a convergent sequence $(x_j,\gamma(x_j))$ in this graph, with limit $(x,y)\in X_1^-\times X_2^-$. We have $x\in\overline{X_1^-}$, $y\in\overline{X_2^-}$, and $d(x,y)\leq\epsilon$. Hence, $x\in U_1$ by the second assumption. Continuity gives $\gamma(x_j)\to\gamma(x)$, while $\gamma(x_j)\to y$; therefore $y=\gamma(x)$. Since $x\in X_1^-$ and $y\in X_2^-$, we have $x\in \tilde{U}_1$. Thus the graph is closed and the quotient is Hausdorff. Finally, it is connected because it is the union of the images of the two connected spaces $X_1^-$ and $X_2^-$, whose intersection is the nonempty image of $\tilde{U}_1$.
\end{proof}

\subsubsection{}
\begin{proposition}\label{bg:prop:family-regluing}
    Let $\pi\colon X=X_1\cup X_2\to B$ be a family of complex manifolds over a finite simplicial complex $B$, put $U=X_1\cap X_2$, and let $V\subseteq H^0(U,T_{U/B})$ be a finite dimensional real vector space. After trimming $X$ near a prescribed compact subset, choose connected open subsets $X_i^-\subseteq X_i$ which are relatively compact over $B$, i.e.\ $\overline{X_i^-}\cap\pi^{-1}(K)$ is compact for a compact neighborhood $K$ of $B$ in $B_\bC$, such that $X_1^-\cap X_2^-\neq\varnothing$ and $\overline{X_1^-}\cap\overline{X_2^-}\subseteq U$. Let $W\subseteq H^0(U,T_{U/B})$ be the complex span of $V$. After shrinking $B_\bC$ to a neighborhood of $B$, for a sufficiently small ball $B_0W\subseteq W$, regluing by the time $1$ maps of the vector fields in $W$ gives a separated submersion $\widetilde X\to B_\bC\times B_0W$ whose restriction to $B_\bC\times\{0\}$ is $X_1^-\cup X_2^-$. After choosing a compact ball $B_0V\subseteq B_0W\cap V$ with a real analytic triangulation with $0$ as a vertex, and a triangulation of $B\times B_0V$ refining the products of simplices, its restriction to $B\times B_0V$ is a family in the sense of \cite[Definition~2.7]{pardon2023}, which restricts to $X_1^-\cup X_2^-$ over the subcomplex $B\times\{0\}$.
\end{proposition}

\begin{proof}
    Fix a metric on $X$ which induces its topology. Since the $\overline{X_i^-}\cap\pi^{-1}(K)$ are compact and $\overline{X_1^-}\cap\overline{X_2^-}$ is contained in $U$, there are open subsets $U_0\Subset U_1\Subset U$ and $\epsilon>0$ such that $\overline{X_1^-}\cap\overline{N}_\epsilon(\overline{X_2^-})\cap\pi^{-1}(K)\subseteq U_0$. We replace $B_\bC$ by the interior of $K$, and $X$ by its preimage.

    By the theorem on holomorphic ordinary differential equations, for a sufficiently small ball $B_0W$ every $w\in B_0W$ has a time $1$ map $\Phi_w$ defined on $U_1$, the map $(x,w)\mapsto\Phi_w(x)$ is holomorphic, and $d(x,\Phi_w(x))<\epsilon$ on $U_1$; see \cite[Chapter~VII]{cartan1995elementary}. Since the vector fields in $W$ are vertical, their trajectories stay in the fibers of $\pi$, so $\pi\circ\Phi_w=\pi$. Shrinking the ball once more, we may assume that each $\Phi_w$ is an isomorphism onto its image. Glue $X_1^-\times B_0W$ and $X_2^-\times B_0W$ by the isomorphism $(x,w)\mapsto(\Phi_w(x),w)$ on the open subset
    \begin{equation*}
        \tilde{\cU}_1=\{(x,w):x\in X_1^-\cap U_1,\ \Phi_w(x)\in X_2^-\}.
    \end{equation*}
    Since this isomorphism commutes with $\pi\times\mathrm{id}$, the resulting space is a complex manifold locally over $B_\bC\times B_0W$. If $(x_j,\Phi_{w_j}(x_j),w_j)$ in the graph of the gluing map converges to $(x,y,w)$, the metric bound gives $d(x,y)\leq\epsilon$, and $x\in\pi^{-1}(K)$, so $x\in U_0\subseteq U_1$ and holomorphic dependence gives $y=\Phi_w(x)$. Thus the graph is closed. As in Proposition~\ref{bg:prop:hausdorff-gluing-deformation}, the quotient map is open, and the total space is Hausdorff. The projection to $B_\bC\times B_0W$ is therefore a separated submersion. Since $\Phi_0$ is the identity and $U_1$ contains $\overline{X_1^-}\cap\overline{X_2^-}$, its restriction to $B_\bC\times\{0\}$ is $X_1^-\cup X_2^-$.

    Choose the compact real ball $B_0V$ small enough that it is contained in $B_0W$, fix a real analytic triangulation with $0$ as a vertex, and a triangulation of $B\times B_0V$ refining the products of simplices. Every simplex of $B\times B_0V$ lies in some $\sigma\times B_0V$ for a simplex $\sigma$ of $B$, and its affine inclusion into $\sigma\times V$ extends complex affinely to a map into $\sigma_\bC\times W$. Pulling the separated submersion back by these extensions gives compatible separated submersions over neighborhoods of the simplices in their complexifications. They glue over $(B\times B_0V)_\bC$ by \cite[Corollary~2.6]{pardon2023}, which proves that the restriction is a family in the stated sense.
\end{proof}

\subsection{CY Preservation}\label{bg:sec:cy-preservation}

\subsubsection{}
In our setting, it is additionally important to check that the Calabi--Yau property of the neighborhood of $C$ is preserved under the regluing deformations we consider. Due to the Bogomolov--Tian--Todorov Theorem \cite{Bogomolov1974, Tian1987, Todorov1989}, this is trivial in the compact case. However, in the non-compact case, this requires proof. First, we examine, which deformations preserve the CY property.

\begin{proposition}\label{gen:prop:cy-preservation-condition}
    Let $X$ be a (possibly non-compact) smooth complex analytic Calabi--Yau $n$-fold. Given a CY volume form $\omega\in H^0(X,\Omega_X^n)$ on $X$, the obstruction to $\omega$ extending to a relative CY volume form $\overline{\omega}\in H^0(\cX,\Omega_{\cX/A}^n)$ over a first-order deformation $\cX/A$ of $X$ is
    \begin{equation*}
        \divergence_\omega: H^1(X,T_X) \to H^1(X,\cO_X),\ v\mapsto \rd(\iota_v\omega),
    \end{equation*}
    where $v\in H^1(X,T_X)$ is the corresponding Kodaira--Spencer class, and $\iota$ is the contraction homomorphism. Write $H^1_\rCY(X,T_X)\coloneqq \ker(\divergence_\omega)$ for the CY-preserving infinitesimal deformations.
\end{proposition}

\subsubsection{}

\begin{proof}
    Let $\pi: \cX\to A=\spec(\bC[\varepsilon]/(\varepsilon^2))$, $c: X\hookrightarrow\cX$ be a first-order deformation of $X$ with Kodaira--Spencer class $v\in \Ext^1(\Omega_X,\cO_X)\cong H^1(X,T_X)$ given by the conormal exact sequence
    \begin{equation*}
        0\to\cO_X\to c^*\Omega_\cX\to\Omega_X\to 0,
    \end{equation*} 
    of the embedding $c$. Flatness of $\pi$ over the Artinian local ring $A$ forces $\cO_\cX$ to be, locally on $X$, a free $\cO_X[\varepsilon]/(\varepsilon^2)$-module of rank $1$. So, we can choose an open Stein cover $\{U_i\}$ of $X$ by coordinate polydiscs such that $\cO_\cX|_{U_i}\cong \cO_{U_i}[\varepsilon]/(\varepsilon^2)$ for each $i$. Note here that as a nilpotent thickening of $X$, $\cX$ has the same underlying topological space.
    
    On $U_{ij}=U_i\cap U_j$, the two trivializations differ by an automorphism of $U_{ij}\times A$ over $A$ restricting to the identity on $U_{ij}\times\{0\}$. This is given by 
    \begin{align*}
        \gamma_{ij}:\cO_{U_i}[\varepsilon]/(\varepsilon^2)\ &\to\ \cO_{U_i}[\varepsilon]/(\varepsilon^2),\\
        f+\varepsilon g\ &\mapsto\ f+\varepsilon\big(g+v_{ij}(f)\big),
    \end{align*}
    where the holomorphic vector fields $v_{ij}\in H^0(U_{ij},T_X)$ are \v Cech representatives of the Kodaira--Spencer class $v$.
    
    Tensoring the ideal sheaf exact sequence of $X$ in $\cX$ by $\Omega_{\cX/A}^n$, we obtain the sequence
    \begin{equation*}
        0\to c_*\Omega_X^{n} \to \Omega_{\cX/A}^n \to c_*\Omega_X^{n} \to 0
    \end{equation*}
    whose connecting homomorphism in cohomology
    \begin{equation*}
        \cdots \to H^0(\cX,\Omega_{\cX/A}^n) \to H^0(X,\Omega_X^n) \xrightarrow{\delta} H^1(X,\Omega_X^n) \to \cdots
    \end{equation*}
    determines whether the CY volume form $\omega\in H^0(X,\Omega_X^n)$ extends to a relative form $\overline{\omega}\in H^0(\cX,\Omega_{\cX/A}^n)$. By the trivializations $\cO_\cX|_{U_i}\cong \cO_{U_i}[\varepsilon]/(\varepsilon^2)$ above, we have $H^0(U_i,\Omega^n_{\cX/A})\cong \Omega^n_X(U_i)\oplus\varepsilon\,\Omega^n_X(U_i)$, which admits trivial local lifts $\alpha+\varepsilon\cdot 0\ \in H^0(U_i,\Omega_{\cX/A}^n)$ of $\alpha \in H^0(U_i,\Omega_X^n)$. Now, $v_{ij}$ is the derivative of an associated local flow $\Phi_{ij}^t$ at $t=0$, so that the pullbacks of the trivial local lifts on overlaps can be computed as
    \begin{equation*}
        \gamma_{ij}^*(\alpha+\varepsilon\cdot 0)= \alpha + \varepsilon \frac{\rd}{\rd t}\Big|_{t=0}\big(\Phi_{ij}^{t}\big)^*\alpha = \alpha +\varepsilon\cdot\cL_{v_{ij}}\alpha,
    \end{equation*} 
    where the second equality is the definition of the Lie derivative. That means the obstruction to gluing the trivial local lifts, and hence the obstruction to lifting a class, is given by the cycle represented by the \v Cech $1$-cochain $(\mathcal L_{v_{ij}}\alpha)_{ij}$, and (since $(v_{ij})$ is a cocycle and $\cL_{(-)}\alpha$ is linear) this is a $1$-cocycle representing exactly $\delta(\alpha)$
    \begin{equation*}
        \delta(\alpha) = \big[\,(\mathcal L_{v_{ij}}\alpha)_{ij}\,\big]  = \big[\,(\rd \iota_{v_{ij}}\alpha)_{ij}\,\big] \in H^1(X,\Omega_X^n),
    \end{equation*}
    where we apply Cartan's magic formula to a top-form $\alpha$ in the second equality. The latter is exactly the \v Cech-cocycle representative of $\divergence_\omega(v)$ obtained by applying $\divergence_\omega$ to the cocycle $(v_{ij})$ representing $v$. Thus,
    \begin{equation*}
        \delta(\omega) = \divergence_\omega(v)\in H^1(X,\cO_X)
    \end{equation*}
    as desired. Since $\cX$ is a nilpotent thickening of $X$, Nakayama's lemma can be used to prove that any such extension $\overline{\omega}\in H^0(\cX,\Omega_{\cX/A}^n)$ with $c^*\overline{\omega}=\omega$ is trivializing if $\omega$ is.
\end{proof}

\subsubsection{}
In order for the regluing deformation above to preserve the CY volume form, we need the slightly stronger following notion.
\begin{definition}
    Let $X/B$ be a (possibly non-compact) smooth complex analytic Calabi--Yau $n$-fold family with CY volume form $\omega\in H^0(X,\Omega_{X/B}^n)$. Define the (non-$\cO_X$-linear) morphism of sheaves  
    \begin{equation}\label{bg:eq:sheaf-divergence}
        \sdiv\coloneqq \rd(\iota_{(-)}\omega):T_{X/B}\to\cO_X
    \end{equation}
    by composing the contraction homomorphism, the differential, and the isomorphism $\Omega_{X/B}^n\cong \cO_X$. A vector field $v\in H^0(U,T_{U/B})$ on an open set $U\subseteq X$ over $B$ is called \textit{divergence-free} if $h^0(\sdiv|_U)(v)=0$. We denote the set of divergence-free vector fields by $H_{\rCY}^0(U,T_{U/B})$.
\end{definition}

\subsubsection{}
\begin{corollary}\label{bg:cor:regluing-cy-preservation}
    The regluing deformation of Proposition~\ref{bg:prop:family-regluing} yields a CY family if the finite dimensional vector space $V$ is contained in $H_{\rCY}^0(U,T_{U/B})$. 
\end{corollary}

\begin{proof}
    In the regluing of Proposition~\ref{bg:prop:family-regluing}, we only need to additionally check that pullback along the flows preserves the volume form. Take some flow $\Phi_t$ of a vector field $v$ and consider
    \begin{equation*}
        \frac{\rd}{\rd t}\Phi_t^*(\omega)=\cL_v\omega=\rd(\iota_v\omega)=0,
    \end{equation*}
    where the first equality is the definition of the Lie derivative, the second is Cartan's magic formula applied to a top-form, and the last comes from choice of $v$ in $H_{\rCY}^0(U,T_{U/B})$.
\end{proof}

\section{Universal Curve Enumeration}\label{sec:unem}

\subsection{Pardon Groups}\label{uce:sec:pardon-ring}

\subsubsection{}
\label{coh:sec:bivariant}

Following \cite[\S 2.2]{pardon2023}, we use Pardon's bivariant
groups.  For a map \(\pi:W\to B\), and with \(a:B\to\ast\), set
\begin{align*}
  H_c^*(W/B)&\coloneqq a_!\pi_!\pi^*a^!\mathbb Z,\text{ and}\\
  H_*^{\mathrm{BM}}(W/B)&\coloneqq a_*\pi_*\pi^!a^*\mathbb Z.
\end{align*}
For two such maps $\pi:W\to B$ and $\pi:W'\to B'$, given a map \(\alpha:B'\to B\), and a proper map
\(q:W\times_B B'\to W'\) over \(B'\), we obtain pushforward and pullback maps
\begin{equation}\label{uce:eq:pushforward-pullback-bivariant-theory}
  \alpha_*q^*:H_c^*(W'/B')\to H_c^*(W/B),\qquad q_*\alpha^*:H_*^{\mathrm{BM}}(W/B)\to H_*^{\mathrm{BM}}(W'/B')
\end{equation}
via the six functor formalism. In particular, a closed inclusion \(i:C\hookrightarrow W\)
over \(B\) induces restriction \(i^*:H_c^*(W/B)\to H_c^*(C/B)\).

\subsubsection{}
\label{subsec:bivariant-formal-properties}

We record the properties used in this work.  We use rational coefficients for
K\"unneth and when averaging over a finite cover.

\begin{enumerate}
\renewcommand{\theenumi}{B\arabic{enumi}}
\renewcommand{\labelenumi}{\textup{(\theenumi)}}

\item\label{biv:functoriality}
The two maps above are compatible with composition and proper base
change.

\item\label{biv:localization}
An open inclusion \(j:U\hookrightarrow W\) gives extension by zero
\(j_!:H_c^*(U/B)\to H_c^*(W/B)\).  If
\(i:C\hookrightarrow W\) is its closed complement, the triangle
\(j_!\mathbb Q_U\to\mathbb Q_W\to i_*\mathbb Q_C\xrightarrow{+1}\)
induces the localization sequence
\(\cdots\to H_c^q(U/B)\xrightarrow{j_!}H_c^q(W/B)
\xrightarrow{i^*}H_c^q(C/B)\xrightarrow{\partial}
H_c^{q+1}(U/B)\to\cdots\). 

\item\label{biv:kunneth}
There are associative K\"unneth maps, commutative in the graded sense,
from \(H_c^*(W/B)\otimes H_c^*(W'/B')\) to
\(H_c^*(W\times W'/B\times B')\).  They are isomorphisms in every
situation used below.

\item\label{biv:compatibilities}
These constructions commute with restriction to faces.  Moreover,
if \(p:\widetilde W\to W\) is a finite cover of degree \(d\), then there is a trace map $\Tr_p$ such that 
\(\operatorname{Tr}_p\circ p^*=d\,\mathrm{id}\).

\item\label{biv:diagonal-gysin}
Let \(Y\to B\) be relatively smooth and let \(W\to Y\times_BY\) be a
space over \(B\).  The complex orientation of the relative diagonal
gives a refined pullback from \(H_k^{\mathrm{BM}}(W/B)\) to
the Borel--Moore group of \(W\times_{Y\times_BY}Y\) in degree
\(k-2\dim_{\mathbb C}(Y/B)\).  It commutes with proper pushforward,
Cartesian base change, exterior products, and restriction to faces.

\item\label{biv:cap-product-pairing}
We have a cap product $\cap:H^q(W)\otimes H_p^\BM(W/B)\to H_{p-q}^\BM(W/B)$ and a pairing
\begin{equation*}
  \langle-,-\rangle:H_p^\BM(W/B)\otimes H_c^p(W/B)\to \bZ
\end{equation*}
for classes of matching degree. This satisfies a push-pull formula with respect to the maps of \eqref{uce:eq:pushforward-pullback-bivariant-theory}.

\end{enumerate}

\subsubsection{}

\begin{definition}
\label{coh:def:grothendieck-group-1-cycles}
The \textit{Pardon group of curve enumeration problems} is
\begin{equation*}
  \cH\coloneqq H_c^*(\mathcal Z/\CY_4)\coloneqq \varinjlim_{X\to B}H_c^*(\mathcal Z(X/B)/B).
\end{equation*}  
The colimit is taken
under pushforward of \eqref{uce:eq:pushforward-pullback-bivariant-theory} over all families of not necessarily proper
Calabi--Yau fourfolds.  A morphism is a Cartesian square up to an open
embedding: the map of bases \(B'\to B\) is injective and the induced
map \(X'\to X\times_BB'\) is an open embedding.  As in
\cite[Lemma~3.2]{pardon2023}, the indexing category may be replaced by
a directed one. Since this is the group we mainly study throughout the paper, we sometimes use the shorthand $\cH$ introduced above.

The \textit{Pardon group of curve enumeration theories} is the inverse limit 
\begin{equation*}
  H_*^{\mathrm{BM}}(\mathcal Z/\CY_4)\coloneqq \varprojlim_{X\to B}H_*^{\mathrm{BM}}(\mathcal Z(X/B)/B)
\end{equation*}
over all families of Calabi--Yau fourfolds, under the pullback maps of \eqref{uce:eq:pushforward-pullback-bivariant-theory}. These two groups come with a
natural pairing
\begin{equation*}
  \langle-,-\rangle:
 H_*^{\mathrm{BM}}(\mathcal Z/\CY_4)\otimes
H_c^*(\mathcal Z/\CY_4)\to\mathbb Q.
\end{equation*}
\end{definition}

\subsubsection{}
\label{subsec:pardon-product}
Suppose
that \(\lambda\) and \(\lambda'\) are represented over \(X\to B\) and
\(X'\to B'\).  Over \(B\times B'\), form the disjoint union
\(X\boxplus X'\coloneqq (X\times B')\sqcup(B\times X')\).  Its cycle space is
canonically
\(\mathcal Z(X/B)\times\mathcal Z(X'/B')\), and the K\"unneth exterior
product defines $\label{eq:pardon-product-representative}
  \lambda\lambda'\in\mathcal H$.
  
Functoriality shows that this depends only on the two universal
classes.  Associativity of disjoint union and symmetry of K\"unneth make
\(\mathcal H\) an algebra which is commutative in the graded sense.  Its
unit is represented by the empty cycle in the empty target over a
point.

Addition of effective cycles gives the coalgebra
operation.  For every \(X\to B\), the map
\(\Sigma:\mathcal Z(X/B)\times_B\mathcal Z(X/B)\to\mathcal Z(X/B)\)
is proper.  Pullback by \(\Sigma\), change from the diagonal base
\(B\) to \(B\times B\), and K\"unneth define
\begin{equation}
\label{eq:pardon-coproduct-representative}
  \Delta:\mathcal H\to\mathcal H\otimes\mathcal H.
\end{equation}
This coproduct is coassociative and cocommutative in the graded sense,
and is compatible with the product.  Together with the unit and the
counit below, it makes \(\mathcal H\) a bialgebra.

\subsubsection{}
\label{subsec:augmentation-ideal}
Restriction to the section $B\to \cZ(X/B)$ of empty cycles is the counit
\begin{equation}\label{uce:eq:counit}
  \epsilon:\mathcal H\to\mathbb Q
\end{equation}
and its kernel is the
\textit{augmentation ideal} \(\mathfrak m\).
\label{eq:augmentation}
The colimit of bivariant groups $H_c^*(B/B)$ of this section of empty cycles is simply \(\mathbb Q\cdot1\), since, after an
embedding of the base into a simplex, its reduced homology vanishes in the
universal colimit. Since the empty cycle is both open and closed, and
writing \(\mathcal Z^+\coloneqq \mathcal Z\setminus\{0\}\), we obtain
\begin{equation}
 \mathcal H=\mathbb Q\cdot1\oplus
 H_c^*(\mathcal Z^+/\CY_4;\mathbb Q),
 \qquad
 \mathfrak m=H_c^*(\mathcal Z^+/\CY_4;\mathbb Q).
 \label{eq:positive-part}
\end{equation}
Thus \(\mathfrak m\) consists precisely of classes supported away from
the empty cycle, and \(\mathfrak m^r\) is generated by products of
\(r\) such classes.  The relation between \(\mathfrak m^2\) and
\textit{cycles with disconnected support} is established in Proposition~\ref{prop:disconnected-square}.

\subsubsection{}\label{uce:sec:convolution-product}
Now we want to define one further operation on sequences of linear forms on $\cH$, which is used later in \S\ref{subsec:convolution-logarithms} to define logarithms of GW and PT theories.

For two sequences $\mathbf f=(f_n)_{n\in\mathbb Z}$ and
$\mathbf h=(h_n)_{n\in\mathbb Z}$ in $\Hom_\bQ(\cH,\bQ)$, define their \textit{convolution} by
\begin{equation}\label{uce:eq:convolution-product}
  (\mathbf f\star\mathbf h)_n(x)
  \coloneqq \sum_{a+b=n}\mu_{\mathbb Q}\bigl((f_a\otimes h_b)\Delta x\bigr).
\end{equation}

\subsubsection{Multiplicity Filtration}
\label{subsec:multiplicity-filtration}

For later use, we recall Pardon's multiplicity filtration
\cite[Definition~2.1]{pardon2023}.  Let
\(\mathrm M\coloneqq \coprod_{r\geq0}\mathbb Z_{\geq1}^{r}/S_r\) be the set of
finite multisets of positive integers.  We declare that
\(\mathbf m=[m_1,\ldots,m_r]\geq
\mathbf m'=[m'_1,\ldots,m'_s]\) if there are integers
\(n_{ij}\geq0\) such that
\begin{equation}
 m_i=\sum_{j=1}^{s}n_{ij}m'_j\quad\text{for every }i,
 \qquad
 \sum_{i=1}^{r}n_{ij}>0\quad\text{for every }j.
 \label{eq:multiplicity-order}
\end{equation}
This order records the possibility that a component of multiplicity
\(m_i\) breaks into several nearby components, and the second condition
says that every nearby component occurs in at least one of these
splittings.

For a relative cycle \(z=\sum_i m_i[C_i]\), put
\(\mathbf m(z)\coloneqq [m_1,\ldots,m_r]\), where the \(C_i\) are distinct
and irreducible.  The resulting map is
upper semicontinuous: if \(z'\) is sufficiently close to \(z\), then
\(\mathbf m(z')\leq\mathbf m(z)\).  For
\(\mathbf m\in\mathrm M\), let \(\mathcal Z_{\leq\mathbf m}\),
\(\mathcal Z_{<\mathbf m}\), and \(\mathcal Z_{=\mathbf m}\) denote the
corresponding loci.  The first two are open, while the last is closed in
\(\mathcal Z_{\leq\mathbf m}\).  The formation of these loci commutes
with  base change, restriction to faces, and restriction to an
open subfamily of the target.  Applying \eqref{biv:localization} and passing to the
universal colimit gives the localization sequence
\begin{equation}
  \cdots\to H_c^q(\mathcal Z_{<\mathbf m}/\CY_4)\to H_c^q(\mathcal Z_{\leq\mathbf m}/\CY_4)\to H_c^q(\mathcal Z_{=\mathbf m}/\CY_4)\xrightarrow{\,\partial_{\mathbf m}\,}H_c^{q+1}(\mathcal Z_{<\mathbf m}/\CY_4)\to\cdots .
 \label{eq:multiplicity-localization}
\end{equation}
We refer to these sequences, and to their versions for arbitrary
downward closed subsets of \(\mathrm M\), as the multiplicity
filtration.

A locus with fixed multiplicity is in general only locally closed.  A class
on it therefore has no canonical extension by zero to the full cycle
space.  This is why the classes introduced below on such loci are
supplied with ambient lifts, which may contain corrections on strata of
smaller multiplicity.

\subsubsection{}\label{uce:sec:total-support-multiplicity}
We shall also use the \textit{total support multiplicity}
\begin{equation*}
  \lambda\big(\textstyle\sum_i m_i[C_i]\big)\coloneqq \textstyle\sum_i m_i.
\end{equation*}
Write
\(\mathcal Z_{\leq N}\) for the open locus on which
\(\lambda\leq N\), and \(F_{\leq N}\mathcal H\) for the image of its
compactly supported cohomology in \(\mathcal H\).  These open loci exhaust
the cycle space, every compact subset is contained in one of them, and
both disjoint union and addition of cycles add total support
multiplicity.  In particular, an \(r\)-fold convolution product of maps
that vanish on the summand $\bQ\cdot 1$ of the section of empty cycles, vanishes on
\(F_{\leq N}\mathcal H\) for \(r>N\). 

Now note that, since the open loci $\cZ_{\leq N}(X/B)$ exhaust the cycle space for every family $X\to B$, we have
\begin{equation}\label{uce:eq:H-total-mult-exhaustion}
  \mathcal H=\varinjlim_N H_c^*(\mathcal Z_{\leq N}/\CY_4;\mathbb Q)
\end{equation}
and define the completion
\begin{equation}
 \hat{H}_{*}^{\mathrm{BM}}(\mathcal Z/\CY_4) \coloneqq \varprojlim_N H_*^{\mathrm{BM}}(\mathcal Z_{\leq N}/\CY_4).
 \label{eq:multiplicity-completions}
\end{equation}

\subsubsection{Restrictions of the Universal Bivariant Group}
In order for collections of locally closed loci $\cZ_\alpha$ of $\cZ$ to have well-defined universal bivariant groups using the construction of \S\ref{uce:sec:pardon-ring}, certain conditions need to be satisfied, see \cite[Def. 3.13]{pardon2023}. For future use, we record here some cases in which these conditions are naturally satisfied.

\begin{definition}\label{uce:def:fiberwise-local}
  Let $\cZ_\alpha\subseteq \cZ$ be the specification of a locally closed subset $\cZ_\alpha(X/B)\subseteq \cZ(X/B)$ for every family $X\to B$ of CY$4$-folds. We say $\cZ_\alpha$ is a \textit{fiberwise sublocus} if for any family $X\to B$ we have $\cZ(X/B)=\bigcup_b \cZ_\alpha(X_b)$. We say it is \textit{fiberwise local} if additionally $\cZ(X^-)\cap\cZ_\alpha(X)=\cZ_\alpha(X^-)$ for any open subset $X^-\subseteq X$.
\end{definition}

\begin{example}
  For example this includes the locus of connected cycles $\cZ_\con$, multiplicity and normalized Euler characteristic conditions, the loci of smooth rational curves $R_m^\sm$ and many more.
\end{example}

\subsubsection{}

\begin{proposition}
\label{rcl:prop:cohomology-fiberwise-cond}
Any locally closed fiberwise local sublocus $\cZ_\alpha\subseteq \cZ$ is compatible with pullbacks and local in the sense of \cite[Def. 3.13]{pardon2023}. Consequently, we can define
the groups
\[
 H_c^*(\cZ_\alpha/\CY_4;\mathbb Q)
 \coloneqq \varinjlim_{X\to B}H_c^*(\cZ_\alpha(X/B)/B;\mathbb Q)
\]
for any such $\cZ_\alpha$ by the construction of \S\ref{uce:sec:pardon-ring}.
\end{proposition}

\begin{proof}
The fiberwise condition makes $\cZ_\alpha$ immediately compatible with base change. Now take an open subset $X^-\subseteq X$. By the fiberwise condition and fiberwise locality, we have
\begin{equation*}
  \cZ_\alpha(X^-/B)=\bigcup_b\cZ(X^-_b)= \bigcup_b \big(\cZ(X^-_b)\cap\cZ_\alpha(X_b)\big)=\cZ(X^-/B)\cap \cZ_\alpha(X/B),
\end{equation*}
so $\cZ_\alpha$ is local.
\end{proof}

\subsubsection{}\label{subsec:support-evaluation-space}
For a family \(X\to B\), consider the cycle space with markings
\begin{equation}\label{uce:eq:cycle-space-marking-maps}
  \begin{tikzcd}[row sep=small,column sep=small]
  & \cZ^\bullet(X/B)\ar[dl,"p"']\ar[dr,"\ev"] & \\
  \cZ^+(X/B) & & X,
  \end{tikzcd}
\end{equation}
which is the incidence space of pairs \((z,x)\),
where \(z\in\mathcal Z^+(X/B)\) and \(x\in|z|\), and $p$ and $\ev$ denote the natural projection maps. Integrations along \(p\) are always taken with the multiplicities of the
universal relative cycle: over \(z=\sum_i m_i[C_i]\), the fiber is the
cycle \(\sum_i m_i[C_i]\), and not its reduction.

Note that the map \(p\) is proper, and the construction commutes with all base
changes of the family. 

We use the notation
\(H_*^{\mathrm{BM}}(\mathcal Z^\bullet/\CY_4)\) for the
inverse system obtained from these spaces. We also use the notation
\( \hat{H}_{*}^{\mathrm{BM}}
(\mathcal Z^\bullet/\CY_4)\)
for the corresponding inverse limit in support multiplicity.

For \(d\geq1\), multiplication of cycles gives proper maps
\([d](z)=dz\) and \([d]^\bullet(z,x)=(dz,x)\).  We denote their
pushforwards on bivariant Borel--Moore theories by 
\begin{equation}\label{uce:eq:mult-BM-pushforwards}
  V_d:\hat{H}_{*}^{\mathrm{BM}}(\mathcal Z/\CY_4)\to \hat{H}^\BM_*(\cZ/\CY_4),\quad\text{and}\quad V_d^\bullet:\hat{H}_{*}^{\mathrm{BM}}(\mathcal Z^\bullet/\CY_4)\to \hat{H}^\BM_*(\cZ^\bullet/\CY_4),
\end{equation}
respectively.  They satisfy \(V_dV_e=V_{de}\), and likewise in the
marked incidence space.  On every part of bounded total support
multiplicity, only finitely many of the \(V_d\) can contribute.

\subsection{GW and PT Enumeration Theories}\label{subsec:normalized-family-theories}
First, we discuss the enumeration theories we work with here. These are given by pushforwards of virtual fundamental classes along the cycle support maps.
\subsubsection{Support Map}
Before discussing virtual classes, we discuss the geometric setup. 

For the moduli of stable maps and stable sheaves\footnote{See \cite[Remark 1.8]{pardon2023} for existence of these spaces in the complex analytic setting.}, we have \textit{support maps}
\begin{equation}\label{vfc:properness}
  \begin{tikzcd}
    \bar{\cM}(X/B)\ar[dr,"\rho_{\GW}"'] & & \cP(X/B)\ar[dl,"\rho_{\PT}"] & \bar{\cM}^\con_{0,1}(X/B)\ar[d,"\rho_{0,1}"]\\
    & \cZ(X/B) & & \cZ^\bullet(X/B)
  \end{tikzcd}
\end{equation}
for every family $X\to B$. These maps record the support (and the image of
the marking), and they are proper on every truncation.  Moreover, for every compact subset \(K\)
of the relative cycle space there is an integer \(n_K\) such that the
inverse image of \(K\) in \(P_n(X/B)\) is empty for
\(n<n_K\).

\subsubsection{GW Enumeration Theory}
As explained in \cite[\S3, App. A]{pardon2023}, the moduli spaces of stable maps $\bar{\cM}(X/B)$ and $\bar{\cM}^\con(X/B)$ have virtual fundamental classes
\begin{equation}\label{vfc:families}
  \big[\bar{\cM}(X/B)\big]^\vir,\qquad \big[\bar{\cM}^\con(X/B)\big]^\vir,
\end{equation}
which form compatible inverse systems under pullback over the category of  families. For a connected stable map write
\(\chi=2\chi(C,\mathcal O_C)\).  Its genus zero and genus one values are
\(2\) and \(0\). Pushing forward the virtual classes of \eqref{vfc:families} above,
and also the class in genus zero with one marking, along the proper support maps \eqref{vfc:properness} above, gives compatible
systems
\begin{equation}
\begin{aligned}
 \mathbf{GW}_0^{\mathrm{raw}}&\coloneqq (\rho_{\mathrm{GW},2}^{\mathrm{con}})_* [\overline{\mathcal M}_{0}^{\mathrm{con}}]^{\mathrm{vir}}\in \hat{H}^\BM_2(\cZ/\CY_4),\\
 \mathbf{GW}_{0,1}^{\mathrm{raw}}&\coloneqq (\rho_{0,1})_*[\overline{\mathcal M}_{0,1}^{\mathrm{con}}]^{\mathrm{vir}}\in \hat{H}^\BM_4(\cZ^\bullet/\CY_4),\\
 \mathbf{GW}_1^{\mathrm{raw}}&\coloneqq (\rho_{\mathrm{GW},0}^{\mathrm{con}})_*[\overline{\mathcal M}_{1}^{\mathrm{con}}]^{\mathrm{vir}}\in \hat{H}^\BM_0(\cZ/\CY_4),
\end{aligned}
 \label{eq:connected-gw-inputs}
\end{equation}
where the Borel--Moore degrees \(2,4,0\) are computed by the virtual dimension $2(1-g+n)$. When we want to stress the Euler characteristic dependence, rather than the genus dependence, we occasionally write
\(\mathbf{GW}_0^{\mathrm{raw}}=\mathbf{gw}_{\chi=2}^{\mathrm{raw}}\) and
\(\mathbf{GW}_1^{\mathrm{raw}}=\mathbf{gw}_{\chi=0}^{\mathrm{raw}}\).

\subsubsection{PT Enumeration Theory}
For the applications of our study of the Pardon group of enumeration problems, we will use the following result
\begin{proposition}\label{prop:family-cy4-virtual-package}
For every analytic $CY4$-family \(X\to B\) over a simplicial complex $B$ as in Definition~\ref{bg:def:admissible-family}, the stable
pair spaces \(P_n(X/B)\) have
canonical relative PT virtual classes
\begin{equation*}
  \big[P_n(X/B)\big]^\vir \in H^\BM_{2n}(P_n(X/B)/B;\bZ),
\end{equation*}
so that these classes form
compatible inverse systems under pullback over the category of families.
\end{proposition}

\begin{proof}
This follows from the two results by the first author. Indeed, the strategy to prove this statement is to first find a \textit{canonical orientation} for our PT-families, and then to apply an analytic version of the result in \cite{Borisov_2017}, which yields the desired virtual class. The first step was done in \cite{karpov-cy4-orientations}, and the second is the subject of the first author's work \cite{karpov-bj}.
\end{proof}

\subsection{Enumeration Problems}
\subsubsection{}\label{subsubsec:primary-insertions}
We recall now, following \cite[\S 3.5]{pardon2023}, the notion of \textit{primary insertions}, and how primary insertions
give curve enumeration problems.

Let \(X\) be a smooth projective
Calabi--Yau fourfold, let \(\beta>0\) be effective, and denote by
\(\mathbf1_\beta\) the characteristic class of the compact, open and
closed locus \(\mathcal Z(X,\beta)\).  For
\(\gamma\in H^m(X;\mathbb Q)\), set
\begin{equation*}
  \widehat\gamma\coloneqq p_!\operatorname{ev}^*\gamma\in H^{m-2}(\mathcal Z^+(X);\mathbb Q)
\end{equation*}
where integration along \(p\)
uses the universal cycle with the multiplicities specified in
\eqref{uce:eq:cycle-space-marking-maps}. Additivity of this
cycle, and its behaviour under multiplication by \(d\), give
\begin{equation}
 \Sigma^*\widehat\gamma
 =\widehat\gamma\boxtimes1+1\boxtimes\widehat\gamma,
 \qquad [d]^*\widehat\gamma=d\widehat\gamma.
 \label{eq:primary-insertion-behaviour}
\end{equation}
For classes \(\gamma_i\in H^{m_i}(X;\mathbb Q)\), Pardon's curve
enumeration problem with primary insertions $(X,\beta;\gamma_\bullet)$ is the image
\begin{equation}
 \mathbf1_\beta\prod_{i=1}^r\widehat\gamma_i
 \in H_c^{q}(\mathcal Z(X);\mathbb Q)
 \longmapsto (X,\beta;\gamma_\bullet)\in H_c^q(\cZ/\CY_4;\bQ),
 \label{eq:primary-insertion-problem}
\end{equation}
where we write $q\coloneqq \sum_i(m_i-2)$. If no insertion is given, we simply write $(X,\beta)$ for the image of the class $\mathbf1_{\beta}$.

Unlike these \textit{primary insertions}, descendent insertions do not, in
general, come from classes on \(\mathcal Z(X)\), and are not considered
in the present paper.

\subsubsection{}
For future reference, we record some basic observations about such primary insertions, which follow from the constraint $\sum_i(m_i-2)=2$. We will use this to reduce to insertions with $m_i\geq 3$.
\begin{lemma}\label{uce:lemma:divisor-insertions}
  If $\gamma\in H^{<2}(X;\bQ)$, then $\widehat{\gamma}$ vanishes. Take insertions $(\gamma_\bullet,\delta)$ with $\delta\in H^2(X;\bQ)$, then
  \begin{equation*}
    (X,\beta;(\gamma_\bullet,\delta))=(\textstyle\int_\beta \delta)\cdot(X,\beta;\gamma_\bullet).
  \end{equation*}
\end{lemma}
\begin{proof}
  If $\gamma\in H^{<2}(X;\bQ)$, then $\widehat{\gamma}$ vanishes as it lies in negative degree. Take insertions $(\gamma_\bullet,\delta)$ with $\delta\in H^2(X;\bQ)$. Then $\widehat{\delta}\in H^0(\cZ^+(X);\bQ)$ is locally constant with value $\int_\beta \delta$ at a point of class $\beta$. This yields $(X,\beta;(\gamma_\bullet,\delta))=(\int_\beta \delta)(X,\beta;\gamma_\bullet)$ by the definition in \eqref{eq:primary-insertion-problem}.
\end{proof}

\subsubsection{}
We first describe its evaluation on the connected genus zero GW class.

Let \(\overline{\mathcal M}_{0,r}^{\mathrm{con}}(X,\beta)\) be the
moduli stack of stable maps
\([f:(C,x_1,\ldots,x_r)\to X]\) from connected curves of genus zero
and class \(\beta\).  The support map
\(\rho_{\mathrm{GW}}:
\overline{\mathcal M}_{0,0}^{\mathrm{con}}(X,\beta)
\to\mathcal Z(X,\beta)\) sends \([f:C\to X]\) to \(f_*[C]\). From \cite[\S3.5]{pardon2023}, by construction of the class $\widehat{\gamma}$, we obtain
\begin{equation}
 \left\langle
 (\rho_{\mathrm{GW}})_*
 [\overline{\mathcal M}_{0,0}^{\mathrm{con}}(X,\beta)]^{\mathrm{vir}},
 \mathbf1_\beta\prod_{i=1}^r\widehat\gamma_i
 \right\rangle = \int_{[\overline{\mathcal M}_{0,0}^{\mathrm{con}}
 (X,\beta)]^{\mathrm{vir}}}
 \prod_{i=1}^r\pi_!\operatorname{ev}^*\gamma_i =
 \int_{[\overline{\mathcal M}_{0,r}^{\mathrm{con}}
 (X,\beta)]^{\mathrm{vir}}}
 \prod_{i=1}^r\operatorname{ev}_i^*\gamma_i,
 \label{eq:primary-insertion-gw-evaluation}
\end{equation}
where the second equality uses an induction argument involving the forgetting tails axiom \cite[Theorem~6, Axiom~IV]{Behrend1997GW}. The real virtual dimension of
\(\overline{\mathcal M}_{0,r}^{\mathrm{con}}(X,\beta)\) is \(2+2r\).
Thus, when \(\sum_i(m_i-2)=2\), the last integral is the numerical
genus zero invariant \(\mathrm{GW}_{0,\beta}(\gamma_\bullet)\) of \eqref{intro:eq:raw-primary-gw}, showing that $\mathbf{GW}_0^{\mathrm{raw}}$ and $(X,\beta;\gamma_\bullet)$ pair to give the invariant $\mathrm{GW}_{0,\beta}(\gamma_\bullet)$. Similarly, $\mathbf{GW}_1^{\mathrm{raw}}$ and $(X,\beta)$ pair to give the invariant $\mathrm{GW}_{1,\beta}$.

\subsubsection{}
We next describe the PT evaluation.  Let
\(P_n(X,\beta)\) be the moduli space of stable pairs
\(s:\mathcal O_X\to \cF\), as in \eqref{intro:eq:stable-pairs-moduli}. The support map of \eqref{vfc:properness} is
\begin{equation*}
  \rho_{\mathrm{PT},n,\beta}:P_n(X,\beta)\to\mathcal Z(X,\beta),\ [s:\mathcal O_X\to \cF]\mapsto[\cF].
\end{equation*}
Let \(\mathbb F\) be the universal sheaf on \(P_n(X,\beta)\times X\),
and let \(\pi_P\) and \(\pi_X\) be the projections to \(P_n(X,\beta)\)
and \(X\), respectively.  The class \(\operatorname{ch}_3(\mathbb F)\)
represents the family of fundamental cycles of the sheaves parametrized
by \(P_n(X,\beta)\).  Consequently,
\begin{equation}
 \rho_{\mathrm{PT},n,\beta}^*\widehat\gamma
 = (\pi_P)_*\bigl(\operatorname{ch}_3(\mathbb F)
       \cup\pi_X^*\gamma\bigr)
 =:\tau(\gamma).
 \label{eq:primary-insertion-pt-pullback}
\end{equation}
The coefficient of \(q^n\) in the PT theory corresponds to
\((\rho_{\mathrm{PT},n,\beta})_*
[P_n(X,\beta)]^{\mathrm{vir}}\).  This virtual class has real degree
\(2n\).  Hence, when
\(\sum_i(m_i-2)=2n\), its evaluation on
\((X,\beta;\gamma_\bullet)\) is the numerical invariant
\begin{equation}
 \left\langle
 (\rho_{\mathrm{PT},n,\beta})_*
 [P_n(X,\beta)]^{\mathrm{vir}},
 \mathbf1_\beta\prod_{i=1}^r\widehat\gamma_i
 \right\rangle
 =\int_{[P_n(X,\beta)]^{\mathrm{vir}}}
       \prod_{i=1}^r\tau(\gamma_i)
 \eqqcolon P_{n,\beta}(\gamma_\bullet).
 \label{eq:primary-insertion-pt-evaluation}
\end{equation}
For $n=0,1$, this is the same expression as \eqref{intro:eq:primary-pt} when the same primary insertions are chosen.

\subsection{GV Type Enumeration Theory}
\subsubsection{}\label{uce:sec:local-laurent-condition}
Throughout this section,
every sequence $A=(A_n)_{n\in\mathbb Z}$ of classes in the
Borel--Moore theory is understood to satisfy the following \textit{local Laurent}
condition.

\begin{definition}
  Let $A=(A_n)_{n\in\mathbb Z}$ be a sequence of classes in $\hat{H}^\BM_*(\cZ/\CY_4)$. We say $A$ is \textit{locally Laurent} if, for every $x\in\mathcal H$, there is an integer $n_x$ such that
$\langle A_n,x\rangle=0$ whenever $n<n_x$. Equivalently,
$\sum_{n\in\mathbb Z}\langle A_n,x\rangle q^n$ is a Laurent series for every
$x\in\mathcal H$.
\end{definition}

\subsubsection{}\label{subsec:connected-gw-inputs}

Recall the GW curve enumeration theories from \eqref{eq:connected-gw-inputs}. Following Klemm--Pandharipande's definition of Gopakumar--Vafa type invariants $n_{g,\beta}$ from $\GW_{g,\beta}$ in \eqref{eq:gv-kp-definitions}, we now want to define similar transformed curve enumeration theories $\bfN_0$ and $\bfN_1$ that evaluate to $n_{0,\beta}(\gamma_\bullet)$ and $n_{1,\beta}$ when paired with the classes $(X,\beta;\gamma_\bullet)$ and $(X,\beta)$ from \eqref{eq:primary-insertion-problem}, respectively.

\begin{remark}
The resulting \textit{transformed classes} defined for families $X\to B$ (which are defined below) will form compatible systems, and we denote the corresponding classes in $\hat{H}^\BM_*(-/\CY_4)$ also simply by the same symbols.
\end{remark}

In \eqref{eq:gv-kp-definitions}, consider the first equation and the second term on the RHS of the second equation. Following these, using the multiplication pushforward $V_d=[d]_*$ from \eqref{uce:eq:mult-BM-pushforwards}, define the two genus zero multiple cover operators and the corresponding primitive
classes by
\begin{equation}
  \begin{aligned}
    \mathsf A_3&\coloneqq \sum_{d\geq1}\frac{1}{d^3}V_d,\qquad &\mathsf A_2&\coloneqq\sum_{d\geq1}\frac{1}{d^2}V_d^\bullet,\\
    \mathbf N_0&\coloneqq \mathsf A_3^{-1}\mathbf{GW}_0^{\mathrm{raw}}\in \hat{H}^\BM_2(\cZ/\CY_4),\qquad &\mathbf N_{0,1}&\coloneqq\mathsf A_2^{-1}\mathbf G_{0,1}^{\mathrm{raw}}\in \hat{H}^\BM_4(\cZ^\bullet/\CY_4).
  \end{aligned}
 \label{eq:primitive-genus-zero-classes}
\end{equation}
The inverses are Dirichlet inverses, i.e.
\(\mathsf A_3^{-1}=\sum_{d\geq1}\mu(d)d^{-3}V_d\), and \(\mathsf A_2^{-1}=\sum_{d\geq1}\mu(d)d^{-2}V_d^\bullet\).
These transformed classes $\bfN_0$ and $\bfN_{0,1}$ are locally Laurent in the sense of Definition~\ref{uce:sec:local-laurent-condition} above.

\begin{lemma}
  For insertions $\gamma_i\in H^{m_i}(X;\bQ)$ with $\sum_i (m_i-2)=2$, we have
  \begin{equation*}
    \langle\bfN_0,(X,\beta;\gamma_\bullet)\rangle=n_{0,\beta}(\gamma_\bullet).
  \end{equation*}
\end{lemma}
\begin{proof}
  We apply the defining identity \(\mathbf G_0^{\mathrm{raw}}=A_3\mathbf N_0\) to $\sum_{\beta>0}(X,\beta;\gamma_\bullet)Q^\beta$. Using the push-pull formula on $V_d=[d]_*$, this yields
\begin{equation*}
  \sum_{d\geq 1}d^{-3}\sum_{\beta>0}\langle\bfN_0,[d]^*(X,\beta;\gamma_\bullet)\rangle Q^\beta.
\end{equation*}
Consider the term in the inner sum. If $\beta$ is not divisible by $d$, it vanishes, but if $\beta= d\beta'$ for some $\beta'>0$, then each primary
insertion contributes a
factor \(d\) under \([d]^*\) by linearity in the support cycle. After pairing with the \(r\) insertions
\(\gamma_\bullet\), the resulting factor is \(d^{r}\), and the summand becomes $d^r\langle\bfN_0,(X,\beta';\gamma_\bullet)\rangle Q^{d\beta'}$. Relabeling, we obtain exactly the RHS of the first equation in \eqref{eq:gv-kp-definitions}, so that we can read off the claim.
\end{proof}

\subsubsection{}
\label{subsec:c2-and-incidence-classes}
Using the multiple cover operators above, we have obtained the transformed genus zero contribution $\bfN_0$. Now, we first want to deal with the rational contributions in genus one, starting with $n_{0,\beta}(c_2(TX))$. Using the maps in \eqref{uce:eq:cycle-space-marking-maps}, define
\begin{equation*}
  \mathbf C_2\coloneqq p_*(\operatorname{ev}^*c_2(T_{X})\cap\mathbf N_{0,1})\in H^\BM_0(\cZ^+(X/B)/B).
\end{equation*}

\begin{lemma}\label{uce:lemma:c_2-class-evaluation}
  We have $\langle \mathbf C_2,(X,\beta)\rangle=n_{0,\beta}(c_2(TX))$.
\end{lemma}
\begin{proof}
  We first note that $\langle \mathbf C_2,(X,\beta)\rangle=\langle\bfN_{0,1},\mathbf1_\beta\ev^*c_2(TX)\rangle$. As above, using the defining identity
\(\mathbf G_{0,1}^{\mathrm{raw}}=A_2\mathbf N_{0,1}\), the claim then reduces to calculating $([d]^\bullet)^*\mathbf1_{d\beta'}\ev^*c_2(TX)$. Since the insertion $\ev^*c_2(TX)$ is independent of the curve class, we simply get $([d]^\bullet)^*\mathbf1_{d\beta'}\ev^*c_2(TX)=\mathbf1_{\beta'}\ev^*c_2(TX)$, from which the claim follows using the first equation in \eqref{eq:gv-kp-definitions}.
\end{proof}

\subsubsection{}
Now we consider the term involving meeting numbers $m_{\beta_1,\beta_2}$ in \eqref{eq:gv-kp-definitions}, and, motivated by the recursive meeting rules of \cite[Section~1.2]{KP}, define classes $\bfM$ that evaluate to the meeting numbers. 

Form the fiber product
\(\mathcal Z^\bullet(X/B)\times_X\mathcal Z^\bullet(X/B)\) using the two evaluation maps, and let
\(p^{(2)}\) be its proper map to
\(\mathcal Z^+(X/B)\times_B\mathcal Z^+(X/B)\). 

Define the \textit{ordered incidence class} as
\begin{equation*}
  \mathbf I\coloneqq (p^{(2)})_*\Delta_{X/B}^{!}(\mathbf N_{0,1}\boxtimes_B\mathbf N_{0,1})\in H^\BM_0(\cZ^+(X/B)\times_B\cZ^+(X/B)/B).
\end{equation*}
(Note that interchanging its two factors preserves the class.) We can now use this incidence class to define a meeting class as follows. On \(\mathcal Z^+(X/B)\times_B\mathcal Z^+(X/B)\), consider the proper maps
\(s_L(a,b)=(a,a+b)\), \(s_R(a,b)=(a+b,b)\), and
\(\Sigma(a,b)=a+b\), together with the diagonal
\(\Delta(a)=(a,a)\).  The \emph{renormalized meeting class} is the
unique class \(\mathbf M\in \hat{H}^\BM_0(\cZ^+(X/B)\times_B\cZ^+(X/B)/B)\) of Borel--Moore degree zero which
satisfies
\begin{equation}
 \mathbf M=\mathbf I+(s_L)_*\mathbf M+(s_R)_*\mathbf M
 +\Delta_*\bigl(\mathbf C_2-\Sigma_*\mathbf M\bigr).
 \label{eq:universal-meeting-recursion}
\end{equation}
Equivalently, if
\(\mathsf K=(s_L)_*+(s_R)_*-\Delta_*\Sigma_*\), then
\(\mathbf M=\sum_{r\geq0}\mathsf K^r
(\mathbf I+\Delta_*\mathbf C_2)\).  Each application of
\(\mathsf K\) raises total support multiplicity, so this second
description is finite after restriction to \(\mathcal Z_{\leq N}\).

\begin{lemma}
  Writing $(X,\beta_1,\beta_2)$ for the image of $\mathbf1_{\cZ^+(X,\beta_1)\times\cZ^+(X,\beta_2)}$ in $H_c^0(\cZ^+\times_{\CY_4}\cZ^+/\CY_4)$, we have
  \begin{equation*}
    \langle\bfM,(X,\beta_1,\beta_2)\rangle=m_{\beta_1,\beta_2}.
  \end{equation*}
\end{lemma}
\begin{proof}
  By slight abuse, write $m_{\beta_1,\beta_2}=\langle\bfM,(X,\beta_1,\beta_2)\rangle$. It suffices to show that these $m_{\beta_1,\beta_2}$ satisfy the recursive ordered meeting rules of \cite[Section~1.2]{KP}. Choose a basis
  \(S_i\) of \(H^4(X;\mathbb Q)\), put
  \(g_{ij}=\int_XS_iS_j\), and denote the inverse matrix by \((g^{ij})\). For dimensional
  reasons, only the component
  \(\sum_{i,j}g^{ij}S_i\boxtimes S_j\) of the K\"unneth decomposition of
  the diagonal contributes to the incidence class \(\mathbf I\), and we can therefore compute its
  degree on this sector as
  \(\sum_{i,j}n_{0,\alpha}(S_i)g^{ij}n_{0,\beta}(S_j)\).
  
  Pair the entire recursion \eqref{eq:universal-meeting-recursion} above with $(X,\alpha,\beta)$. 
  
  First of all, the only sector
  mapped to the sector of class $(\alpha,\beta)$ by \(s_L\) is \((\alpha,\beta-\alpha)\), and the only one
  mapped to it by \(s_R\) is \((\alpha-\beta,\beta)\).

  Second, the last term is
  supported on the diagonal by definition. On the sector \((\beta,\beta)\), this diagonal
  term contributes \(n_{0,\beta}(c_2(TX))\) by Lemma~\ref{uce:lemma:c_2-class-evaluation} above.
  
  Finally, the coefficient of
  \(\Sigma_*\mathbf M\) in class \(\beta\) is the sum of meeting numbers over the ordered
  decompositions \(\beta_1+\beta_2=\beta\). 
  
  Hence, we obtain
  \begin{equation}
   \begin{aligned}
   m_{\alpha,\beta}
   &={}
   \sum_{i,j}n_{0,\alpha}(S_i)g^{ij}n_{0,\beta}(S_j)
   +m_{\alpha,\beta-\alpha}+m_{\alpha-\beta,\beta},
   &&\alpha\ne\beta,\\
   m_{\beta,\beta}
   &={}
   \sum_{i,j}n_{0,\beta}(S_i)g^{ij}n_{0,\beta}(S_j)
   +n_{0,\beta}(c_2(TX))
   -\sum_{\substack{\beta_1+\beta_2=\beta\\
                      \beta_1,\beta_2>0}}m_{\beta_1,\beta_2}.
   \end{aligned}
   \label{app:eq:kp-meeting-rules}
  \end{equation}
  (As usual, meeting number is understood to vanish if either entry is not positive
  and effective.)  The displayed equations determine the numbers uniquely
  by induction on  degree, and are exactly the ordered meeting rules
  of \cite[Section~1.2]{KP}, finishing the proof.
\end{proof}

\subsubsection{}\label{app:subsec:numerical-specialization}
Now that we have lifted all genus $0$ contributions in \eqref{eq:gv-kp-definitions}, we want to define $\bfN_1$, which evaluates to $n_{1,\beta}$. We write 
\begin{equation*}
  \mathbf M_{\mathrm{tot}}\coloneqq \Sigma_*\mathbf M.
\end{equation*}
For \(d\geq1\), write
\(\sigma(d)=\sum_{a\mid d}a\) and
\(\sigma_2(d)=\sum_{a\mid d}a^2\).  On Borel--Moore classes with
positive support, define the continuous operators
\begin{equation*}
  \mathsf E\coloneqq \sum_{d\geq1}\sigma(d)d^{-1}V_d,\quad \mathsf P\coloneqq \sum_{d\geq1}d^{-1}V_d,\quad \mathsf C\coloneqq \sum_{d\geq1}\sigma_2(d)d^{-1}V_d.
\end{equation*}
The coefficient of \(V_1\) in \(\mathsf E\) is one, so
\(\mathsf E\) has a continuous Dirichlet inverse.  Define
\(\mathbf N_1\) and the transformed connected GW class by
\begin{equation}
\begin{aligned}
 \mathbf N_1&\coloneqq \mathsf E^{-1}\!\left(
 \mathbf G_1^{\mathrm{raw}}+\frac1{24}\mathsf P\mathbf C_2
 -\frac1{24}\mathsf P\mathbf M_{\mathrm{tot}}\right)\in \hat{H}_0^\BM(\cZ/\CY_4),
 \\
 \mathbf{gw}^{\mathrm{tr}}(q)
 &\coloneqq q\mathbf N_0+\mathsf C\mathbf N_1.
\end{aligned}
 \label{eq:universal-gw-transformation}
\end{equation}

\begin{proposition}\label{uce:prop:transformed-theory-coefficients}
We have $\langle \bfN_1,(X,\beta)\rangle=n_{1,\beta}$. Consequently,
the coefficients of charge one and charge zero
of \(\sum_{\beta>0}\langle \mathbf{gw}^{\mathrm{tr}}(q),(X,\beta)+(X,\beta;\gamma_\bullet)\rangle\) are respectively
\[
 \sum_{\beta>0}n_{0,\beta}(\gamma_\bullet)Q^\beta,\text{ and}
 \quad
 \sum_{\beta>0}n_{1,\beta}\log M(Q^\beta).
\]
\end{proposition}

\begin{proof}
We first consider the genus one contribution $\langle \bfN_1,(X,\beta)\rangle$. Rewriting the definition of \(\mathbf N_1\) gives
\begin{equation}\label{uce:eq:n1-def-rewriting}
  \mathbf{GW}_1^{\mathrm{raw}}=\mathsf E\mathbf N_1-\frac1{24}\mathsf P\mathbf C_2+\frac1{24}\mathsf P\mathbf M_{\mathrm{tot}},
\end{equation}
which we pair with $\sum_{\beta>0}(X,\beta)Q^\beta$ to obtain the left-hand side of \eqref{eq:gv-kp-definitions} on the left-hand side. It remains to compare the RHS. 

Identifying the first terms is a simpler version of \S\ref{subsec:connected-gw-inputs} without insertions. We can also read off
\begin{equation*}
  \left\langle \mathsf P (-) , \textstyle\sum_{\beta>0}\mathbf1_\beta Q^\beta\right\rangle=-\left\langle-,\textstyle\sum_{\beta>0}\mathbf1_\beta\log(1-Q^\beta)\right\rangle
\end{equation*}
via the argument in \S\ref{subsec:connected-gw-inputs}. By definition, pairing \(\mathbf M_{\mathrm{tot}}\) with $(X,\beta)$ gives 
\(\sum_{\beta_1+\beta_2=\beta}m_{\beta_1,\beta_2}\), with ordered
summation. Hence, pairing \eqref{uce:eq:n1-def-rewriting} with $\sum_{\beta>0}(X,\beta)Q^\beta$ gives us exactly the second equation in \eqref{eq:gv-kp-definitions}, proving $\langle \bfN_1,(X,\beta)\rangle=n_{1,\beta}$.

Using $\log(M(x))=\sum_{d\geq 1}\sigma_2(d)\frac{x^d}{d}$, we find
\begin{equation*}
  \left\langle \mathsf C (-) , \textstyle\sum_{\beta>0}\mathbf1_\beta Q^\beta\right\rangle=\left\langle-,\textstyle\sum_{\beta>0}\log(M(Q^\beta))\right\rangle.
\end{equation*} 
Hence, pairing
\(\mathbf{gw}^{\mathrm{tr}}(q)
=q\mathbf N_0+\mathsf C\mathbf N_1\) with $(X,\beta)+(X,\beta;\gamma_\bullet)$ yields the claimed coefficients of charge $0$ and $1$.
\end{proof}

\subsection{Disconnected Theories \& Logarithms}
\subsubsection{}\label{subsec:disconnected-theories}
We now pass from the connected transformed GW class to the normalized
theory which allows disconnected domains. For some motivation for the constructions here, see for example \cite[\S3]{pt07} for a numerical ($3$-fold) discussion of connected and disconnected theories, and how exponentials and logarithms relate the two.

If
\(\Sigma_r:\mathcal Z(X/B)^{\times_Br}\to\mathcal Z(X/B)\) denotes
addition, set
\begin{equation}
\begin{aligned}
 \mathbf{GW}_{X/B}^{\mathrm{tr}}(q)
 &\coloneqq \mathbf1_0+\sum_{r\geq1}\frac1{r!}(\Sigma_r)_*
 \bigl(\mathbf{gw}_{\mathrm{GW},X/B}^{\mathrm{tr}}(q)
 \bigr)^{\boxtimes_Br},\\
 \mathbf{PT}_{X/B}(q)
 &\coloneqq \mathbf1_0+\sum_{n\in\mathbb Z}q^n(\rho_{\mathrm{PT},n})_*
 [P_n(X/B)]^{\mathrm{vir}},
\end{aligned}
 \label{eq:transformed-support-addition-exponential}
\end{equation}
where we write \(\mathbf 1_0\) for the Borel--Moore vacuum class on the section of
empty cycles, and use \eqref{vfc:properness} and Proposition~\ref{prop:family-cy4-virtual-package} in the second line.

Both classes defined above are ``disconnected'' theories.  On
the GW side this disconnectedness is introduced by the 
exponential type formula using the cycle addition map, while PT theory ``natively'' allows disconnected supports. 

These systems define classes
which are locally Laurent and take values in
\( \hat{H}_{*}^{\mathrm{BM}}
(\mathcal Z/\CY_4)\llbracket q^{\pm 1}\rrbracket \).  Pairing them with
\(\mathcal H\) defines 
\(\boldsymbol\Phi=(\Phi_n)_{n\in\mathbb Z}\in \Hom_\bQ(\cH,\bQ(\!(q)\!))\) and
\[\boldsymbol\Psi=(\Psi_n)_{n\in\mathbb Z}\in \Hom_\bQ(\cH,\bQ(\!(q)\!))\] by
\[\boldsymbol\Phi(x)\coloneqq \sum_n\Phi_n(x)q^n
=\langle\mathbf{GW}^{\mathrm{tr}}(q),x\rangle,\] and
\[\boldsymbol\Psi(x)\coloneqq \sum_n\Psi_n(x)q^n
=\langle\mathbf{PT}(q),x\rangle.\]
The series \(\boldsymbol\Phi\) and \(\boldsymbol\Psi\) are  unital characters, i.e.
their coefficients satisfy
\(\Phi_n(xy)=\sum_{a+b=n}\Phi_a(x)\Phi_b(y)\) and the analogous identity
for \(\Psi_n\), while the value on the unit is $1$.  

\subsubsection{Convolution Logarithms}
\label{subsec:convolution-logarithms}
Recall the convolution product of \eqref{uce:eq:convolution-product}, and write $\boldsymbol{\epsilon}$ for $\epsilon\cdot q^0$ as a function of $q$. By definition in \eqref{eq:transformed-support-addition-exponential}, both
\(\boldsymbol\Phi-\boldsymbol\epsilon\) and
\(\boldsymbol\Psi-\boldsymbol\epsilon\) vanish on the summand of the
empty cycle.  Their \(r\)-fold convolution powers consequently vanish on
every part of total support multiplicity less than \(r\).  We may
therefore take the convolution logarithm 
\begin{equation*}
  \log^\star(\mathbf Z)\coloneqq
\sum_{r\geq1}(-1)^{r-1}
(\mathbf Z-\boldsymbol\epsilon)^{\star r}/r
\end{equation*}
of the whole coefficient
system and then extract a charge
\begin{equation}
\begin{aligned}
 \mathbf L_\Phi&\coloneqq \log^\star\boldsymbol\Phi,
 &\mathbf L_\Psi&\coloneqq \log^\star\boldsymbol\Psi,\\
 L_{\Phi,n}&\coloneqq [q^n]\mathbf L_\Phi,
 &L_{\Psi,n}&\coloneqq [q^n]\mathbf L_\Psi.
\end{aligned}
 \label{eq:chargewise-convolution-logarithms}
\end{equation}
Note that it the order of these
operations is important. Indeed, an individual coefficient \(\Phi_n\) or
\(\Psi_n\) is generally not a character, and we are \textit{not} considering its logarithm.

It remains to study the difference
\begin{equation}\label{uce:eq:delta-n-log-GWPT-difference}
  \delta_n\coloneqq L_{\Phi,n}-L_{\Psi,n}.
\end{equation}
The following lemma is now a matter of bookkeeping.

\begin{lemma}
\label{lem:normalized-theories-characters}
The systems \(\boldsymbol\Phi\) and \(\boldsymbol\Psi\) are 
unital characters.  Each coefficient \(L_{\Phi,n}\) and
\(L_{\Psi,n}\) of their convolution logarithms is infinitesimal, meaning
\(L(xy)=L(x)\epsilon(y)+
\epsilon(x)L(y)\). Consequently, every
\(\delta_n\) vanishes on \(1\) and on \(\mathfrak m^2\), and hence
induces a linear functional on \(\mathfrak m/\mathfrak m^2\).  Moreover,
\(\mathbf L_\Phi\) is obtained by pairing with the connected class
\(\mathbf g_{\mathrm{GW}}^{\mathrm{tr}}(q)\).
\end{lemma}

\subsubsection{}\label{uce:sec:convolution-log-evaluation}
We now want to compare the evaluation of convolution logarithms of $\mathbf \Phi$ and $\mathbf \Psi$ on elements in $\cH$ to taking logarithms of the evaluations of $\mathbf \Phi$ and $\mathbf \Psi$.

\begin{definition}
  Let $R$ be a graded commutative ring. For some fixed projective CY $4$-fold $X$ and fixed ample divisor $H$, we denote by
  $
    R\llbracket Q\rrbracket
  $
  the ring of series in effective curve classes $\beta\in H_2(X,\bZ)_\eff$, completed in the degree $H\cdot \beta$. Now take $R$ to be a $\bQ$-algebra. We say an element $y=1+\sum_{\beta>0}y_\beta Q^\beta\in \cH\widehat{\otimes} R\llbracket Q\rrbracket$ is \textit{group-like} if $\Delta y=y\otimes y$, understood as the completed tensor product.
\end{definition}

\begin{lemma}\label{uce:lemma:conv-log-eval}
  Let $y=1+\sum_{\beta>0}y_\beta Q^\beta\in \cH\widehat{\otimes} R\llbracket Q\rrbracket$ be a group-like element. Then evaluation at $y$ gives a ring homomorphism
  \begin{equation*}
    \left(\Hom_\bQ(\cH,\bQ(\!(q)\!)),\star,\epsilon\right)\to \left(R(\!(q)\!)\llbracket Q\rrbracket,\cdot,1\right),\ f\mapsto f(y),
  \end{equation*}
  where $f$ is extended $R\llbracket Q\rrbracket$-linearly.
  Consequently, the evaluation of the convolution logarithms of $\boldsymbol\Phi$ and $\boldsymbol\Psi$ on $y$ satisfies
  \begin{equation*}
    \log^\star(\boldsymbol\Phi)(y)=\log(\boldsymbol\Phi(y)),\text{ and}\quad \log^\star(\boldsymbol\Psi)(y)=\log(\boldsymbol\Psi(y)).
  \end{equation*} 
\end{lemma}

\begin{proof}
  Given $f$ and $g$ in $\Hom_\bQ(\cH,\bQ(\!(q)\!))$, the definition \eqref{uce:eq:convolution-product} of the convolution product gives us
  \begin{equation*}
    (f\star g)_n(y)=\sum_{a+b=n}\mu_\bQ((f_a\otimes g_b)\Delta y)=\sum_{a+b=n}\mu_\bQ((f_a\otimes g_b)(y\otimes y))=\sum_{a+b=n}f_a(y)g_b(y),
  \end{equation*}
  where we used the group-like property for the second equality. By \eqref{uce:eq:H-total-mult-exhaustion} every $y_\beta$ lies in some finite $F_{\leq N}\cH$, so that only finitely many terms contribute to each $Q^\beta$ in the logarithm, so we can evaluate it term-wise. The logarithm identity then immediately follows from the convolution product statement, using that \(\boldsymbol\Phi\) and \(\boldsymbol\Psi\) are unital characters.
\end{proof}

\subsubsection{}\label{uce:sec:log-insertion-evaluation}
Now we want to check the properties above for our classes $(X,\beta;\gamma_\bullet)$ from \eqref{eq:primary-insertion-problem}. Fix $X$ a projective CY $4$-fold, and $\gamma_i\in H^{m_i}(X,\bQ)$ with $\sum_i(m_i-2)=2$ and $1\leq i\leq r$. Importantly, assume $m_i>2$ for all $i$, which we will later reduce to using Lemma~\ref{uce:lemma:divisor-insertions}, and assume that the specialization of $\boldsymbol{\Phi}$ and $\boldsymbol{\Psi}$ to $X$ have no terms in negative charge $n<0$. Since $\sum_{\beta>0}(X,\beta;\gamma_\bullet)Q^\beta$ is not grouplike, we introduce instead
\begin{equation*}
  y_{(X,\gamma_\bullet)}\coloneqq \sum_{\beta\geq 0}Q^\beta\mathbf1_\beta\prod_{i=1}^r(1+s_i\hat{\gamma}_i)\in \cH\widehat{\otimes}R\llbracket Q\rrbracket,
\end{equation*}
where $R=\bQ[s_1,\dots,s_r]/(s_i^2)$ with $s_i$ of degree $2-m_i$ used to track degrees of primary insertions. Note that the coefficient of $s_1\cdots s_r$ is $\sum_{\beta>0}(X,\beta;\gamma_\bullet)Q^\beta$, and the coefficient of $s_1^0\cdots s_r^0$ is $1+\sum_{\beta>0}(X,\beta) Q^\beta$.

$\Sigma^{-1}\mathcal Z(X,\beta)=\bigsqcup_{\beta_1+\beta_2=\beta}\mathcal Z(X,\beta_1)\times\mathcal Z(X,\beta_2)$ is a finite disjoint union of open and closed pieces, so $\Sigma^*\mathbf 1_\beta=\sum_{\beta_1+\beta_2=\beta}\mathbf 1_{\beta_1}\boxtimes\mathbf 1_{\beta_2}$. Using this and \eqref{eq:primary-insertion-behaviour}, we find that $y_{(X,\gamma_\bullet)}$ is group-like and hence
\begin{equation*}
  \log^\star(\boldsymbol\Phi)(y_{(X,\gamma_\bullet)})=\log(\boldsymbol\Phi(y_{(X,\gamma_\bullet)})),\qquad \log^\star(\boldsymbol\Psi)(y_{(X,\gamma_\bullet)})=\log(\boldsymbol\Psi(y_{(X,\gamma_\bullet)})).
\end{equation*}
However, both sides are expressions in $s_i$ and $q$, so we have to take coefficients of the correct degrees to extract our quantities of interest. The following argument is the same for $\boldsymbol{\Phi}$ and $\boldsymbol{\Psi}$, so we write it here only for $\boldsymbol{\Phi}$.

We take the coefficients of $q^0s_1^0\cdots s_r^0$ and $q^1s_1\cdots s_r$, respectively. On the left-hand side that gives us
\begin{equation*}
  \begin{aligned}
    [q^0s_1^0\cdots s_r^0]\log^\star(\boldsymbol\Phi)(y_{(X,\gamma_\bullet)})&=\sum_{\beta>0}L_{\Phi,0}(X,\beta)Q^\beta,\\
    [q^1s_1\cdots s_r]\log^\star(\boldsymbol\Phi)(y_{(X,\gamma_\bullet)})&=\sum_{\beta>0}L_{\Phi,1}(X,\beta;\gamma_\bullet)Q^\beta.
  \end{aligned}
\end{equation*}
Using $m_i>2$ and $\sum_i(m_i-2)=2$, we get 
\begin{equation*}
  \boldsymbol\Phi(y_{(X,\gamma_\bullet)})=1+\sum_{\beta>0}Q^\beta\big(\Phi_0(X,\beta)+qs_0\cdots s_r\Phi_1(X,\beta;\gamma_\bullet)\big)
\end{equation*}
with $s^2=0$. Hence, we get
\begin{equation*}
  \log\boldsymbol\Phi(y_{(X,\gamma_\bullet)})=\log\big(1+\textstyle\sum_{\beta>0}Q^\beta\Phi_0(X,\beta)\big) +qs_0\cdots s_r\frac{\sum_{\beta>0}Q^\beta \Phi_1(X,\beta;\gamma_\bullet)}{1+\sum_{\beta>0}Q^\beta\Phi_0(X,\beta)},
\end{equation*}
and we can read off the associated coefficients from the above equations to get
\begin{equation}\label{uce:eq:log-coefficient-identities}
  \begin{aligned}
    \log\bigg(1+\sum_{\beta>0}Q^\beta \Phi_0(X,\beta)\bigg)&=\sum_{\beta>0}L_{\Phi,0}(X,\beta)Q^\beta,\\
    \frac{\sum_{\beta>0}Q^\beta\Phi_1(X,\beta;\gamma_\bullet)}{\sum_{\beta\geq0}Q^\beta\Phi_0(X,\beta)}&=\sum_{\beta>0}Q^\beta L_{\Phi,1}(X,\beta;\gamma_\bullet).
  \end{aligned}
\end{equation}
\section{Regularizing Deformations}\label{sec:regdef}

In this section, we explain a proof of Pardon's generic transversality result. Our methods mimic Pardon's in a lot of ways. However, we modified the proof to be more accessible to us as algebraic geometers. This version of the proof was developed while working on an orbifold version in \cite{bryan-thimm}. 

\subsection{Definitions \& Tools}

\subsubsection{}
First, we will consider deformations of a neighborhood of a curve $C$ in $X$ that regularizes an associated stable map given by normalization. Afterwards, we consider the more complicated problem of deforming a family of manifolds $X\to B$ so that an entire family of curves gets regularized. Before fully considering either problem, we fix some definitions and recall some of Pardon's tools.

The first step when eventually applying this to $H_c^*(\cZ/\CY_4)$ is to use Pardon's graph trick, which essentially extends the base of the family to include a copy of the cycle space, such that we may regularize each cycle individually rather than all of them at once.

\begin{lemma}[Graph Trick\ {{\cite[Lemma 5.14]{pardon2023}}}]\label{lemma:graph-trick}
    Every class in $H_c^*(\cZ/\CY_4)$ is in the image of $H_K^*(\cZ(X/B)/B)$ for some family of CY$4$-folds $X\to B$ over a finite simplicial complex $B$ and a compact subanalytic set $K\subseteq \cZ(X/B)$ whose map to $B$ is injective.
\end{lemma}

\subsubsection{}
Now, we define the type of regularity that we try to achieve by deformation. First, we consider regularity of maps.
\begin{definition}\label{gt:def:regular}
    Let $f:C\to X$ be a holomorphic map from a proper smooth curve $C$. We call a point $x\in C$ \textit{special} (for $f$) if $\rd f(x)=0$ or $\#f^{-1}(f(x))>1$. Let $A\subseteq C$ be the set of special points. If $f$ is simple (nowhere a multiple cover), then $A$ is finite.
    
    We now assume $f$ is simple. Consider the problem of deforming the triple $(C,A,f)$ subject to the constraint that the points $A$ remain special with the same discrete data, meaning that all conditions $f(x)=f(x')$ and $(D^rf)(x)=0$ which hold for $f$ are preserved. We say that the map $f$ is \textit{regular} if this deformation problem is unobstructed.

    Regularity is also defined for curves in fibers of a family $X\to B$, meaning that a map $f_b:C_b\to X_b$ to a fiber is called \textit{regular} if $T_bB$ surjects onto the obstruction space and higher obstructions vanish.
\end{definition}

\begin{remark}
    Regularity is effectively unobstructedness of $f$ viewed as a stable map, but subject to additional point constraints.
\end{remark}

\subsubsection{}
Since deformations of cycles are not generally well-behaved, we associate a ``stable map'' to every cycle and then use the notion of regularity above. This leads to the following notion of semi-regularity, which is first defined for manifold bases and then generalized to simplicial bases.
\begin{definition}\label{gt:def:sr-complex}
    Let $X\to B$ be a family of complex manifolds over $B$.
    \begin{enumerate}[label=(\alph*)]
        \item For a $1$-cycle $z=\sum_i m_i[C_i]$ in $\cZ(X/B)$, with distinct and irreducible $C_i$, define the \textit{associated embedded curve} $u(z):C\coloneqq \bigcup_i C_i\hookrightarrow X$, where $C$ is equipped with its reduced structure. Define the \textit{associated stable map} $f(z):\tilde{C}\coloneqq \sqcup_i \tilde{C}_i\to X$ using normalizations of the components. When we want to record the fiber $X_b$ in which $z$ sits, we write $u_b(z):C_b\hookrightarrow X_b$ and $f_b(z):\tilde{C}_b\to X_b$ instead.
        \item If $B$ is a complex (or real) manifold, we call a $1$-cycle $z=\sum_im_iC_i\in\cZ(X/B)$ \textit{semi-regular} if the associated stable map $f(z):\bigsqcup_i\tilde C_i\to X$ is regular (in the family $X\to B$) in the sense of Definition~\ref{gt:def:regular}.
        \item If $B$ is a simplicial complex, we call a $1$-cycle $z=\cZ(X/B)$ \textit{semi-regular} when it is semi-regular (in the sense of the definition for smooth bases just above) as a point of $\cZ(X_{\sigma^\circ}/\sigma^\circ)$ where $\sigma\subseteq B$ is the smallest real simplex containing the image of $z$, i.e. the unique simplex whose interior contains the image of $z$.
    \end{enumerate}
\end{definition}

\subsubsection{}
First, we discuss the following lemma, which helps us work with the constrained deformation problem above. Intuitively, since the special constraints in Definition~\ref{gt:def:sr-complex} come from normalizing $C$, its obstructions lift all obstructions to deforming $\tilde{C}$ under these constraints.
\begin{lemma}
    Consider the associated embedded curve $u:C\to X$ of a cycle. Then we have a surjection
    \begin{equation}\label{gt:eq:emb-curve-ob-surjects-onto-constrained-surjection}
        H^1(C,u^*T_X)\twoheadrightarrow \rOb_A(f)
    \end{equation}
    onto the constrained deformation problem of the associated stable map in Definition~\ref{gt:def:sr-complex} above.
\end{lemma}
\begin{proof}
    Let us study the deformation theory of the associated stable map $f$ given by composition with the normalization $\nu:\tilde{C}\to C$. Let $A\subset \tilde{C}$ denote the special points, which lie over the singular points of $C$. Now the deformation theory of the marked map $(\tilde{C},A,f)$ is determined by the complex $\big[T_{\tilde{C}}(-A)\xrightarrow{\rd f} f^*T_X\big]$ or equivalently $\big[\nu_*T_{\tilde{C}}(-A)\to \nu_*f^*T_X\big]$ as $\nu_*$ is exact.
    
    To account for the additional constraints in Definition~\ref{gt:def:sr-complex}, let $F\subseteq\nu_*\cO_{\tilde C}$ be the coherent $\cO_C$-submodule of functions $h$ such that for every $y\in\mathrm{Sing}(C)$ there is a constant $c\in\bC$ with
    \begin{equation*}
        h-c\in\fm_x^{\ord_x(\rd f)+1}\quad\text{for all }x\in\nu^{-1}(y).
    \end{equation*}
    That is, $h$ takes the same value at all points over $y$, and its derivatives of order $1,\dots,\ord_x(\rd f)$ vanish at each such $x$. These are the analogues for $h$ of the conditions $f(x)=f(x')$ and $(D^rf)(x)=0$ imposed on $f$. By slight abuse of notation, we label the $x$ in $\nu^{-1}(y)$ by $\{x_1,\dots,x_m\}$ with orders $r_i$. We can then describe the cokernel $Q$ of the inclusion $F \subseteq \nu_*\cO_{\tilde{C}}$ as a sum of skyscraper sheaves
    \begin{align*}
        Q&\cong \bigoplus_{y\in \mathrm{Sing}(C)}\left[\bC^{\abs{\nu^{-1}(y)}-1}\oplus \bigoplus_{x\in \nu^{-1}(y)} \bC^{\ord_x(\rd f)}\right],\\ 
        [h]&\mapsto \Big(\big(h(x_1)-h(x_j)\big)_{j>1},\big(h^{(k)}(x_i)\big)_{1\leq i\leq m,1\leq k\leq r_i}\Big).
    \end{align*}
    Tensoring the associated short exact sequence by $u^*T_X$ and using the projection formula gives
    \begin{equation*}
        0\to G_A \to \nu_*f^*T_X \xrightarrow{\lambda} J\coloneqq Q\otimes_{\cO_C} u^*T_X \to 0, 
    \end{equation*}
    where $J$ is again just a skyscraper sheaf at the singular points of $C$ and $G_A$ is the kernel. The quotient map amounts to linearizing the constraint equation along local sections of $\nu_*f^*T_X$.
    
    Since $\rd f$ vanishes at $x_i$ to $r_i$th order and sections of $T_{\tilde{C}}(-A)$ vanish at points in $\nu^{-1}(\mathrm{Sing}(C))$, so that $\nu_*\rd f$ lands in $G_A$. The complex controlling the deformations of the constrained problem is hence $\big[\nu_*T_{\tilde{C}}(-A)\to G_A\big]$, which gives us a surjection
    \begin{equation*}
        H^1(C,G_A)\twoheadrightarrow \rOb_A(f)
    \end{equation*} 
    by using cohomological dimension $1$ of $\tilde{C}$.
    
    Now let $g\in \cO_{C,y}$ be a function on $C$ near $y\in\mathrm{Sing}(C)$. This is the restriction of a function $G\in \cO_{X,y}$ on $X$ near $y$, so that $g\circ\nu=G\circ f$ near each $x\in\nu^{-1}(y)$. The differential $\rd(G\circ f)$ vanishes to order at least $\ord_x(\rd f)$ at $x$, and $(G\circ f)(x)=g(y)$ by construction. So, $g\circ\nu-g(y)\in\fm_x^{\ord_x(\rd f)+1}$ for every $x\in\nu^{-1}(y)$, which exactly tells us that $g\circ\nu$ is in $F$ with $c_y=g(y)$. Hence, we proved that $\cO_C\subseteq F\subseteq \nu_*\cO_{\tilde{C}}$, so that $u^*T_X=\cO_C\otimes_{\cO_C}u^*T_X\subseteq F\otimes_{\cO_C}u^*T_X=G_A$ with cokernel $Q'$ of $0$-dimensional support at the singular points of $C$, since elsewhere $\nu$ is an isomorphism. Hence, the long exact sequence of $0\to u^*T_X\to G_A\to Q'\to 0$ gives us
    \begin{equation*}
        \dots \to H^1(C,u^*T_X)\to H^1(C,G_A)\to H^1(C,Q')=0.
    \end{equation*}
    Composing the two surjections proves the lemma.
\end{proof}

\subsubsection{}
In order to later control the chosen points for the neighborhoods $U_C$ of Lemma~\ref{gen:lemma:cohomological-dimension-1-neighborhood} in families of curves, we need to recall the notion of controlling divisors from \cite{pardon2023}.

\begin{definition}
    Let $D\subseteq X\times_B U$ be a relative divisor over an open subset $U$ in $B$. Given a cycle $z$ in $\cZ(X/B)$ in a fiber $X_b$ with $b\in U$, we say $D$ \textit{controls} $z$ if $D_b$ intersects $z$ in a finite set and meets every irreducible component of $z$.
\end{definition}

This is an open condition and Pardon shows that we can always obtain such divisors.

\begin{proposition}[Enough Divisors\ {{\cite[Prop. 5.10]{pardon2023}}}]\label{gt:prop:enough-divisors}
    Let $X\to B$ be a family of complex manifolds over a finite simplicial complex, and let $K\subseteq \cZ(X/B)$ be a compact analytic set for which the map $K\to B$ is injective. After possibly removing a closed subset from $X$ disjoint from $K$, there exists a finite collection of open subsets $U_i\subseteq B$ and disjoint relative divisors $D_i\subseteq X\times_B U_i$ such that every $z\in K$ is $D_i$-controlled for some $i$.
\end{proposition}

\subsection{Single Cycle Regularization}
\label{gt:sec:inf-def}

\subsubsection{}
The strategy to show generation by local curve elements is to first prove that we can restrict to so-called semi-regular elements. We will now do a ``warm-up'' for one (not necessarily irreducible) curve, which we should view as the union of all components of a cycle without multiplicity. The first key technical result is the following, showing that deformations of a neighborhood of such a curve surject onto obstructions of the map from the curve.

\begin{lemma}\label{gen:lemma:infinitesimal-obstructions-surjection}
    Let $u:C\to X$ be an embedding of a proper curve $C$ into a complex manifold $X$. Then there exists a neighborhood $U_C$ of the curve, such that the map
    \begin{equation*}
        H^1(U_C,T_X) \to H^1(C,u^*T_X)
    \end{equation*}
    is surjective.
\end{lemma}

We diverge here from Pardon's approach in \cite[Lemma 5.6]{pardon2023}, since the following presentation seemed more natural to us as algebraic geometers. Additionally, the importance of cohomological dimension in this version of the Lemma is interesting.

\subsubsection{}

The proof consists of two main steps:
\begin{enumerate}
    \item Find a neighborhood $U_C$ of the curve of cohomological dimension $1$.
    \item Examine the natural exact sequences associated to the embedding $u$ and use the cohomological dimension $1$ to prove the desired surjection.
\end{enumerate}

\begin{lemma}\label{gen:lemma:cohomological-dimension-1-neighborhood}
    Let $u:C\to X$ be an embedding of a proper curve $C$ into a complex manifold $X$. Given a finite collection $P=\{p_1,\dots,p_k\}$ of distinct points on the curve $C$, such that each irreducible component of $C$ contains at least one of them, we obtain a neighborhood $U_C$ of $C$ in $X$ of cohomological dimension $1$.
\end{lemma}

\begin{figure}[!ht]
  \centering
  \includegraphics[scale=1.1]{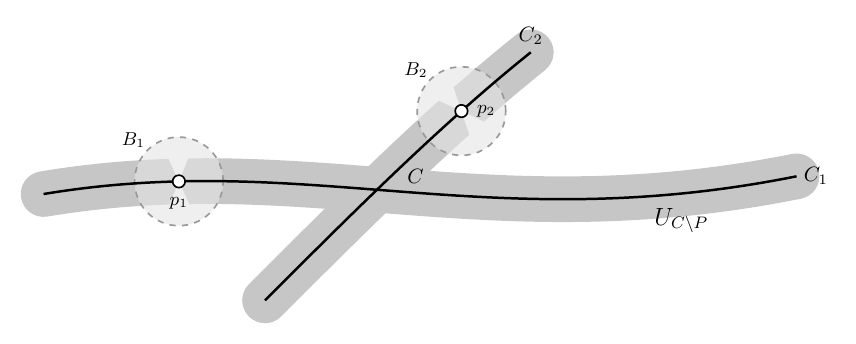}
  \caption{A cohomological dimension $1$ neighborhood of a curve $C$ given by the union of a Stein neighborhood $U_{C\setminus P}$ of the curve with a point removed per irreducible component, and open balls around these points.}
  \label{gen:fig:cohomological-dimension-1-nbd}
\end{figure}

\begin{proof}
    We construct $U_C$ as the union of two Stein open subsets, i.e. two open sets of cohomological dimension $0$ by Cartan's Theorem B. See Figure~\ref{gen:fig:cohomological-dimension-1-nbd} above for an illustration of the process. First note, that $C\setminus P$ is Stein by \cite[Theorem $\mathrm{(A)_1}$]{siu-analytic-sheaf-cohom}. By \cite{Siu1976}, $C\setminus P$ then has an open Stein neighborhood $U_{C\setminus P}$ in $X$.

    Since $P$ contains distinct points $p_i$, we can find sufficiently small open balls $B_i$ around each point in $P$, so that they don't intersect. The union $U_P$ of those balls is Stein, because any open ball is Stein.

    The union $U_C\coloneqq U_{C\setminus P}\cup U_P$ is then a neighborhood of the curve $C$ and has cohomological dimension $1$ as a union of two Stein open subsets.
\end{proof}

\subsubsection{}\label{gen:sec:proof-infinitesimal-obstruction-surjection}

\begin{proof}[Proof of Lemma~\ref{gen:lemma:infinitesimal-obstructions-surjection}]
    We consider the short exact sequence given by the ideal sheaf exact sequence of $u$ tensored by $T_X$
    \begin{equation*}
        0 \to I_C\otimes T_X \to T_X \to u_*u^*T_X \to 0,
    \end{equation*}
    where we used the projection formula. We can restrict this to the open neighborhood $U_C$ of $C$ of cohomological dimension $1$ from Lemma~\ref{gen:lemma:cohomological-dimension-1-neighborhood} above. Then consider the associated long exact sequence.
    \begin{equation}\label{gen:eq:cohomological-dimension-1-surjection}
        \cdots \to H^1(U_C,I_C\otimes T_X)\to H^1(U_C,T_X) \to H^1(C,u^*T_X) \to H^2(U_C,I_C\otimes T_X)=0,
    \end{equation}
    where the last term vanishes due to $U_C$ being cohomological dimension $1$, making the map $H^1(U_C,T_X) \to H^1(C,u^*T_X)$ surjective. 
\end{proof}

\subsubsection{}

\begin{lemma}\label{gen:lemma:cy-preservation-cohomological-dim-1}
    Let $X\to B$ be a family of CY manifolds of relative dimension $n\geq 3$ with a cover by two open subsets with Stein fibers $U_1$, $U_2$ over $B$ such that their intersection $U$ is Stein. Let $u:C\to X_b$ be an embedding of a proper reduced curve $C$ into a fiber. Then we have a surjection
    \begin{equation*}
        H^0_\rCY(U,T_{U/B}) \twoheadrightarrow H^1(C,u^*T_{X_b}).
    \end{equation*}
\end{lemma}

\begin{proof}
    Now first that we may view $H^0_\rCY(U,T_{U/B})$ as simply the closed forms in $H^0(U,\Omega^{n-1}_{U/B})$ via the isomorphism $T_{U/B}\cong \Omega_{U/B}^{n-1}$ induced by $\iota_{(-)}\omega$. We are going to realize a surjection onto $H^1(C,u^*T_{X_b})$ from a space of exact forms instead, which implies the desired surjection. First, we study the curve $C$ in the fiber $X_b$. Consider the map
    \begin{equation*}
        \mu_b:\cN^\vee_{C/X_b}\otimes \Omega_{X_b}^{n-2}|_C\to \Omega_{X_b}^{n-1}|_C,\ [g]\otimes \beta|_C\mapsto (\rd g\wedge \beta)|_C.
    \end{equation*}
    While this contains a differential operator, the restriction to $C$ makes this well-defined and $\cO_C$-linear, since the first term in $\mu_b(h[g]\otimes \beta|_C)=\big((g\rd h+h\rd g)\wedge \beta\big)|_C$ has coefficients in $I_C$ and hence vanishes. Around a smooth point of $C$, we can pick coordinates $z_1,\dots, z_n$ such that $C$ is cut out by $z_2=\cdots =z_n=0$. In these coordinates, we can lift every element $\bigwedge_{k\neq i}\rd z_k$ of the standard basis of $\Omega_{X_b}^{n-1}|_C$ to $\pm [z_j]\otimes \bigwedge_{k\neq i,j}\rd z_k$ if by choosing $j\neq i$ using $n\geq 3$. Hence, $\mu_b$ is surjective at smooth points of $C$ and the inclusion of
    \begin{equation*}
        \cI\coloneqq \im(\mu_b)\subseteq \Omega_{X_b}^{n-1}|_C\cong u^*T_{X_b}
    \end{equation*}
    has $0$-dimensional cokernel supported at the singular points. Taking long exact sequences of the corresponding quotient short exact sequence hence yields surjectivity of the bottom map in the commutative diagram
    \begin{equation}\label{gt:eq:fiber-curve-surjection-diagram}
        \begin{tikzcd}[row sep=small]
            H^0(C\cap U_b,\cI|_{U_b})\ar[r]\ar[d,twoheadrightarrow] & H^0(C\cap U_b,u^*T_{X_b})\ar[d,twoheadrightarrow]\\
            H^1(C,\cI)\ar[r,twoheadrightarrow] & H^1(C,u^*T_{X_b}),
        \end{tikzcd}
    \end{equation}
    where the vertical arrows are the ones induces by the Stein \v Cech cover $U_{i,b}$. 

    Now, we want to restrict to the fiber $U_b$ used above from the total space $U$. Writing $I_C$ for the ideal of $C$ in $X_b$ and $J_C$ for the ideal of $C\cap U_b\xhookrightarrow{j}U$, consider the following map
    \begin{equation*}
        \rho: J_C\Omega_{U/B}^{n-2}\to j_*(\Omega_{U_b}^{n-1}|_{C\cap U_b}),\ \beta\mapsto \rd\beta|_{C\cap U_b},
    \end{equation*}
    where the restriction denotes the adjunction map of $j$. Similar to before, the restriction to $C\cap U_b$ makes this map well-defined and $\cO_U$-linear. Because restriction to open subsets and $j_*$ are exact, we have
    \begin{equation*}
        j_*(\cI|_{U_b})=\im(j_*(\mu_b|_{U_b}))=\im(j_*(\cN^\vee_{C\cap U_b/U_b}\otimes \Omega^{n-2}_{U_b}|_{C\cap U_b})\to j_*(\Omega_{U_b}^{n-1}|_{C\cap U_b})).
    \end{equation*}
    Since $J_C$ is simply the pre-image of $I_C|_{U_b}$ under the surjection $\cO_U\twoheadrightarrow \cO_{U_b}$, we get a surjection $J_C\twoheadrightarrow j_*(\cN^\vee_{C\cap U_b/U_b})$, which tensored with $\Omega_{U/B}^{n-2}$ yields the surjection in the following commutative diagram
    \begin{equation*}
        \begin{tikzcd}
            J_C\Omega_{U/B}^{n-2}\ar[d,twoheadrightarrow]\ar[dr,"\rho"] & \\
            j_*(\cN^\vee_{C\cap U_b/U_b}\otimes \Omega^{n-2}_{U_b}|_{C\cap U_b})\ar[r,"j_*(\mu_b|_{U_b})"] & j_*(\Omega_{U_b}^{n-1}|_{C\cap U_b}).
        \end{tikzcd}
    \end{equation*}
    Hence, we get $\im(\rho)=j_*(\cI|_{U_b})$. Since $U$ is Stein, the cohomology long exact sequence associated to $\rho:J_C\Omega_{U/B}^{n-2}\twoheadrightarrow j_*(\cI|_{U_b})$ yields surjectivity of the left vertical map in the diagram
    \begin{equation}\label{gt:eq:open-exact-surjection}
        \begin{tikzcd}[row sep=small]
            H^0(U,J_C\Omega_{U/B}^{n-2})\ar[r,"\rd_{U/B}"]\ar[d,twoheadrightarrow,"\rho"] & H^0_\rCY(U,T_{U/B})\ar[r] & H^0(U_b,T_{U_b})\ar[d]\\
            H^0(C\cap U_b,\cI|_{U_b})\ar[rr] & & H^0(C\cap U_b,u^*T_{X_b}),
        \end{tikzcd}
    \end{equation}
    which commutes by construction. Composing this diagram with \eqref{gt:eq:fiber-curve-surjection-diagram} yields the desired surjection.
\end{proof}

\subsubsection{}\label{gt:sec:single-curve-ks-strategy}
Let $u:C\to X$ be an embedding of a proper curve $C$ into a complex CY manifold $X$, let $n:\tilde{C}\to C$ be its normalization, and write $f\coloneqq u\circ n$. As motivation for the family version below, we now construct a deformation of $X$ in which $f$ is unobstructed subject to the given point constraints.

Using \eqref{gt:eq:emb-curve-ob-surjects-onto-constrained-surjection} and Lemma~\ref{gen:lemma:cy-preservation-cohomological-dim-1} (with $B=\pt$), we can take a finite dimensional (complex vector space viewed as a) real vector space $V\subseteq H^0_\rCY(U_{C\setminus P}\cap U_P,T_X)$ which surjects onto $\rOb_A(f)$. Using the stronger version \cite[Thm. 1 ($q=1$)]{Demailly1990} of \cite{Siu1976}, we may arbitrarily shrink the Stein pieces of our neighborhood $U_C$, so that Proposition~\ref{bg:prop:family-regluing} gives us a CY (by Corollary~\ref{bg:cor:regluing-cy-preservation}) family
\begin{equation*}
    \widetilde X\to B_0V
\end{equation*}
whose Kodaira--Spencer map $V\to \rOb_A(f)$ is exactly the surjection above by construction using regluing along flows of the vector fields in the complex span of $V$. Hence, $f$ is unobstructed in $\widetilde X/B_0V$ subject to the constraints imposed by Definition~\ref{gt:def:sr-complex}.

\subsection{Generic Transversality}\label{gt:sec:gt}
Now we want to use Proposition~\ref{bg:prop:family-regluing} to prove the main technical Theorem~\ref{gt:thm:generic-transversality-in-families} of this Section, which is a version of the procedure in \S\ref{gt:sec:inf-def} for families of cycles. Specifically, we work in the situation produced by the graph trick above.

\subsubsection{}
While the divisors supplied by Pardon's ``Enough Divisors'' Proposition~\ref{gt:prop:enough-divisors} control the points chosen in the construction of the cohomological dimension $1$ neighborhood above, it is difficult to control to Stein neighborhoods of the complement given by \cite{Siu1976} in families. Instead, we note that the proof of Lemma~\ref{gen:lemma:cy-preservation-cohomological-dim-1} above works verbatim only using weaker conditions, which puts us into a position to simply use Pardon's Sard--Smale argument, since he only works with his analogue of the balls around the chosen points.

\begin{lemma}\label{gt:lemma:cy-preserving-surjection-only-overlap}
    Let $X\to B$ be a family of CY manifolds of relative dimension $n\geq 3$ with a cover by two open subsets $U_1$, $U_2$ over $B$ such that their intersection $U$ is Stein. Let $u:C\to X_b$ be an embedding of a proper reduced curve $C$ into a fiber, such that $C\cap U_{i,b}$ is Stein for both $i$. Then we have a surjection
    \begin{equation*}
        H^0_\rCY(U,T_{U/B}) \twoheadrightarrow H^1(C,u^*T_{X_b}).
    \end{equation*}
\end{lemma}
\begin{proof}
    Lemma~\ref{gen:lemma:cy-preservation-cohomological-dim-1} above works verbatim. We note here, where the conditions are used. The condition that $C\cap U_{i,b}$ is Stein for both $i$ gives us the vertical surjections in \eqref{gt:eq:fiber-curve-surjection-diagram}. The condition that $U$ is Stein gives the surjectivity of the left vertical map in \eqref{gt:eq:open-exact-surjection}.
\end{proof}

\subsubsection{}
\begin{theorem}[Generic Transversality in Families\ {{\cite[Prop. 5.9]{pardon2023}}}]\label{gt:thm:generic-transversality-in-families}
    Let $X\to B$ be a family of CY manifolds of dimension $n\geq 3$ over a finite simplicial complex. Let $D_i$ be a finite collection of disjoint relative divisors over open subsets $U_i\subseteq B$, together with sets $A_i\subseteq D_i$ (our ``divisor cycle intersections'') whose map to $U_i$ is proper with finite fibers. After ``trimming'' $X$ by removing closed subsets contained in $D_i\setminus A_i$, and subdividing $B$, there exists a Banach space $\cJ_D^\rCY(X/B)$ of CY preserving deformations of $X\to B$, a generic point $\tilde{X}\to B$ of which has the property that any stable map associated to a $D_i$-controlled $1$-cycle of $\tilde{X}/B$ is regular.

    More generally, suppose $D_i$ are disjoint relative divisors over subcomplexes $V_i\subseteq B$ and suppose $U_i\subseteq V_i\subseteq B$ are constructible open sets. Then after trimming $X$ over all $V_i$ and subdividing $B$, there exists a Banach space $\cJ_D^\rCY(X/B)$ of CY preserving deformations of $X\to B$, vanishing over the complement of $\bigcup_i U_i$, a generic point $\tilde{X}\to B$ of which has the property that any stable map over $\bigcup_i U_i$ associated to a $D_i$-controlled $1$-cycle of $\tilde{X}/B$ is regular.
\end{theorem}

\begin{proof}
With the technical result of Lemma~\ref{gt:lemma:cy-preserving-surjection-only-overlap} and \eqref{gt:eq:emb-curve-ob-surjects-onto-constrained-surjection} (in place of Pardon's \cite[Lemma 5.6]{pardon2023}), Pardon's proof of generic transversality goes through without changes. We sketch it here for convenience. First, we find charts of polydiscs
\begin{equation}\label{gt:eq:gt-charts}
    \begin{tikzcd}
        (\bD,0)\times \bD^{n-1} \times M_{ij} \ar[r,hookrightarrow,"\phi_{ij}"]\ar[d] & (X|_{U_i},D_i)\ar[d]\\
        M_{ij}\ar[r,hookrightarrow] & U_i,
    \end{tikzcd}
\end{equation}
where $M_ij$ are subcomplexes and $\bD$ is the closed unit disc, so that $0\times (\bD^{n-1})^{\circ}\times M_{ij}^\st$ jointly cover $A_i$.\footnote{Note that we can avoid Pardon's more complicated condition that the boundary strips of the charts avoid $A_i$ since this assumption is not needed for our surjection onto the constrained obstruction space.} Since the $D_i$ are disjoint, we may assume that these charts are also disjoint on the overlaps of various $U_i$. Here, $M_{ij}^\st$ is Stein, by Lemma~\ref{bg:lemma:constructible-stein}, and by Proposition~\ref{bg:prop:bump-functions}, we may assume that there are bump functions vanishing on the boundary of $M_{ij}^\st$ after subdivision. We then trim $X$ by removing the part of each $D_i$ outside the charts $0\times (\bD^{n-1})^\circ$, which requires the divisors to be disjoint.

With these charts in place, the regluing space $\cJ_D^\rCY(X/B)$ is now built from $H_\rCY^0((\bD^\circ\setminus 0)\times \bD^{n-1}\times M_{ij}^\st,T_{X/B})$. Using the bump functions on $M_{ij}^\st$ above and that the $D_i$ are disjoint, we can add vector fields from the different parts and reglue using Proposition~\ref{bg:prop:family-regluing} and Corollary~\ref{bg:cor:regluing-cy-preservation}. $\cJ_D^\rCY(X/B)$ is made into a Banach space by equipping the pieces an $L^2$-norm exponentially weighted near the punctures, and taking direct sums. 

Now, any $D_i$-controlled cycle $z$ with associated embedded curve $C\subset X_b$ satisfies that $C\cap (X\setminus D_i)$ is Stein by \cite[Theorem $\mathrm{(A)_1}$]{siu-analytic-sheaf-cohom}, $C\cap (\bD^{n-1})^\circ$ is Stein as a closed subset of a Stein space $(\bD^{n-1})^\circ$, and $\bD^\circ\setminus 0\times \bD^{n-1}\times M_{ij}^\st$ is Stein. Lemma~\ref{gt:lemma:cy-preserving-surjection-only-overlap} (using $n\geq 3$) and \eqref{gt:eq:emb-curve-ob-surjects-onto-constrained-surjection} now give the surjection 
\begin{equation}\label{gt:eq:composed-surjection-constrained-obstructions-gt-proof}
    H^0_\rCY((\bD^\circ\setminus 0)\times \bD^{n-1}\times M_{ij}^\st,T_{X/B}) \twoheadrightarrow H^1(C,u^*T_{X_b}) \twoheadrightarrow \rOb_A(f)
\end{equation}
onto the constrained obstruction space of the associated stable map. Now note that any simple $D_i$-controlled map is the associated stable map of its image cycle, which is $D_i$-controlled.

Consider a compact smooth (possibly disconnected) surface $C$ together with point constraints $\gamma$, given by a finite set $A_\gamma\subseteq C$, a partition of $A_\gamma$ and orders $e_p\geq 0$ for $p\in A_\gamma$. Take $\cJ(C)$ to be a finite-dimensional smooth manifold parametrizing a family of almost complex structures on $C$ which agree near $A_\gamma$ and contain all small deformations of the marked curve $(C,A_\gamma)$, and fix holomorphic coordinates $z$ near the points of $A_\gamma$. For $k\geq 2+\max_p e_p$ (so that $W^{k,2}\subseteq C^{e_p}$), the space $W^{k,2}(C,X)$ is a smooth Banach manifold, see \cite[\S3.1]{wendl-symplectic-lectures}, and the maps satisfying the constraints
\begin{equation*}
    u(p)=u(p')\ \text{ for $p,p'$ in the same part},\qquad \partial_z^r u(p)=0\ \text{ for } 1\leq r\leq e_p,
\end{equation*}
computed in local holomorphic coordinates on $X$, form a Banach submanifold $W^{k,2}_\gamma(C,X)$ of finite codimension. Let $W^{k,2}_\gamma(C,X)^+$ be the closed and open subset of maps whose restriction to every component of $C$ has positive algebraic intersection with $D_i$ for a fixed $i$ (holomorphic maps are in here if and only if they are $D_i$-controlled). Restricting the smooth Banach bundle of \cite[(5.3)]{pardon2023} gives a bundle over $\cJ(C)\times W^{k,2}_\gamma(C,X)^+\times\cJ^\rCY_D(X/B)$ with a section $\bar\partial$ measuring the failure of $(j,u,J)$ to be holomorphic. At a holomorphic triple, the linearization of $\bar\partial$ in the $(j,u)$-directions is the linearized deformation problem of $(C,A_\gamma,u)$ subject to $\gamma$, so its cokernel is the constrained obstruction space $\rOb_\gamma(u)$. 

If $u$ is simple, it is the associated stable map of its image cycle, so its special points are $A=u^{-1}(\mathrm{Sing}\,u(C))$ and the constraints of Definition~\ref{gt:def:regular} at $A$ are exactly those of the normalization construction in \eqref{gt:eq:emb-curve-ob-surjects-onto-constrained-surjection}. Since $u\in W^{k,2}_\gamma(C,X)$, the constraints $\gamma$ are among these, and forgetting the others gives a surjection $\rOb_A(u)\twoheadrightarrow\rOb_\gamma(u)$. Since $W^{k,2}_\gamma(C,X)$ does not depend on $J$, the $\cJ^\rCY_D(X/B)$-directions act on this cokernel through the Kodaira--Spencer map, so the surjection \eqref{gt:eq:composed-surjection-constrained-obstructions-gt-proof} shows that $\bar{\partial}$ is transverse to $0$ at every $D_i$-controlled simple holomorphic triple $(j,u,J)$, making
\begin{equation*}
    \bar{\partial}^{-1}(0)_{\mathrm{simple}}\to \cJ_D^\rCY(X/B)
\end{equation*}
a Fredholm map from a smooth Banach manifold, see \cite[(5.5)]{pardon2023}. Sard--Smale \cite{Smale1965} now implies that the fibers of this map over generic points are regular for the constraints $\gamma$. Covering all curves and all possible special data by countably many triples $(C,\cJ(C),\gamma)$, and taking $\gamma$ to be the special data of the map in question, we find that for generic deformations $\tilde{X}/B$ in $\cJ_D^\rCY(X/B)$ all stable maps associated to a $D_i$-controlled $1$-cycle of $\tilde{X}/B$ are regular.
\end{proof}

\subsection{Generic Transversality and the Pardon Ring}
\subsubsection{The Semi-Regular Group}\label{gt:subsubsec:semi-regular-group}
A semi-regular version of the Pardon group is defined in \cite[Def. 5.13]{pardon2023}. As discussed there, this applies to the interior semi-regular locus $\cZ^{\sr\circ}$, but also $\cZ^{\sr\circ}_\alpha$ for any locally closed sublocus $\cZ_\alpha\subseteq \cZ$ satisfying the conditions of \cite[Def. 3.13]{pardon2023}. In particular, by Proposition~\ref{rcl:prop:cohomology-fiberwise-cond}, this applies to fiberwise local subloci, which includes the locus $\cZ_\con$ of connected cycles and the locus $R_m^\sm$ of smooth embedded rational curves of multiplicity $m$, which is defined in Definition~\ref{rnd:def:rational-smooth-locus} below. For such a fiberwise local sublocus $\cZ_\alpha$, \cite[Def. 5.13]{pardon2023} then defines a group
\begin{equation*}
    H_c^*(\cZ_\alpha^{\sr\circ}/\CY_4)\to H_c^*(\cZ_\alpha/\CY_4),
\end{equation*}
with the map being an extension by zero. This is the map, which is shown to be an isomorphism by the following CY$4$ version of Pardon's Theorem \cite[Theorem 5.15]{pardon2023}.

\subsubsection{}
\begin{theorem}[{\cite[Theorem 5.15]{pardon2023}}]\label{gt:thm:pardon-ring-sro-isom}
    The map 
    \begin{equation*}
        H_c^*(\cZ^{\sr\circ}/\CY_4)\to H_c^*(\cZ/\CY_4)
    \end{equation*}
    is an isomorphism. Additionally, the same holds for the corresponding map $H_c^*(\cZ^{\sr\circ}_\alpha/\CY_4)\to H_c^*(\cZ_\alpha/\CY_4)$ for any fiberwise local sublocus $\cZ_\alpha\subseteq \cZ$ in the sense of Definition~\ref{uce:def:fiberwise-local} (or more generally satisfying \cite[Def. 3.13]{pardon2023}).
\end{theorem}

\begin{proof}
    Using generic transversality in Theorem~\ref{gt:thm:generic-transversality-in-families}, Pardon's proof goes through without changes. We sketch here how the different parts come together. Take a class in $H_c^*(\cZ/\CY_4)$, and represent it by a class $\gamma\in H_K^*(\cZ(X/B)/B)$ for some family of CY$4$-folds $X\to B$ over a finite simplicial complex $B$ and some compact $K\subseteq \cZ(X/B)$. By the graph trick \ref{lemma:graph-trick}, we may assume $K$ is subanalytic with $K\to B$ injective.

    Now enough divisors Proposition~\ref{gt:prop:enough-divisors} is applied to find disjoint relative divisors $D_i$ over open subsets $U_i$ in $B$, which jointly control all cycles in $K$. This requires trimming $X$ away from $K$, which leaves the image of $\gamma$ in the Pardon group unchanged. For the universal cycle $\cU(X/B)$, write $\im \cU(X/B)$ for its image inside $\cZ(X/B)\times_B X$ viewed as the reduced induced complex analytic subspace. Write $\cC \coloneqq \pi_X(\im \cU(X/B)\times_{\cZ(X/B)}K)$ and
    \begin{equation*}
        \cC_{L}\coloneqq \pi_X(\im \cU(X\times_B L/L)\times_{\cZ(X\times_B L/L)}(K\times_B L))
    \end{equation*}
    for any subset $L\subseteq B$. Consider now an open subset $U_i$. Since every $z\in K$ is $D_i$-controlled for some $i$, the intersection 
    \begin{equation*}
        A_i\coloneqq \cC_{U_i}\cap D_i
    \end{equation*}
    maps to $U_i$ with proper finite fibers. Generic transversality now deforms $X/B$ in neighborhoods of the $A_i$, excluding $D_i$. Since the $D_i$ are disjoint, the neighborhoods can be chosen disjoint as well. The regluing neighborhood is roughly pictured in Figure~\ref{gen:fig:pure-subdivision}.
    \begin{figure}[h]
        \centering
        \includegraphics[scale=0.8]{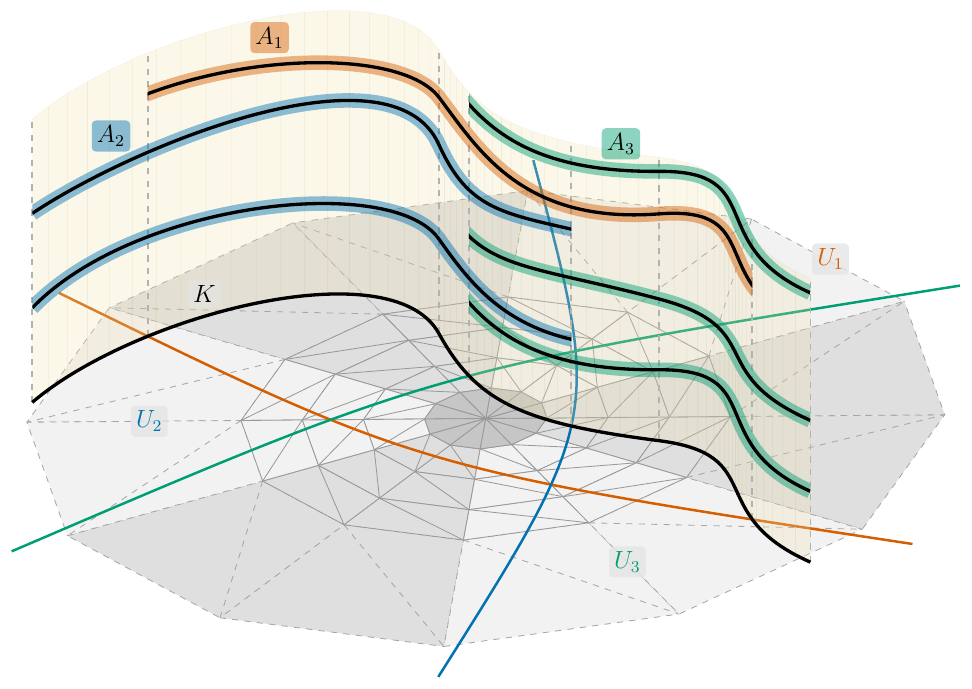}
        \caption{A subdivision of the base $B$ subordinate to a finite cover $\{U_1, U_2, U_3\}$ (in orange, blue, and green respectively), with shaded constructible opens (and their intersections). The base contains a compact set $K\subseteq \cZ(X/B)$ above which lie the intersections $A_i$ of the respective curves in $K$ with the divisors $D_i$ in every fiber, together with small regluing neighborhoods nearby.}
        \label{gen:fig:pure-subdivision}
    \end{figure}

    Now, instead of applying generic transversality directly, Pardon applies it to the stabilization $(X/B)\times \bR_{\geq 0}$ over the constructible open subsets $U_i\times \bR_{>0}$, in order to perform a retraction argument which allows him to identify the resulting classes in the Pardon group. Since the regularizing deformation does not change connectedness of the cycles under consideration, we obtain the statement for connected loci by the same proof.
\end{proof}

\subsection{Finite-Dimensional Regularization of Smooth Rationals}\label{gt:sec:fixing-subcomplex}

\subsubsection{}
Now that we have explained the generic transversality formalism, and used it to establish the isomorphism in Theorem~\ref{gt:thm:pardon-ring-sro-isom} above, we discuss two extensions of generic transversality used in the proofs of our generation statements in Section~\ref{rat:sec:rational-normal-data} and Section~\ref{inc:sec:proof} below. For both statements, we will focus on the adjustments necessary in the proof of Theorem~\ref{gt:thm:generic-transversality-in-families}.

For the smooth rational locus of Section~\ref{rat:sec:rational-normal-data}, we can recover a useful finite-dimensional formulation of generic transversality similar to \S\ref{gt:sec:single-curve-ks-strategy}. We now make this precise.

\subsubsection{}
\begin{definition}\label{rnd:def:rational-smooth-locus}
    For any CY$4$ family $X\to B$, we define the locally closed subset $R_m^\sm(X/B)\subset \cZ_\con(X/B)$ (and $R_m^{\sm,\sr\circ}(X/B)\subseteq \cZ_\con^{\sr\circ}(X/B)$) of cycles $z=m[C]$ where $C\subseteq X_b$ is a smooth embedded rational curve and the multiplicity is exactly \(m\). We call this entire collection, denoted by $R_m^\sm$ (resp. $R_m^{\sm,\sr\circ}$), the \textit{smooth rational locus}, and since this is a fiberwise local sublocus, we can construct the universal bivariant groups
    \begin{equation*}
        H_c^*(R_m^\sm/\CY_4),\quad\text{and}\quad H_c^*(R_m^{\sm,\sr\circ}/\CY_4)
    \end{equation*}
    by Proposition~\ref{rcl:prop:cohomology-fiberwise-cond} and \S\ref{gt:subsubsec:semi-regular-group}.
\end{definition}

\subsubsection{}
For a smooth embedded rational cycle $u:C\hookrightarrow X_b$, no normalization or special point data are required in the definition of semi-regularity and the obstructions are simply given by
\begin{equation*}
    \rOb(u)= H^1(C,N_{C/X_b}).
\end{equation*}
In light of the finite-dimensional formulation in \S\ref{gt:sec:single-curve-ks-strategy}, this motivates the following definition.
\begin{definition}\label{gt:def:v-regular}
    Let $V$ be a complex vector space, and let $X\to B\times B_0V$ be a family of CY$4$-folds, where $B_0V$ is a sufficiently small compact ball around $0$ in $V$ with a real analytic triangulation having $0$ as a vertex, such that each $X_b\coloneqq X|_{\{b\}\times B_0V}$ is the restriction of a holomorphic family over a neighbourhood of $B_0V$ in $V$, as for the families of Proposition~\ref{bg:prop:family-regluing}. Consider a cycle $z=m[C]\in R_m^\sm(X/B\times B_0V)$ over $(b,v)$ with $v$ in the interior of $B_0V$. Then the Kodaira--Spencer map of $C$ in $X_b\to B_0V$ is given by
    \begin{equation*}
        \kappa_z:V\to H^1(C,N_{C/X_{(b,v)}}).
    \end{equation*}
    We call $z$ \textit{$V$-regular} if the Kodaira--Spencer map above is surjective.
\end{definition}

\begin{remark}
    Equivalently, regarding $X\to B\times B_0V$ as a family over $B$ with fibres $X_b$, a cycle is $V$-regular if and only if it is regular in its fibre. In particular, semi-charts relative to $B$ at $V$-regular cycles are smooth, and we always form them over the simplices of $B$.

    $V$-regularity is not a stratum-wise condition as in Definition~\ref{gt:def:sr-complex}(c): $\kappa_z$ uses all directions of $V$, regardless of the smallest simplex of $B\times B_0V$ containing $(b,v)$, and in particular also over the vertex $0$. Hence, it does not depend on the triangulation of $B_0V$. We take $0$ to be a vertex so that $B\times\{0\}\subseteq B\times B_0V$ is a subcomplex, over which a deformation $\tilde X\to B\times B_0V$ of $X\to B$ restricts to $X$.
\end{remark}

\subsubsection{}
Now we address the main lemma, which allows us to recover the finite-dimensional formulation in \S\ref{gt:sec:single-curve-ks-strategy} for smooth embedded rational curves.

\begin{lemma}\label{gt:lemma:v-reg-open}
    The locus of $V$-regular cycles in $R_m^\sm(X/B\times B_0V)$ is open.
\end{lemma}
\begin{proof}
    Take a smooth embedded rational cycle $z=m[C]$ over $(b,v)$ and consider the exact sequence of normal bundles
    \begin{equation*}
        0\to N_{C/X_{(b,v)}}\to N_{C/X_b}\to V \otimes \cO_C \to 0.
    \end{equation*}
    Taking the associated long exact sequence in cohomology recovers the Kodaira--Spencer map
    \begin{equation*}
        \cdots \to V \xrightarrow{\kappa_z} H^1(C,N_{C/X_{(b,v)}}) \to H^1(C,N_{C/X_b}) \to 0.
    \end{equation*}
    This shows that $z$ is $V$-regular if and only if $H^1(C,N_{C/X_b})=0$. 
    
    Now consider a space $S$  of deformations of the associated stable map of $z$, which we will properly introduce in Definition~\ref{inc:def:semi-chart}. By \cite[Lemma 6.2]{pardon2023}, this is an open embedding in $R_m^\sm(X/B\times B_0V)$ around $z$. If $S$ is sufficiently small, all elements are smooth embedded rational curves in $X$, viewed as cycles with the same fixed multiplicity $m$. This makes a sufficiently small $S$ a Douady space of curves in $X$ relative to $B\times B_0V$, with a proper smooth universal family $\pi:\cC\to S$ with fibers $\bP^1$. Since maps from compact connected curves to a neighborhood of $B_0V$ in the affine $V$ must be constant, $S$ is equivalently just a Douady space of curves in $X$ relative to $B$. Since $\pi$ is smooth, $N_{\cC/X\times_B S}$ is locally free with fiber $N_{C'/X_{b'}}$ over $[C']$ with $C'$ over $(b',v')$. As the fibers of $\pi$ are $1$-dimensional, we have $R^2\pi_*=0$, so that
    \begin{equation*}
        R^1\pi_*\big(N_{\cC/X\times_B S}\big)
    \end{equation*}
    is coherent by the Direct Image Theorem \cite[Ch.10,\S4.6]{Grauert1984} and commutes with base change by cohomology and base change. In particular, the fibers
    \begin{equation*}
        R^1\pi_*\big(N_{\cC/X\times_B S}\big)|_{[C']}=H^1(C',N_{C'/X_{b'}})
    \end{equation*}
    control $V$-regularity. Precisely, the support of $R^1\pi_*\big(N_{\cC/X\times_B S}\big)$ consists of the non-$V$-regular points in $S$. By coherence, Nakayama's lemma shows that $V$-regularity is open.
\end{proof}

\begin{remark}
    This matches the behavior of semi-regularity in Pardon's work. Regularity itself is an open condition. However, since we forget multiplicity of cycles for its definition, the notion of semi-regularity is not necessarily open. Restricting to this smooth rational locus with fixed multiplicity allows us to recover openness.
\end{remark}

\subsubsection{}
\begin{proposition}\label{cor:maincor}
Let $X\to B$ be a family of CY $n$-folds over a finite simplicial complex with $n\geq 3$. Fix $K\subseteq R^\sm_m(X/B)$ a compact subanalytic subset for which the map $K\to B$ is injective. After possibly removing a closed subset from $X$ disjoint from $K$ and subdividing $B$, there exists a CY deformation 
\begin{equation*}
    \tilde{X}\to \tilde{B}\coloneqq B\times B_0V
\end{equation*}
of $X\to B$,  where $V$ is a finite dimensional complex vector space and $B_0V$ is as in Definition~\ref{gt:def:v-regular}, such that every $z$ in a neighborhood of $K$ inside $R^\sm_m(\tilde{X}/\tilde{B})$ is $V$-regular.
\end{proposition}

\begin{proof}
By the enough divisors Proposition~\ref{gt:prop:enough-divisors}, we have disjoint relative divisors $D_i$ over $U_i$ jointly controlling $K$. Taking $A_i$ to be the intersection of $D_i$ with the underlying curves of the cycles in $K$ puts us into the situation of Theorem~\ref{gt:thm:generic-transversality-in-families}. After trimming $X$ and subdividing $B$, its proof constructs the space of CY deformations $\cJ_D^\rCY(X/B)$ from the vector spaces $H^0_\rCY((\bD^\circ\setminus 0)\times \bD^{n-1}\times M_{ij}^\st,T_{X/B})$ associated to the charts \eqref{gt:eq:gt-charts}. For a finite-dimensional subspace $V$ of their direct sum, Proposition~\ref{bg:prop:family-regluing} and Corollary~\ref{bg:cor:regluing-cy-preservation} give a CY deformation $\tilde{X}\to \tilde{B}\coloneqq B\times B_0V$ of $X\to B$. By construction using regluing flows, its Kodaira--Spencer map $\kappa_z:V\to \rOb(u)=H^1(C,N_{C/X_{(b,v)}})$ at a cycle $z$ with underlying curve $u:C\hookrightarrow X_{(b,v)}$ is exactly
the surjection \eqref{gt:eq:composed-surjection-constrained-obstructions-gt-proof} restricted to $V$, whose target is $\rOb(u)=H^1(C,N_{C/X_{(b,v)}})$ here as there are no normalizations or special points. Since \eqref{gt:eq:composed-surjection-constrained-obstructions-gt-proof} is surjective, every $z\in K$ has a finite dimensional subspace $V_z$ for which $\kappa_z$ is surjective, making $z$ $V_z$-regular. By Lemma~\ref{gt:lemma:v-reg-open} this $V_z$ makes nearby points in $\tilde{R}_m^\sm(\tilde{X}/B\times B_0V_z)$ $V_z$-regular. As $K$ is compact, we can now take $V$ to be the sum of finitely many such $V_z$, and the resulting deformation makes every $z$ in a neighborhood of $K$ inside $R^\sm_m(\tilde{X}/\tilde{B})$ $V$-regular.
\end{proof}

\subsection{Branch Regularization}\label{rbr}

\subsubsection{}
On the locus of meetings of two rational curves discussed in Section~\ref{inc:sec:proof}, we require an extension of generic transversality, which allows us to simultaneously regularize both rational components as well as the ``common image problem'', which is the deformation problem of the full meeting curve subject to preserving its intersection point. We make this precise using the following definition.

\begin{definition}\label{gt:def:ordinary-rational-meeting-cycle}
    Consider a cycle $z\in \cZ(X/B)$ of the form 
    \begin{equation*}
        z=m_1C_1+m_2C_2
    \end{equation*}
    with two embedded smooth rational curves $u_i:C_i\hookrightarrow X_b$ and points $p_i\in C_i$ such that
    \begin{equation*}
        u_1(p_1)=u_2(p_2)=q.
    \end{equation*}
    If this is their only intersection points, and the two tangent lines at $q$ are distinct, then we call $z$ an \textit{ordinary rational meeting cycle}.
\end{definition}

\subsubsection{}
We want to regularize such cycles in a stricter way, which we define now.

\begin{definition}\label{gt:def:branch-semi-regular}
    Take an ordinary rational meeting cycle $z=m_1C_1+m_2C_2$. We say $z$ is \textit{branch semi-regular} if $z$ is semi-regular in the sense of Definition~\ref{gt:def:sr-complex}, and the deformation problems of the marked maps $(C_i,p_i,u_i)$ are unobstructed (in the family $X\to B$).
\end{definition}

\subsubsection{}

\begin{proposition}\label{gt:prop:rational-rational-gt}
    Let $X/B$ be a family of complex manifolds over a finite simplicial complex. Let $K\subseteq \cZ(X/B)$ be a finite subset of points, such that $K\to B$ is injective and each point in $K$ is an ordinary rational meeting cycle. After possibly removing a closed subset of $X$ disjoint from $K$ and subdividing $B$, there exists a CY deformation
    \begin{equation*}
        \tilde{X}\to \tilde{B}\coloneqq B\times S_0V
    \end{equation*}
    of $X\to B$, i.e. its pullback along $b\mapsto(b,0)$ is $X$, where $V$ is a finite dimensional complex vector space and $S_0V\subseteq B_0V$ is a $\dim_\bR(V)$-dimensional simplex with $0$ in its interior, over which we restrict the family of Proposition~\ref{bg:prop:family-regluing}, such that every $z\in K$, lying over $b\in B$, is branch semi-regular in the slice family $\tilde X|_{\{b\}\times S_0V^\circ}\to S_0V^\circ$.
\end{proposition}

\subsubsection{}
\begin{proof}
    Consider an ordinary rational meeting cycle $z=m_1C_1+m_2C_2$ with $u:C_1\cup C_2 \hookrightarrow X_b$. The components $u_i:C_i\to X_b$ of the associated stable map are the inclusions $j_i:C_i\hookrightarrow C$ composed with $u$. The adjunction $\id\mapsto (j_i)_*j_i^*$ then yields a surjection
    \begin{equation*}
        H^1(C,u^*T_{X_b})\twoheadrightarrow H^1(C_i,u_i^*T_{X_b})
    \end{equation*}
    using cohomological dimension $1$ of $C$. Since the deformations of the marked maps $(C_i,p_i,u_i)$ is controlled by the complex $[T_{C_i}(-p_i)\to u_i^*T_{X_b}]$, we get a surjection 
    \begin{equation*}
        H^1(C_i,u_i^*T_{X_b}) \twoheadrightarrow \rOb(C_i,p_i,u_i).
    \end{equation*}
    So, similar to \eqref{gt:eq:emb-curve-ob-surjects-onto-constrained-surjection}, we find that $H^1(C,u^*T_{X_b})$ surjects onto the obstructions of the associated stable map with fixed special point data as well as $\rOb(C_i,p_i,u_i)$.

    Pardon constructs the controlling divisors in \cite[Prop. 5.10]{pardon2023} as a union of transverse disks at a finite collection of smooth embedded points of the cycles, so that the controlling divisors never meet the intersection point $q$ of any ordinary rational meeting cycle $z$ in $K$. Hence, we obtain disjoint relative divisors $D_k$ over disjoint open subsets $U_k\subseteq B$, for any point $k\in K$, controlling the corresponding cycle $k$. We can then pick disjoint collections of disjoint open balls $U_{A_k}$ in $X$ around the intersection points $A_k$ of the cycle $k$ and $D_k$. Now we trim $X$ by removing $D_k\setminus U_{A_k}$ over $U_k^\Delta$. 
    
    If we write $U_k:C_k\to X_{b_k}$ for the associated embedded curve of every cycle $k\in K$, we get
    \begin{equation*}
        H^0_\rCY(U_{A_k}\setminus D_k,T_{X/B}) \twoheadrightarrow H^1(C_k,u_k^*T_{X_{b_k}})
    \end{equation*}
    by Lemma~\ref{gt:lemma:cy-preserving-surjection-only-overlap}. Hence, we can pick a finite dimensional $V_k$ in every $H^0_\rCY(U_{A_k}\setminus D_k,T_{X/B})$ surjecting onto the finite dimensional $H^1(C_k,u_k^*T_{X_{b_k}})$ and hence all three obstruction spaces as shown above. Since we chose the $U_k$ and $U_{A_k}$ disjoint, we can sum these vector spaces to obtain a finite dimensional subspace of $H^0_\rCY\big(\cup_k U_{A_k}\setminus D_k,T_{X/B} \big)$, which yields the desired deformation by Proposition~\ref{bg:prop:family-regluing} and Corollary~\ref{bg:cor:regluing-cy-preservation}.
\end{proof}

\section{Semi-Charts and Cycles with Connected Support}
\label{inc:sec:incidence}
After applying Theorem~\ref{gt:thm:generic-transversality-in-families},
we may work near a compact subset of the interior semi-regular locus.

We stratify this neighbourhood, and then estimate the dimension of each stratum.

For cycles with connected support, the estimate leaves only three
possibilities which may give a non-negative contribution: an embedded rational curve, an
embedded elliptic curve, or two embedded rational curves meeting at
one point.

\subsection{Semi-Charts}
\label{inc:sec:semi-charts}

\subsubsection{}
Let
\[
z=\sum_{i=1}^{\ell}m_i[C_i]\in\cZ(X/B)
\]
lie over $b\in B$, where the $C_i$ are distinct, irreducible, and
reduced. Write
\[
f_z\colon\widetilde C_z=\bigsqcup_i\widetilde C_i\longrightarrow X_b
\]
for the disjoint union of their normalization maps. We consider
deformations of $f_z$, allowing the complex structures of the domains
and the point of the base to vary, while keeping the multiplicities
$m_i$.

\begin{definition}
\label{inc:def:semi-chart}
The resulting deformation germ $S_z$, with its map
\begin{equation}
(S_z,z)\longrightarrow\bigl(\cZ(X/B),z\bigr),
\qquad
\bigsqcup_i f_i'\longmapsto\sum_i m_i(f_i')_*[\widetilde C_i'],
\label{inc:eq:semi-chart}
\end{equation}
is the \textit{semi-chart at $z$}.
\end{definition}

If $B$ is simplicial and $b\in\sigma^\circ$, we form the analytic
deformation space in a holomorphic extension of the family near
$\sigma$ and restrict it to the real open simplex $\sigma^\circ$.
Regularity is tested using the parameter directions in
$\sigma^\circ$.

\subsubsection{}
Recall the multiplicity stratification of
Subsection~\ref{subsec:multiplicity-filtration}, and put
\[
\widetilde\chi(z)=\sum_i\chi_{\mathrm{top}}(\widetilde C_i).
\]

\begin{lemma}[Semi-chart stratification
{\cite[Lemma~6.2]{pardon2023}}]
\label{inc:lem:semi-chart-stratification}
On each multiplicity stratum, $\widetilde\chi$ is upper
semicontinuous. The semi-chart at $z$ is an open embedding into the
locally closed locus $\cZ_{\bfm,\tilchi}$ on which
$\bfm=\bfm(z)$ and $\widetilde\chi=\widetilde\chi(z)$.
\end{lemma}

\begin{proof}
We recall Pardon's argument. A nearby cycle in the same multiplicity
stratum has a unique expression $\sum_i m_i[C_i']$ with each $C_i'$
close to $C_i$. We may therefore treat one reduced irreducible
curve $C$ at a time.

Choose small coordinate balls around the singular points of $C$.
Outside these balls, a nearby curve is diffeomorphic to $C$.
Inside each ball, the normalization of $C$ is a union of disks;
the nearby normalization has the same number of boundary circles
and no closed component. A compact oriented surface with $r$
boundary circles and no closed component has Euler characteristic
at most $r$, with equality precisely for $r$ disks. Thus the Euler
characteristic cannot increase. Equality holds precisely when the
nearby normalization map is a deformation of that of $C$.

For several components, Euler characteristics add, and equality
of the sum forces equality for each component. This proves both
assertions. The argument also applies when the target varies and,
over a simplicial base, on each real open simplex.
\end{proof}

Thus the ambient locus $\cZ_{\bfm,\tilchi}$ fixes only the
multiplicity multiset and the total normalization Euler
characteristic. Within each semi-chart, the topology of every
normalization component is fixed.

\subsection{Special Points}
\label{inc:sec:incidence-strata}

\subsubsection{}
Let $S$ be a sufficiently small representative of a semi-chart,
with universal family and map
\[
\pi\colon\widetilde{\mathcal C}
=\bigsqcup_i\widetilde{\mathcal C}_i\longrightarrow S,
\qquad
F\colon\widetilde{\mathcal C}\longrightarrow X_S.
\]
Write $f_s\colon\widetilde C_s\to X_{b(s)}$ for the map over
$s\in S$. As in Definition~\ref{gt:def:regular}, its special
points form the set
\begin{equation}
A_s=\bigl\{p\in\widetilde C_s:\
\rd f_s(p)=0\ \text{or}\
f_s(p)=f_s(q)\ \text{for some }q\ne p\bigr\}.
\label{inc:eq:special-points}
\end{equation}
We write $\cA(F)=\bigcup_{s\in S}A_s\subset\widetilde{\mathcal C}$.
After shrinking $S$, the component maps are simple and nonconstant,
and different components have distinct images. In particular, each
$A_s$ is finite.

\begin{definition}
\label{scs:def:discrete-data}
The \textit{discrete data} of $A_s$ record which normalization
component contains each point, which points have the same image,
and the vanishing orders
\[
e_p=\operatorname{ord}_p(\rd f_s),\qquad p\in A_s.
\]
An \textit{exact locus} is a locus on which these data, including
the full set of special points, remain constant up to relabelling.
\end{definition}

At a semi-regular cycle, the deformation problem preserving the
prescribed special-point constraints of
Definition~\ref{gt:def:regular} has a smooth local deformation
space. Within this space, the locus where the complete special-point
data remain unchanged maps to the exact locus in the semi-chart.
The marked points are allowed to move on their domains.

\subsubsection{}
We use subanalytic Whitney stratifications with real analytic
strata. Over a simplicial base, we require each stratum to lie over
the interior of a single simplex.

\begin{definition}
\label{inc:def:incidence-stratum}
Such a stratification of $U\subset\cZ(X/B)$ is an
\textit{incidence stratification} if, for each stratum $T$,
\begin{itemize}
\item the multiplicity and normalization Euler characteristic are
constant, and $T$ is locally contained in a semi-chart;
\item after passing to a finite cover of $T$, the normalization
components can be labelled and all special points are given by
disjoint analytic sections with constant discrete data.
\end{itemize}
The sections in the second condition are sections of the
special-point families $\cA(F)|_T\to T$ in the local semi-charts.
These families agree on overlaps by uniqueness of normalization.
\end{definition}

\begin{lemma}
\label{inc:lem:special-point-correspondence}
\label{inc:lem:finite-incidence-stratification}
For a sufficiently small semi-chart $S$, the map $\cA(F)\to S$
is finite. Every compact set $K\subset\cZ^{\sr\circ}(X/B)$ has an
open neighbourhood in the interior semi-regular locus admitting a
finite incidence stratification. For a simplicial base, the
stratification may be chosen compatibly on all faces.
\end{lemma}

\subsubsection{}
\begin{proof}
\textit{Step 1.}
We first show that the special points form a finite analytic
family. Work in the holomorphic extension of a semi-chart and set
\[
D=\widetilde{\mathcal C}\times_{X_S}\widetilde{\mathcal C}.
\]
Let $D^{\mathrm{res}}$ be the reduced union of the irreducible
components of $D$ not contained in the relative diagonal. Its
projection to $\widetilde{\mathcal C}$ is proper, so the image is
analytic by Remmert's theorem \cite[\S10.6.1]{Grauert1984}.
Away from the diagonal, $D^{\mathrm{res}}$ consists of pairs of
distinct points with the same image. If such pairs approach the
diagonal at $(p,p)$, then $p$ must be critical: near a point where
$\rd f_s(p)\ne0$, the map is an embedding and the fibre product
consists only of the diagonal. It follows that
\[
\cA(F)=\operatorname{pr}_1(D^{\mathrm{res}})
\cup\operatorname{Crit}(F/S).
\]
Thus $\cA(F)$ is a closed analytic subset of
$\widetilde{\mathcal C}$, proper over $S$. Since its fibres are
finite, the map $\cA(F)\to S$ is finite.

\medskip\noindent
\textit{Step 2.}
We next divide a neighbourhood of $K$ into loci over which the
normalizations form smooth families. It is enough to work in
finitely many analytic charts of the cycle space, using the
holomorphic extensions when the base is real. The multiplicity
and normalization Euler characteristic are analytically
constructible by \cite[\S6.1]{pardon2023}. The partition we
construct will refine their level sets.

Let $V$ be a closed irreducible analytic subset of one of these
charts. Take the reduced universal curve $\mathcal C_V$ over $V$
and normalize its total space:
\[
n\colon N\longrightarrow\mathcal C_V,
\qquad p\colon N\longrightarrow V.
\]
There is a dense open subset $V^\circ\subset V$ over which the
base is smooth, $p$ is smooth, and $n$ restricts to the
normalization of each fibre. To obtain it, remove the singular
locus of $V$ and the images of components of $N$ that do not
dominate $V$. The
singular locus of the remaining normal space has codimension
at least two, so its image is a proper analytic subset of $V$.

We also remove the critical values of $p$ and the points over
which the locus where $n$ is not an isomorphism contains a
complete component of the curve. These removals are analytic (because of the
properness and of the upper semicontinuity of the fiber dimension). The
remaining maps $n_s$ are normalizations: indeed, they are finite and birational, and have a smooth
source.

We will shrink $V^\circ$ once more to make the data at the
special points constant and then repeat the argument on its
closed complement.

\medskip\noindent
\textit{Step 3.}
The construction in Step~1 applies to the family over
$V^\circ$. To keep track of the orders of vanishing, write
$F\colon N\to X_V$ for the universal map and put
\[
\mathcal R=\operatorname{coker}
\bigl(F^*\Omega_{X_V/V}\longrightarrow\Omega_{N/V}\bigr).
\]
For $s\in V^\circ$, its restriction to $N_s$ is the cokernel
of $f_s^*\Omega_{X_{b(s)}}\to\Omega_{N_s}$, since relative
differentials and cokernels commute with base change. In a
local coordinate $t$ centred at $a\in N_s$, this cokernel is
isomorphic to $\mathcal O_{N_s,a}/(t^{e_a})$. Its length is
therefore exactly the order $e_a$ that we wish to record.

Since $p$ is proper, $p_*\mathcal R$ is coherent. We may shrink
$V^\circ$ so that this sheaf is locally free and the reduced
spaces parametrizing special points and pairs of special
points with the same image are finite unramified covers of
$V^\circ$. The fibre product $N\times_{X_V}N$ and the support
of $\mathcal R$ are defined over all of $V$ and proper over it.
Thus the sets removed can
be taken to be closed analytic subsets of $V$.

Locally label the special points by sections $a_1,\ldots,a_r$.
Restricting $\mathcal R$ to disjoint neighbourhoods of these
sections gives a decomposition $p_*\mathcal R=\bigoplus_j E_j$.
Each $E_j$ is locally free, being a direct summand of
$p_*\mathcal R$. Since the support is finite over $V^\circ$,
direct image commutes with taking fibres. Thus the fibre of
$E_j$ at $s$ is the vector space $(\mathcal R|_{N_s})_{a_j(s)}$,
of dimension $e_{a_j(s)}$. This order therefore equals the
locally constant rank of $E_j$.
The cover parametrizing pairs likewise makes the partition
by common image locally constant. The normalization components
and their genera are locally constant because $p$ is smooth
and proper.

We have now made all the discrete data constant, up to
relabelling, on $V^\circ$. Repeat the construction on its
closed complement. The dimension decreases, so induction
gives a locally finite partition into analytically
constructible subsets with constant discrete data.

\medskip\noindent
\textit{Step 4.}
Restrict these subsets and maps to the real base. The Whitney
stratification theorem for subanalytic sets and proper maps
\cite[Part~I, \S1.7]{Goresky1988} gives a refinement by Whitney
strata, together with compatible stratifications of the finite
maps parametrizing special points. For a simplicial base,
include the inverse images of all faces among the subsets to
be stratified. The restrictions to faces and their
intersections then agree. By local finiteness and compactness,
we may choose a neighbourhood of $K$ meeting only finitely
many strata.

Over each stratum $T$, the special points form a finite
covering with locally analytic sections. Passing to the
finite cover of $T$ parametrizing their labellings makes
these sections global. We may also label the normalization
components after a finite cover. The discrete data remain
constant under these labels, so this is an incidence
stratification.
\end{proof}

\subsection{Dimension Estimate}
\label{inc:sec:connected-support-index}

\begin{lemma}
\label{inc:lem:connected-support-index}
Let $T$ be an incidence stratum in the semi-regular locus over
$\sigma^\circ$, and let
\[
f\colon\bigsqcup_{i=1}^{\ell}\widetilde C_i\longrightarrow X_b
\]
represent a point of $T$. Write $g_i=g(\widetilde C_i)$,
$d_\sigma=\dim_{\mathbb R}\sigma$, and let $A$ be the set of
special points. The problem of deforming $f$ while preserving
the prescribed conditions at these points has relative complex
index
\begin{equation}
\nu(f)=\sum_{i=1}^{\ell}(1-g_i)+|A|
-4\left(|A|-|f(A)|+\sum_{p\in A}e_p\right).
\label{inc:eq:connected-support-index}
\end{equation}
If the reduced support is nonempty and connected, then
\begin{equation}
\nu(f)\leq2-\ell-\sum_i g_i\leq1.
\label{inc:eq:connected-support-bound}
\end{equation}
In this case, $\nu(f)=1$ precisely when the support is a smooth
embedded rational curve, and $\nu(f)=0$ precisely when it is a
smooth embedded elliptic curve or the union of two smooth embedded
rational curves meeting at exactly one point. The tangent lines
at their meeting point are allowed to coincide.

In all cases, the dimension of $T$ satisfies
\begin{equation}
\dim_{\mathbb R}T\leq d_\sigma+2\nu(f).
\label{inc:eq:stratum-dimension-bound}
\end{equation}
\end{lemma}

\begin{proof}
\textit{Step 1.}
We first compute the index. Without marked points, Riemann--Roch
gives $\sum_i(1-g_i)$, allowing the complex structures of the
domains to vary. Allowing the points of $A$ to move adds $|A|$
parameters. For $q\in f(A)$, write $A_q=A\cap f^{-1}(q)$ and
$t_q=|A_q|$. Requiring the points of $A_q$ to have a common image
gives $4(t_q-1)$ equations. Requiring $\rd f$ to vanish to order
at least $e_p$ at $p$ gives another $4e_p$ equations. Summing
these contributions proves \eqref{inc:eq:connected-support-index},
compare \cite[\S5.1, eq.~(5.1)]{pardon2023}.

By semi-regularity, these deformations form a smooth local space
of real dimension $d_\sigma+2\nu(f)$. Labelling the special points
embeds $T$ locally in this space: forgetting the marked points
recovers its inclusion in the semi-chart. This proves
\eqref{inc:eq:stratum-dimension-bound}.

\medskip\noindent
\textit{Step 2.}
Now assume that the support is nonempty and connected.
If $t_q\geq2$, the contribution over $q$ to the index is
\[
t_q-4(t_q-1)-4\sum_{p\in A_q}e_p
\leq-2(t_q-1).
\]
If $t_q=1$, the unique point is critical, so its
contribution is at most $1-4=-3$.

Form the graph whose vertices are the normalization components
and the points of $f(A)$. Each $p\in A$ gives an edge joining
its component to $f(p)$. This graph is connected, so its
$|A|$ edges and $\ell+|f(A)|$ vertices
therefore satisfy
\[
\sum_{q\in f(A)}(t_q-1)=|A|-|f(A)|\geq\ell-1.
\]
Combining this with the preceding estimates gives
\[
\nu(f)\leq\ell-\sum_i g_i-2(\ell-1)
=2-\ell-\sum_i g_i\leq1.
\]

\medskip\noindent
\textit{Step 3.}
It remains to determine when $\nu(f)\geq0$.
The bound gives $\ell\leq2$.

If $\ell=1$, any special points would lower the index by at
least two. Hence $A$ is empty and $g_1$ is zero or one.
The normalization map is injective and immersive, hence an
embedding by properness. Its index is one in genus zero and
zero in genus one.

If $\ell=2$, both genera are zero and $\nu(f)=0$. The contribution
over each $q\in f(A)$ is at most $-2$, with equality only for
two unramified preimages. Thus there is exactly one
such image, with one preimage on each component, and no other
special points. Both components are therefore embedded rational
curves meeting at one point. Conversely, this configuration has
index zero regardless of whether the tangent lines coincide.
\end{proof}

\subsection{The Cohomological Range}
\label{inc:sec:cohomological-range}

Let \(\cZ_{\mathrm{con}}^{\mathrm{sr}\circ}(X/B)\) be the locus of
cycles with connected reduced support in the interior of the
semi-regular locus.  We use the same notation after passing to the
directed colimit of Definition~\ref{coh:def:grothendieck-group-1-cycles}.

\begin{proposition}
\label{inc:prop:cohomological-range}
With rational coefficients,
\begin{equation}
 H_c^j\left(\cZ_{\mathrm{con}}^{\mathrm{sr}\circ}/\CY_4\right)=0
 \qquad\text{for }j>2.
 \label{inc:eq:cohomological-range}
\end{equation}
More precisely, after filtering a finite Whitney stratification by
closed unions, its associated graded in nonnegative bivariant degree
can receive contributions only as follows:
\begin{equation}
\begin{array}{c|c|c}
 \text{degree}&\text{locus}&\text{representative on a stratum}\\ \hline
 2&\text{smooth embedded rational curves}&\text{point}\\
 1&\text{smooth embedded rational curves}&\text{oriented loop}\\
 0&\text{smooth embedded rational curves}&\text{oriented surface}\\
 0&\text{smooth embedded elliptic curves}&\text{point}\\
 0&\text{two embedded rational curves meeting at one point}&\text{point}.
\end{array}
\label{inc:eq:associated-graded-generators}
\end{equation}
The last row includes the case in which the two tangent lines agree.  The
Poincar\'e duals and their ambient lifts will be constructed in the
following sections.
\end{proposition}

\begin{proof}
Represent a compactly supported class over a family \(X\to B\), and
choose the finite Whitney stratification of
Lemma~\ref{inc:lem:finite-incidence-stratification} near its support.
If \(T\) is a stratum over \(\sigma^\circ\), let \(\nu(T)\) denote the
constant value of \(\nu(f)\) on \(T\).  The dualizing complex of
the oriented real manifold \(\sigma^\circ\) gives, up to its orientation
local system,
\(H_c^j(T/\sigma^\circ)=H_c^{j+d_\sigma}(T)\).  The dimension estimate
in Lemma~\ref{inc:lem:connected-support-index} therefore implies that
this group vanishes for \(j>2\nu(T)\), and hence for \(j>2\).
The frontier axiom of the stratification allows us to construct a filtration by closed unions of strata. Using the resulting
localization sequences proves the same vanishing near the entire
compact support. Passage to the directed colimit gives
\eqref{inc:eq:cohomological-range}; this is the connected analogue of
the argument in \cite[Lemma~6.3 and Propositions~6.4 and~6.9]{pardon2023}.

A stratum with \(\nu(T)<0\) contributes nothing in nonnegative degree.
If \(\nu(T)=0\), only degree zero can occur, and if \(\nu(T)=1\), only
degrees zero, one, and two can occur.  The last assertion of
Lemma~\ref{inc:lem:connected-support-index} gives the loci in
\eqref{inc:eq:associated-graded-generators}.  Strata of smaller dimension
give no additional possibilities.  On a stratum of maximal dimension,
Poincar\'e duality with its orientation local system identifies these
compactly supported classes with homology in dimensions zero, one, and
two. 
\end{proof}

\subsection{The Elliptic Locus}\label{sec:elliptic-locus}
We first consider smooth elliptic curves, whose deformation problem
has index zero. Smooth rational curves and pairs of rational curves
meeting at one point are studied in
\S\ref{rat:sec:rational-normal-data} and \S\ref{inc:sec:proof}, respectively.

\subsubsection{}
\begin{definition}\label{def:elliptic-smooth-locus}
    For any family $X\to B$, we define the locally closed subset $E_m^\sm(X/B)\subset \cZ_\con^{\sr\circ}(X/B)$ of smooth embedded elliptic curves of multiplicity $m$. We call this entire collection, denoted by $E_m^\sm$, the \textit{smooth elliptic locus}, and since this is a fiberwise local sublocus, we can construct the universal bivariant group
    \begin{equation*}
        H_c^*(E_m^\sm/\CY_4)
    \end{equation*}
    by Proposition~\ref{rcl:prop:cohomology-fiberwise-cond}.
\end{definition}

\subsubsection{}\label{sec:elliptic-point-thom-classes}
Let $X\to B$ be a family, and consider a cycle $z\in E_m^\sm(X/B)$ in the smooth elliptic locus. By semi-regularity, it is contained in a smooth semi-chart $S$, open in the smooth elliptic locus. By the index formula \eqref{inc:eq:connected-support-index}, the semi-chart is $0$-dimensional over the base. Hence, the Poincar\'e dual of the class $[z]$ defines
\begin{equation*}
    \tau_z \in H_c^0(S/B),
\end{equation*}
which we call the \textit{Thom representative} of the point $z$. Via extension by zero, this may also be seen as a class in $H_c^0(E_m^\sm/\CY_4)$. 

\subsubsection{}\label{gt:subsec:elliptic-regularization}
\begin{lemma}\label{gt:lem:elliptic-point-comparison}
    Let $z=m[E]\in E_m^\sm(X/B)$ be a cycle in the smooth elliptic locus of a family $X\to B$, and let
$\tau_z$ be its Thom representative. Choose line bundles
$L_i\in\operatorname{Pic}^0(E)$ of infinite order, with
$L_1\otimes L_2\otimes L_3\simeq K_E$, and put
$N^\circ=L_1\oplus L_2\oplus L_3$. Then $m[E]$ is an isolated
semi-regular cycle in $\operatorname{Tot}_E(N^\circ)$, and $\tau_z$ agrees with its
 Thom representative in $H_c^0(E_m^\sm/\CY_4)$.
\end{lemma}

\begin{proof}
\textit{Step 1: The split model.}
Each $L_i$ is nontrivial of degree zero, so
$H^0(E,N^\circ)=H^1(E,N^\circ)=0$.
Let $C$ be a compact irreducible curve in
$\operatorname{Tot}_E(N^\circ)$. It cannot lie in an affine
fibre, so its normalization projects
finitely onto $E$, say by $f:\widetilde C\to E$ of degree $d$.
The fibre coordinates are sections of $f^*L_i$.
These degree $0$ bundles are nontrivial: otherwise the norm
map would give $L_i^d\simeq\mathcal O_E$.

Thus all coordinates under consideration vanish, or, equivalently, $C=E$. Therefore, every compact effective cycle is a multiple of $E$.
Moreover, clearly, $m[E]$ is isolated and lies in $E_m^\sm$.

\medskip\noindent
\textit{Step 2: A family joining the two curves.}
Put $N=N_{E/X_b}$, with $\det N\simeq K_E$ induced by the
volume form. Work simplex-wise in holomorphic extensions of
the base. If $V$ is the tangent space to this extension and
$\mathcal N$ the normal bundle of $E$ in the total space, then
\[
0\longrightarrow N\longrightarrow\mathcal N
\longrightarrow\mathcal O_E\otimes V\longrightarrow0.
\]
Deforming $E$ and the base point $b$ to their normal cones
joins the given family to
$\operatorname{Tot}_E(\mathcal N)\to V$.
Multiplying a cocycle for this extension by $1-s$ joins this to
$\operatorname{Tot}_E(N)\times V$.

The smooth bundles underlying $N$ and $N^\circ$ are isomorphic
compatibly with their determinants. On this bundle, set
\[
\bar\partial_t=(1-t)\bar\partial_N+t\bar\partial_{N^\circ}.
\]
These operators are integrable because $E$ is a curve, their total spaces, with the
product factor $V$, complete the
desired family of Calabi--Yau fourfolds with embedded elliptic
curves. 

\medskip\noindent
\textit{Step 3: Comparison of the point classes.}
Apply the infinitesimal deformation construction of
\S\ref{gt:sec:inf-def} along this compact family.
Openness of surjectivity and compactness result in finitely many
additional parameters which, together with the existing base
directions, make the linearized deformation operator surjective. 

Choose the subdivision transverse to these smooth deformation
spaces, relative to the endpoint faces. Apply the relative
Sard--Smale argument of \cite[Lemma~5.9]{pardon2023} to the
Calabi--Yau deformations of \S\ref{gt:sec:inf-def}, fixing
the endpoint families.
Choose the perturbation sufficiently small. Continuous right
inverses along the compact interval and the parametric
implicit function theorem give a nearby path of embedded
elliptic curves, with the same endpoints. Choose these paths
compatibly on common faces, using transversality there.
Generic transversality
makes an open neighbourhood of their multiples $m[E_t]$
in the full cycle space
semi-regular. Thus this path $\gamma$ lies in $E_m^\sm$.
The relation $\partial[\gamma]=[z_1]-[z_0]$ finishes the proof.
\end{proof}

\subsubsection{}
The same proof allows $E$ to vary in a family of elliptic curves.
Since their moduli space is connected, the resulting class is
also independent of $E$.

\begin{definition}\label{def:elliptic-standard-class-Em}
    For any $N^\circ$ as in Lemma~\ref{gt:lem:elliptic-point-comparison} above, we denote the point Thom representative of $m[E]$ in $\operatorname{Tot}_E(N^\circ)$ by
    \begin{equation*}
        e_m\coloneqq \tau_{m[E]}\in H_c^0(E_m^\sm/\CY_4).
    \end{equation*}
    By Lemma~\ref{gt:lem:elliptic-point-comparison}, this is independent of $N^\circ$ and hence well-defined.
\end{definition}

\begin{corollary}
    The group $H_c^0(E_m^\sm/\CY_4)$ is generated by $e_m$.
\end{corollary}

\section{Rational Normal Data and Poincar\'e-Dual Representatives}
\label{rat:sec:rational-normal-data}
The chart \eqref{inc:eq:associated-graded-generators} gives us explicit types of strata to study. In this section we study the locus $R_m^\sm$ of smooth embedded rational curves, which is precisely defined in Definition~\ref{rnd:def:rational-smooth-locus}. Our goal is to prove that all surfaces inside smooth rational locus can be reduced to a single (and very concrete) family over $\mathbb P^1$ via deformations. 

The main technique to do so is analytification of rational normal data. This means we take a purely topological version of a family of rational curves with rank $3$ normal bundle, see Definition~\ref{rat:def:rational-normal-datum}, and analytify it in \S\ref{rat:sec:relative-analytic-realization} using tools from \S\ref{rat:sec:analytification-tools}. Composed with the generic transversality results in \S\ref{gt:sec:fixing-subcomplex}, this allows us, in the desired degrees, to define a surjection onto $H_c^*(R_m^\sm/\CY_4;\bQ)$ from a certain bordism group of rational normal data, see \S\ref{subsec:rational-normal-data-map}. The latter group is computable, see \S\ref{rat:sec:degree-two-normal-generator} giving us explicit generators.

\subsection{Analytification Tools}\label{rat:sec:analytification-tools}
First, we collect some tools used in the proof of our relative analytic realization theorem.
\subsubsection{Approximation Relative to Faces}
\label{rat:sec:exact-face-approximation}

The first tool in the proof of relative analytic realization is an
approximation result which leaves fixed all previously constructed faces.
Here, a smooth map on a product with a simplex is understood to be smooth
up to the boundary, and a real-analytic map on a face $F$ is understood to
extend real-analytically to a neighborhood of that face in its affine
span, which we write as $\Aff(F)$.

The material here may be standard, but we record it since we do not know a precise reference.

\begin{lemma}
\label{rat:lem:exact-face-approximation}
Let $Z$ be a point or $S^2$, and let $\Delta^k$ have barycentric
coordinates
\(
    \lambda_0,\ldots,\lambda_k.
\)
Suppose that
\(
    f:Z\times\Delta^k\longrightarrow\mathbb R^q
\)
is smooth and that the restriction of $f$ to every proper face of
$\Delta^k$ is real-analytic. Then, for every finite $r$ and every
$\epsilon>0$, there is a real-analytic map $g$, defined on a neighborhood
of $Z\times\Delta^k$ in $Z\times \bR^k$, such that $g=f$  on $Z\times\partial\Delta^k$, and
\(
    \lVert g-f\rVert_{C^r}<\epsilon.
\)

More generally, the values on the faces may be replaced by any
compatible values sufficiently close to those of $f$ in a
sufficiently high finite $C^s$-norm. The resulting $g$ may still be chosen
arbitrarily close to $f$ in $C^r$.
\end{lemma}

\begin{proof}
Write
\(
    F_i=\{\lambda_i=0\}\subseteq\Delta^k,
\)
\(
    b_i=\prod_{j<i}\lambda_j,
\)
where $b_0=1$. We first construct a real-analytic map $H$ whose restriction
to the boundary agrees with that of $f$.

Order the faces as $F_0,\ldots,F_k$. Suppose inductively that an analytic
map $H_{i-1}$ has been constructed which agrees with $f$ on
$F_0,\ldots,F_{i-1}$. On $F_i$, set
\[
    d_i=f|_{Z\times F_i}-H_{i-1}|_{Z\times F_i}.
\]
For every $j<i$, the function $d_i$ vanishes on
\(
    Z\times(F_i\cap F_j).
\)
Successive applications of the  Hadamard lemma, inside the affine
span of $F_i$, therefore give
\(
    d_i=b_iq_i
\)
for a real-analytic map $q_i$ defined on a neighborhood of
$Z\times F_i$ in $Z\times\operatorname{Aff}(F_i)$.

Choose an affine retraction from a neighborhood of $\Delta^k$ to a
neighborhood of $F_i$ which restricts to the identity on $F_i$. It gives
an analytic extension $\widetilde q_i$ of $q_i$. Define
\[
    H_i=H_{i-1}+b_i\widetilde q_i.
\]
We corrected the restriction to $F_i$, while none of the already corrected faces $F_j$, $j<i$ was changed.
Continuing in this fashion, we get a real-analytic map $H=H_k$ satisfying $H=f$ on $Z\times\partial\Delta^k$.

Note that the smooth map $f-H$ vanishes on every face. Thus, applying successively
the smooth Hadamard lemma, we get a smooth map $s$ such that
\(
    f-H=(\lambda_0\cdots\lambda_k)s.
\)
Set
\[
    b=\lambda_0\cdots\lambda_k.
\]

Now, if $Z=S^2$, embed $S^2$ as the unit sphere in $\mathbb R^3$ and
extend $s$ smoothly to a neighborhood of
\(
    S^2\times\Delta^k\subseteq\mathbb R^3\times\mathbb R^k.
\)
After multiplying by a cutoff outside a slightly larger neighborhood,
standard simultaneous polynomial approximation gives a polynomial map
$p$ which is arbitrarily close to $s$ in $C^r$ on
$S^2\times\Delta^k$. 

When $Z$ is a point, we may apply the same argument
only in the simplex variables.

Now define
\(
    g=H+bp.
\)
The map $g$ is real-analytic on a neighborhood of
$Z\times\Delta^k$. Since $b$ vanishes on every proper face,
\[
    g|_{Z\times\partial\Delta^k}
    =
    H|_{Z\times\partial\Delta^k}
    =
    f|_{Z\times\partial\Delta^k}.
\]
Moreover,
\(
    g-f=b(p-s),
\)
so choosing $p$ sufficiently close to $s$ gives
\(
    \lVert g-f\rVert_{C^r}<\epsilon.
\)

The rest is clear.
\end{proof}

\subsubsection{Analytic Stable Presentations}
\label{rat:sec:analytic-stable-presentations}

\begin{definition}
    A \textit{stable analytic presentation} of a compatible simplexwise
real-analytic complex vector bundle $E\to Y$ is a compatible
simplexwise real analytic complex-linear injection
\[
    \iota:E\hookrightarrow Y\times\mathbb C^M
\]
for some finite $M$.
\end{definition}

\begin{lemma}
\label{rat:lem:analytic-stable-presentation}
Let $K$ be a finite simplicial complex. Let $Y$ be either
\(
    Y=|K|
\)
or the unit sphere bundle
\(
    Y=S(V)\longrightarrow |K|
\)
of a compatible simplexwise real-analytic oriented Euclidean vector
bundle $V$ of rank three. For a subcomplex $L\subseteq K$, write
\(
    Y_L=Y\times_{|K|}|L|.
\)

Now, let $E\to Y$ be a compatible simplexwise real-analytic complex vector
bundle. Then $E$ admits a stable analytic presentation.

More generally, suppose that a stable analytic presentation
\[
    \iota_L:E|_{Y_L}\hookrightarrow Y_L\times\mathbb C^{M_0}
\]
has already been chosen over a subcomplex $L\subseteq K$. Then there are
$M\geq M_0$ and a stable analytic presentation
\[
    \iota:E\hookrightarrow Y\times\mathbb C^M
\]
whose restriction to $Y_L$ is
\[
    \iota_L\oplus 0:
    E|_{Y_L}\hookrightarrow
    Y_L\times\bigl(\mathbb C^{M_0}\oplus\mathbb C^{M-M_0}\bigr).
\]
The analogous statement holds for real vector bundles.
\end{lemma}

\begin{proof}
We start with the
presentation prescribed on $L$, and argue inductively over skeleta of $K$.

We suppose, hence, that a presentation has already
been constructed over the union of $L$ and of the simplices of dimension
less than $k$, and let $\sigma$ be a $k$-simplex not contained in $L$.

The presentation is therefore already defined over
\(
    Y_{\partial\sigma}
    =
    Y\times_{|K|}\partial\sigma.
\)

In the sphere bundle case, over $\sigma$, the oriented Euclidean bundle $V$ is analytically trivial.
Indeed, it is smoothly trivial because $\sigma$ is contractible, a smooth
orthonormal frame may be approximated by an analytic frame, and analytic
Gram--Schmidt process makes it orthonormal. Thus, in this case, we
may identify
\(
    Y_\sigma\cong S^2\times\sigma
\)
real-analytically.

(In the $Y=\abs{K}$ case, we, of course, have $Y_\sigma=\sigma$ instead.)

For the sphere factor, we will use the affine complexification
\[
    Q=\left\{
        (z_1,z_2,z_3)\in\mathbb C^3:
        z_1^2+z_2^2+z_3^2=1
      \right\}.
\]
Its real locus is $S^2$, and it admits arbitrarily small Stein
neighborhoods in $Q$. Let $T$ denote such a Stein neighborhood, and choose also a bounded
convex Stein neighborhood $D$ of $\sigma$ in its affine
complexification.

After shrinking $T$ and $D$, the real-analytic bundle
$E|_{Y_\sigma}$ extends to a holomorphic vector bundle
$\mathcal E$ on $T\times D$. The already constructed presentation on
the faces extends holomorphically to the closed analytic subset
\[
    A=
    \bigcup_{i=0}^k
    T\times\bigl(D\cap\{\lambda_i=0\}\bigr)
    \subseteq T\times D.
\]

Suppose the presentation on $Y_{\partial\sigma}$ takes values in
$\mathbb C^m$. Its coordinate functions give sections
\(
    s_1,\ldots,s_m\in
    H^0\bigl(A,\mathcal E^*|_A\bigr)
\)
which span every fiber of $\mathcal E^*|_A$. Since $T\times D$ is Stein,
Cartan's Theorem B shows that $H^0(T\times D,\mathcal E^*)$ surjects onto $H^0\bigl(A,\mathcal E^*|_A\bigr)$, and hence, these section extend to holomorphic sections,
still denoted
\(
    s_1,\ldots,s_m\in H^0(T\times D,\mathcal E^*)
\)
by abuse of notation.

By Cartan's Theorem A and compactness of $Y_\sigma$, there are finitely
many further holomorphic sections
\[
    t_1,\ldots,t_\ell\in H^0(T\times D,\mathcal E^*)
\]
which span $\mathcal E^*$ on a neighborhood of $Y_\sigma$. Set
\(
    b_\sigma=\lambda_0\cdots\lambda_k.
\)
We obtain a holomorphic bundle map defined by the sections
\[
    (s_1,\ldots,s_m,\,
    b_\sigma t_1,\ldots,b_\sigma t_\ell):\mathcal E\longrightarrow
    (T\times D)\times\mathbb C^{m+\ell}.
\]

We have therefore extended the presentation across $\sigma$ without
changing it on any face.
Proceeding over all simplices and enlarging the common target whenever
necessary gives a finite $M$, since $K$ is finite. This proves both the
absolute and relative complex statements.

For a real vector bundle, apply the complex statement to its
complexification and then restrict to the real bundle, regarding
$\mathbb C^M$ as $\mathbb R^{2M}$.
\end{proof}

\subsubsection{Analytification Relative to Faces}
\label{rat:sec:analyticization-faces}

\begin{definition}
    A \textit{concordance relative to} $Y_L$ means a compatible simplexwise smooth
    family over $Y\times[0,1]$ which is constant over
    \(
        Y_L\times[0,1].
    \)
\end{definition}

\begin{lemma}
\label{rat:lem:analyticization-faces}
Let $K$ be a finite simplicial complex and $L\subseteq K$ a subcomplex.
Let the space $Y\to |K|$ be either $|K|$ itself or the
unit sphere bundle
\(
    Y=S(V)\longrightarrow |K|
\)
of a compatible simplexwise real-analytic oriented Euclidean bundle
$V$ of rank three. 

Let $E\to Y$ be a complex vector bundle which is continuous globally and
smooth over every simplex. Suppose that $E|_{Y_L}$ is equipped with a
compatible simplexwise real-analytic structure.

After a finite subdivision compatible with $L$, the bundle $E$ is
concordant relative to $Y_L$ to a compatible simplexwise real-analytic
bundle $E_{\mathrm{an}}\to Y$.

The same assertion holds simultaneously for a finite collection of real
or complex vector bundles and for the following data between them:
\[
    \text{sections, bundle maps, bundle isomorphisms, metrics,
    and connections}.
\]

Furthermore, the concordance may be chosen to preserve any condition
which is open in the corresponding finite-dimensional bundle, including
nonvanishing, injectivity, surjectivity, invertibility, positivity of a
metric, and preservation of orientation.
\end{lemma}

\begin{proof}
We first treat a single complex vector bundle. The same construction applies simultaneously to finitely many bundles. Carrying the additional data through this construction, and is left to the reader as an exercise.

By Lemma~\ref{rat:lem:analytic-stable-presentation}, there is an analytic
injection
\[
    \iota_L:E|_{Y_L}\hookrightarrow
    Y_L\times\mathbb C^{M_0}.
\]
Using \cite[Lemma 10.12]{Lee2012} simplex-wise, extend $\iota_L$ to a compatible simplexwise smooth bundle map
\(
    \widetilde\iota:E\longrightarrow Y\times\mathbb C^{M_0}.
\)

The same proof as \cite[Prop. 1.4]{Hatcher-VBKT} allows us to choose a global compatible simplexwise smooth injection
\(
    \iota_0:E\hookrightarrow Y\times\mathbb C^{M_1}.
\)
After a subdivision compatible with $L$, choose a nonnegative compatible
simplex-wise smooth function
\(
    \rho:|K|\longrightarrow[0,\infty)
\)
whose zero locus is precisely $|L|$. Pulling $\rho$ back to $Y$, define
\[
    \iota=
    \bigl(\widetilde\iota,\rho\iota_0\bigr):
    E\longrightarrow
    Y\times\bigl(\mathbb C^{M_0}\oplus\mathbb C^{M_1}\bigr).
\]
Over $Y_L$, this is $\iota_L\oplus0$. Away from $Y_L$, its second
component is a nonzero multiple of the injection $\iota_0$. It is
therefore injective everywhere.

Since vector bundles are cumbersome to deal with, we now shift perspective to instead consider certain operator-valued functions. 

Set $M=M_0+M_1$ and identify $E$ with the image of $\iota$ in the trivial
Hermitian bundle $Y\times\mathbb C^M$. The bundle of Hermitian operators on $Y\times\mathbb C^M$ is the trivial bundle $Y\times\Herm_M$, where $\Herm_M$ is the vector space of Hermitian $M\times M$-matrices. Let
\[
    P_E:Y\longrightarrow\operatorname{Herm}_M,\ y\mapsto \proj_{E_y}^\perp
\]
be the section defined by the orthogonal projection onto the image of $E$ in $\bC^M$. Then
\(
    P_E^2=P_E=P_E^*,
\)
viewing the images as operators. By the properties of $\iota$, the map $P_E$ is compatible and smooth over every simplex, and it is real-analytic over $Y_L$.

We now make $P_E$ analytic one simplex at a time. Suppose the construction
has already been made over $L$ and over the simplices of dimension less
than $k$, and let $\sigma$ be a $k$-simplex not contained in $L$. As in \S\ref{rat:sec:analytic-stable-presentations}, an
analytic orthonormal trivialization of $V|_\sigma$ identifies
$Y_\sigma$ with either
\(
    \sigma
\) or \(
    S^2\times\sigma.
\)
Apply Lemma~\ref{rat:lem:exact-face-approximation} to the matrix entries
of $P_E|_{Y_\sigma}$, prescribing on the boundary the analytic values
already constructed there. Taking the approximation in the real vector
space of Hermitian matrices gives a compatible simplexwise real-analytic
Hermitian matrix
\(
    A_\sigma:Y_\sigma\longrightarrow\operatorname{Herm}_M
\)
which agrees exactly with the previously constructed projection on
$Y_{\partial\sigma}$ and is arbitrarily close to $P_E$.

For a Hermitian matrix sufficiently close to an orthogonal projection,
define
\[
    \Pi(A)=
    \frac{1}{2\pi i}
    \int_{|\zeta-1|=1/2}
       (\zeta I-A)^{-1}\,d\zeta.
\]
The contour separates the part of the spectrum near $1$ from the part
near $0$. Consequently, $\Pi(A)$ is an orthogonal projection of the same
rank. The map $A\mapsto\Pi(A)$ is real-analytic on a neighborhood
of any compact set of projections, and fixes all elements of such a family pointwise.

We can thus replace $A_\sigma$ by $\Pi(A_\sigma)$. Since $\Pi$ fixes projections, this new
projection agrees with the already constructed projection on
every face. Proceeding over the skeleta, we get thus a compatible simplexwise
real-analytic projection
\(P_{E,\mathrm{an}}:Y\longrightarrow\operatorname{Herm}_M
\)
such that
\(
    P_{E,\mathrm{an}}=P_E
\)
on $Y_L$ and $P_{E,\mathrm{an}}$ is uniformly as close to $P_E$ as desired.

Choose the approximation sufficiently close that the same spectral
contour is valid for
\(
    (1-t)P_E+tP_{E,\mathrm{an}}
\)
for every $t\in[0,1]$. Then
\(
    P_{E,t}=
    \Pi\bigl((1-t)P_E+tP_{E,\mathrm{an}}\bigr)
\)
is a compatible simplexwise smooth path of projections. It satisfies
\(
    P_{E,0}=P_E,P_{E,1}=P_{E,\mathrm{an}},
\)
and is constant over $Y_L$: the images
\(
    E_t=\operatorname{im}(P_{E,t})
\)
form, hence, the required concordance.
\end{proof}

\subsubsection{Reduction of Sphere Bundles}
\label{rat:sec:relative-smale-reduction}

The analytification argument of
Lemma~\ref{rat:lem:analyticization-faces} applies directly to unit sphere
bundles of oriented Euclidean bundles of rank three. We now explain why an
arbitrary oriented smooth $S^2$-bundle can be brought into this form,
relative to any locus on which such a presentation has already been
chosen.

\begin{lemma}
\label{rat:lem:relative-smale-reduction}
Let $(W,L)$ be a finite CW pair, and let
\[
    \pi:C\longrightarrow W
\]
be an oriented smooth $S^2$-bundle. Suppose that over $L$ a reduction of
its structure group to $SO(3)$ has been chosen. Equivalently, suppose
that there are an oriented Euclidean vector bundle $V_L\to L$ of rank
three and an identification
\(
    C|_L\cong S(V_L).
\)

Then there are an oriented Euclidean vector bundle $V\to W$ extending
$V_L$ and an oriented $S^2$-bundle
\[
    \mathcal C\longrightarrow W\times[0,1]\qquad\text{with}\qquad \mathcal C|_{W\times\{0\}}\cong C,
    \quad
    \mathcal C|_{W\times\{1\}}\cong S(V),
\]
and, over $L\times[0,1]$, the concordance is the constant bundle
\(
    S(V_L)\times[0,1].
\)
\end{lemma}

\begin{proof}
Put
\(
    G=\operatorname{Diff}^+(S^2),
    \
    H=SO(3).
\)
Smale's theorem \cite{Smale1959} states that the natural inclusion
\(
    H\longrightarrow G
\)
is a homotopy equivalence. The rest is standard.
\end{proof}

\subsubsection{Parametric Dolbeault Operators}
\label{rat:sec:parametric-dolbeault}

The preceding arguments produce a real-analytic complex vector bundle on
a real-analytic family of rational curves. To obtain a holomorphic bundle
over a complexification of the parameter space, we use the following
elementary parametric form of the Dolbeault lemma.

\begin{definition}
    Let $\pi:C\longrightarrow S$  be a real-analytic family of Riemann surfaces, and let $E\to C$ be a real-analytic complex vector bundle. Writing $\cO(E)$ for the real analytic sections of $E$, a \emph{real analytic vertical Dolbeault operator} on $E$ is a $\mathbb{C}$-linear map $D \colon \cO(E) \to \Omega^{0,1}_{C/S}\otimes E$ satisfying the Leibniz rule
    \begin{equation*}
        D(fs) = \left(\bar\partial_{C/S} f\right)\,s + f\,Ds
    \end{equation*}
    for real-analytic functions $f$ and sections $s$. Such an operator is equivalent to prescribing a holomorphic structure on each fiber of $E$. A local frame $\{e_i\}$ is \emph{$D$-holomorphic along the fibres of} $\pi$ if $De_i = 0$ for all $i$.
\end{definition}

\begin{lemma}
\label{rat:lem:parametric-dolbeault}
Let
\(
    \pi:C\longrightarrow S
\)
be a real-analytic family of Riemann surfaces, and let $E\to C$ be a
real-analytic complex vector bundle. Suppose that $E$ is equipped with a
real-analytic vertical Dolbeault operator
\(
    D:\cO(E)\longrightarrow
    \Omega^{0,1}_{C/S}\otimes E.
\)
Then every point of $C$ has a neighborhood on which $E$ admits a frame
which is $D$-holomorphic along the fibers of $\pi$ and real-analytic in
all variables.
\end{lemma}

\begin{proof}
This is standard.
\end{proof}

\subsection{Relative Analytic Realization}
\label{rat:sec:relative-analytic-realization}

\subsubsection{}

Let $\iota:C\hookrightarrow X_T$ be a smooth family over $T$ of embedded rational curves in a 
CY\(4\)-family \(X_T\to T\), and let
$N=N_{C/X_T}$
be its rank-three complex normal bundle.  The relative tangent sequence
and the volume form give a canonical determinant isomorphism $\theta_\Omega:
    \det N\xrightarrow{\sim}K_{C/T}$.

This motivates the following definition in the topological setting. 

\begin{definition}
\label{rat:def:rational-normal-datum}
    A \textit{rational normal datum} over a space \(T\) consists of

    \begin{enumerate}[label=(\alph*)]
        \item an oriented smooth \(S^2\)-bundle
        $\pi:C\longrightarrow T$
        with a fiberwise complex structure compatible with its
        orientation;

        \item a rank $3$-complex vector bundle
        $N\longrightarrow C$;

        \item an isomorphism
        \[
            \theta:\det N\xrightarrow{\sim}K_{C/T}.
        \]
    \end{enumerate}

    Isomorphisms of rational normal data are required to preserve
    \(\theta\).
\end{definition}

Note that for holomorphic rational normal data, i.e. when the fiberwise complex structure can be upgraded to a holomorphic structure of the family making the vector bundle holomorphic, the total space
\[
    p:\operatorname{Tot}_C(N)\longrightarrow T
\]
is naturally a CY\(4\)-family, since
\[
    K_{\operatorname{Tot}_C(N)/T}
    \cong
    p^*\left(
        K_{C/T}\otimes(\det N)^{-1}
    \right)
    \xrightarrow[\theta]{\sim}
    \mathcal O.
    \label{rat:eq:local-volume-form}
\]

\subsubsection{}
\begin{definition}\label{def:N0-classifying-space}
Let $\mathscr N_0$ be the topological groupoid whose objects are triples
\[
 (J,N,\theta),\qquad
 \theta\colon\det N\xrightarrow{\sim}K_{(S^2,J)},
\]
where $J$ is a fiberwise complex structure compatible with the orientation of
$S^2$ and $N$ is a rank $3$-complex vector bundle on $S^2$.

An isomorphism is an orientation-preserving isomorphism of the complex
spheres, together with a complex bundle isomorphism preserving $\theta$.
We define
\[
    \mathcal N_0:=B\mathscr N_0 
\]
to be the classifying space of rational normal data, which is equivalently the homotopy
quotient of the classifying groupoid, with fixed parametrization of $S^2$, by
\(\operatorname{Diff}^+(S^2)\).

Equivalently, after fixing the standard complex sphere $\mathbb P^1$, put
\[
 \mathcal F_0=
 \operatorname{hofib}_{[K_{\mathbb P^1}]}\!
 \left(
   \operatorname{Map}_{-2}(\mathbb P^1,BU(3))
   \xrightarrow{\det}
   \operatorname{Map}_{-2}(\mathbb P^1,BU(1))
 \right).
\]
Here the subscript denotes the component of bundles of first Chern
number $-2$.  

Then
\(
             \mathcal N_0\simeq
             \mathcal F_0\times_{\operatorname{Aut}(\mathbb P^1)}E\operatorname{Aut}(\mathbb P^1).
\)
\end{definition}

\subsubsection{}
Suppose one has a \(d\)-parameter family of rational normal data, hence a map
to \(\mathcal N_0\).  To replace this representative by a geometric local
curve family, we must realize its smooth normal data by compatible
simplexwise holomorphic data.  We will prove the corresponding low-dimensional  statements below.

\begin{theorem}
\label{rat:thm:relative-analytic-realization}
    Let \(W\) be a compact smooth manifold of dimension at most three,
    possibly with boundary, equipped with a finite triangulation \(K\).
    Let \(L\subseteq K\) be a subcomplex, and let
    \[
        \pi:C\longrightarrow W,
        \qquad
        N\longrightarrow C,
        \qquad
        \theta:\det N\xrightarrow{\sim}K_{C/W}
        \label{rat:eq:smooth-normal-data}
    \]
    be smooth rational normal data.

    Suppose that, over a neighborhood of \(|L|\), these data are already
    the underlying smooth data of compatible simplexwise holomorphic
    rational normal data.  Then, after
    \begin{enumerate}
        \item a homotopy through smooth rational normal data, fixed over a
        neighborhood of \(|L|\), and

        \item a finite subdivision of \(K\), relative to \(L\),
    \end{enumerate}
    the following data may be chosen for every simplex
    \(\sigma\subseteq K\):
    \begin{enumerate}[label=(\alph*)]
        \item an open neighborhood
        \(
            U_\sigma\subseteq\sigma_{\mathbb C}
        \)
        of \(\sigma\) in its affine complexification;

        \item a proper holomorphic family of rational curves
        \(
            \pi_\sigma:C_\sigma\longrightarrow U_\sigma;
        \)

        \item a holomorphic rank-three vector bundle
        \(
            N_\sigma\longrightarrow C_\sigma;
        \)

        \item a holomorphic determinant isomorphism
        \[
            \theta_\sigma:
            \det N_\sigma\xrightarrow{\sim}
            K_{C_\sigma/U_\sigma}.
            \label{rat:eq:holomorphic-normal-data}
        \]
    \end{enumerate}
    These choices satisfy the following conditions.
    \begin{enumerate}[label=(\roman*)]
        \item Their restrictions to the real simplex \(\sigma\) are
        homotopic to the original smooth data
        \eqref{rat:eq:smooth-normal-data}, through the homotopy chosen
        above.

        \item If \(\tau\subseteq\sigma\) is a face, the restriction of the
        data over \(U_\sigma\) to a neighborhood of \(\tau\) agrees with
        the data over \(U_\tau\).  The identifications may be chosen so
        that all chains of face restrictions agree.

        \item Over \(L\), the original holomorphic data and all their
        chosen face identifications are unchanged.

        \item The total space
        \(
            \operatorname{Tot}_{C_\sigma}(N_\sigma)
            \longrightarrow U_\sigma
            \label{rat:eq:analytic-local-fourfold}
        \) 
        is a holomorphic CY\(4\)-family, and its zero section
        is a proper holomorphic family of rational curves.
    \end{enumerate}
\end{theorem}

\subsubsection{}

\begin{proof}[Proof of
Theorem~\ref{rat:thm:relative-analytic-realization}]
Let $U\subseteq W$ be a neighbourhood of $|L|$ on which the
holomorphic data are prescribed. After subdivision relative to
$L$, choose a
subcomplex $P$ with
\[
    |L|\subseteq\operatorname{Int}_W|P|\subseteq |P|\subseteq U.
\]
We make all constructions relative to $P$, so the resulting
homotopy is fixed on a neighbourhood of $|L|$.

\medskip\noindent
\textit{Step 1.}
Over $P$, the prescribed family has structure group
$PGL_2(\mathbb C)$. The bundle of its reductions to
$PSU(2)\cong SO(3)$ has contractible fibre
$PGL_2(\mathbb C)/PSU(2)$, so it admits a smooth section
\cite[\S9.4]{Steenrod1951}. After subdividing, apply
Lemma~\ref{rat:lem:exact-face-approximation} to this section
in local analytic trivializations. This gives an analytic
reduction on each simplex, compatible on common faces.

By Lemma~\ref{rat:lem:relative-smale-reduction}, we have a concordance,
fixed over $P$, from $C\to W$ to the unit sphere bundle $S(V)$
of an oriented Euclidean bundle $V$ of rank three.
The space of complex structures compatible with an orientation
on a real plane is contractible. Thus we may also deform the
fibrewise complex structures to those of the round spheres,
fixed over $P$. Transport $N$ and $\theta$ through these
concordances.

Apply Lemma~\ref{rat:lem:analyticization-faces} to $V$ and its
metric, relative to $P$. We may now choose simplexwise analytic
orthonormal frames. Their transition maps and face
identifications extend holomorphically, with values in
$SO(3,\mathbb C)\cong PGL_2(\mathbb C)$, and
the standard action on $\mathbb P^1$ gives proper holomorphic
families $C_\sigma\to U_\sigma$ over sufficiently small complex
neighbourhoods of the simplices. By the holomorphic uniqueness, we may extend the cocycle and face
composition identities as well.

\medskip\noindent
\textit{Step 2.}
Choose a compatible simplexwise real-analytic Hermitian metric
on the resulting bundle over $P$, and extend its Chern
connection smoothly to $N$.
Apply Lemma~\ref{rat:lem:analyticization-faces} simultaneously
to $N$, $\theta$, and this connection, relative to $P$.
They are now simplexwise real-analytic, and $\theta$ remains
an isomorphism.

Let $D_0$ be the vertical $(0,1)$-part of the connection, and
let $D_{\det,0}$ be its induced operator on $\det N$.
To make $\theta$ holomorphic along the fibres, put
\[
    a=\theta^{-1}\circ\bar\partial_{K_{C/W}}\circ\theta
       -D_{\det,0}
    \in\Omega^{0,1}_{C/W},
\]
and set
\begin{equation}
    D=D_0+\frac13a\,\operatorname{id}_N.
    \label{rat:eq:trace-correction}
\end{equation}
This operator $D$ induces the required holomorphic
structure on $\det N$, moreover, over $P$, nothing changes since $a=0$.

The equality $D^2=0$ holds since the fibres are curves, and
there are no vertical $(0,2)$-forms.

\medskip\noindent
\textit{Step 3.}
By Lemma~\ref{rat:lem:parametric-dolbeault}, $N$ has local
frames that are $D$-holomorphic along the fibres and real-analytic in all variables. Their transition matrices are
therefore holomorphic in the curve coordinate and extend
holomorphically in the parameter coordinates.
The same holds for $\theta$.
After reducing the neighbourhoods $U_\sigma$, one gets therefore holomorphic bundles $N_\sigma$ and isomorphisms
\[
    \theta_\sigma:\det N_\sigma
    \xrightarrow{\sim}K_{C_\sigma/U_\sigma}.
\]

The rest should be clear by now.
\end{proof}

\subsection{The Rational Normal Data Map}
\label{subsec:rational-normal-data-map}

\subsubsection{Rational Semi-Charts \& Thom Representatives}
For a smooth rational cycle \(z=m[C_b]\in R_m^{\mathrm{sm}}(X/B)\) in the interior semi-regular
locus, let \((S_z,z)\) denote the semi-chart of
Definition~\ref{inc:def:semi-chart}, obtained by deforming
\(C_b\cong\mathbb{P}^1\hookrightarrow X_b\) in the family \(X/B\), with
multiplicity \(m\) retained. If \(B\) is an oriented real \(n\)-manifold, then \(S_z\) is an
oriented real manifold of dimension \(n+2\).

\begin{lemma}
\label{rat:lem:degrees-of-thom-representatives}
    Let \(W\subseteq S\) be a compact oriented \(d\)-dimensional
    submanifold.  The ordinary Thom class of \(W\subseteq S\) defines a
    bivariant class
    \begin{equation}
        \operatorname{PD}_S[W]
        \in H_c^{2-d}(S/B),
        \label{rat:eq:thom-degree}
    \end{equation}
    which we call the \textit{Thom representative of the Poincar\'e dual}. 
\end{lemma}

\begin{proof}
    This is standard dimension counting.
\end{proof}

\begin{remark}
    Since Thom classes are compatible with restriction to open neighbourhoods of $W$, the image of $\operatorname{PD}_S[W]$ in $H_c^{2-d}(R_m^\sm/\CY_4)$ does not depend on the open neighbourhood $S$ of $W$ in $R_m^\sm(X/B)$, as long as $S$ is a smooth manifold, e.g. any $V$-regular locus or smooth deformation space glued from semi-charts.
\end{remark}

\subsubsection{}
\begin{lemma}\label{rat:lem:thom-restriction}
    Let $X\to B\times B_0V$ be a family of CY$4$-folds, with $B$ an oriented manifold, and
    let $X'\coloneqq X|_{B'\times B_0V}$ be the restricted family over a subcomplex $\iota\colon B'\hookrightarrow B$ which is an oriented manifold. Let $S$ and $S'$ be the open sets of
    $V$-regular cycles in $R_m^\sm(X/B\times B_0V)$ and $R_m^\sm(X'/B'\times B_0V)$
    (see Lemma~\ref{gt:lemma:v-reg-open}). Then $S'=S\times_BB'$, and for every compact
    oriented $d$-manifold $W\subseteq S'$,
    \[
        \iota_*\operatorname{PD}_{S'}[W]=\operatorname{PD}_S[W]
        \in H_c^{2-d}(S/B\times B_0V).
    \]
    In particular, the Thom representatives of $W$ in $X'$ and in $X$ define the same
    class in $H_c^{2-d}(R_m^\sm/\CY_4)$.

    More generally, if $B_0,B_1\subseteq B$ are two such subcomplexes and
    $W_i\subseteq S\times_BB_i$ are compact oriented $d$-manifolds which are homologous in $S$,
    then the Thom representatives of $W_0$ in $X|_{B_0\times B_0V}$ and of $W_1$ in
    $X|_{B_1\times B_0V}$ define the same class in $H_c^{2-d}(R_m^\sm/\CY_4)$.
\end{lemma}

\begin{proof}
    $V$-regularity of  a cycle over $(b,v)$ depends only on the restriction of the family over
    $\{b\}\times B_0V$, which coincides for $X$ and $X'$, giving us $S'=S\times_BB'$.
    At $z=m[C]\in S$ over $(b,v)$, the tangent space $T_zS$ maps onto the kernel of the
    Kodaira--Spencer map $T_bB\oplus V\to H^1(C,N_{C/X_{(b,v)}})$. Its restriction to $V$ is
    surjective by $V$-regularity, so this kernel projects onto $T_bB$, showing that $S\to B$ is a submersion.
    Thus, $S'$ is a submanifold of $S$ with normal bundle pulled back from that of
    $B'\subseteq B$, and the Gysin map $\iota_*$ of $B'\subseteq B$ pulls back to the
    Gysin map of $S'\subseteq S$. Under Poincar\'e duality, the latter is the
    pushforward $H_d(S')\to H_d(S)$, which sends $[W]$ to $[W]$. 
    
    The last statement follows, since both classes are images of
    $\operatorname{PD}_S[W_0]=\operatorname{PD}_S[W_1]$.
\end{proof}

\subsubsection{}
A family of rational normal data
$(C\to T,N,\theta)$ determines a classifying map
$T\to\mathcal N_0$, well-defined up to homotopy.

In particular, a closed oriented manifold $W$ equipped with rational
normal data defines a bordism class
\[
    [W,C,N,\theta]\in\Omega_*^{\mathrm{SO}}(\mathcal N_0).
\]

\begin{proposition}
\label{prop:rational-normal-data-bordism}
For $m\geq1$ and $0\leq d\leq2$, analytic realization and
generic transversality define a surjective
homomorphism
\[
 \Phi_{m,d}\colon
 \Omega_d^{\mathrm{SO}}(\mathcal N_0)\otimes\mathbb Q
 \twoheadrightarrow
 H_c^{2-d}(R_m^{\mathrm{sm}}/\CY_4;\mathbb Q)\xrightarrow{\sim} H_c^{2-d}(R_m^{\sm,\sr\circ}/\CY_4;\mathbb Q),
\]
where the last isomorphism is due to Theorem~\ref{gt:thm:pardon-ring-sro-isom}.
\end{proposition}

\begin{proof}
Represent a bordism class by rational normal data $(C\to W,N,\theta)$.
Choose the analytic realization
\begin{equation*}
    C\subset\operatorname{Tot}_C(N)\to W
\end{equation*}
supplied by Theorem~\ref{rat:thm:relative-analytic-realization}. Consider $W$ as a subspace of the cycle space via
\begin{equation*}
    W\hookrightarrow\cZ(\Tot_C(N)/W),\ w\mapsto m[C_w].
\end{equation*}
$W$ comes with a finite triangulation $K$ and the map to the base is the identity and hence injective. Thus, we may apply Proposition~\ref{cor:maincor}. After trimming, it gives a deformation over $W\times B_0V$ for which all cycles $w\in W$ are $V$-regular. Regarding $\tilde X\to W\times B_0V$ as a family over $W$ of complex manifolds $\tilde X_w\coloneqq\tilde X|_{\{w\}\times B_0V}$, $V$-regularity says that $C_w$ is regular in $\tilde X_w$, i.e.\ $H^1(C_w,N_{C_w/\tilde X_w})=0$, see the proof of Lemma~\ref{gt:lemma:v-reg-open}. Hence the semi-charts relative to $W$ are smooth near $W$, and they glue over the simplices of $W$ using \cite[Lemma~2.8]{pardon2023} to a smooth chart $S\to W$ around $W$, namely the $V$-regular locus near $W$. Lemma~\ref{rat:lem:degrees-of-thom-representatives} gives us a corresponding Thom representative of $W$ in this space, which defines
$\Phi_{m,d}[W,C,N,\theta]$ via extension by $0$, and has bivariant degree $2-d$.

This class is well-defined, i.e. independent of choices, in the colimit. It does not depend on the deformation: for two such deformations $\tilde{X}\to W\times B_0V$ and $\tilde{X}'\to W\times B_0V'$, the deformation to the normal cone of $C$ replaces them by $\Tot(E)$ and $\Tot(E')$, where $E,E'$ are the normal bundles of $C$ relative to $W$, keeping $C$ regular, and both are restrictions of $\Tot$ of the pushout of $E\oplus E'$ along $N\oplus N\to N$. So all four Thom representatives agree by Lemma~\ref{rat:lem:thom-restriction}. This argument is more precisely given in the surjectivity proof below.

For bordism invariance, let $M$ be an oriented bordism with rational normal data between $W_0$ and $W_1$, extended by collars $W_i\times[0,1]$ along its boundary over which the data are constant, so that $W_0$ and $W_1$ lie in its interior. Theorem~\ref{rat:thm:relative-analytic-realization}, with $L$ the two collars, gives an analytic realization over $M$ restricting to the given ones, and Proposition~\ref{cor:maincor} gives us a deformation of it over $M\times B_0V$ in which $C$ is $V$-regular. As $W_0$ and $W_1$ are homologous via $M$ in the $V$-regular locus, Lemma~\ref{rat:lem:thom-restriction} shows that the restrictions of this deformation to $W_0\times B_0V$ and $W_1\times B_0V$ define the same class $\Phi_{m,d}[W_0,C,N,\theta]=\Phi_{m,d}[W_1,C,N,\theta]$. The case $M=W\times[0,1]$ can be used to show independence of the analytic realization.

\subsubsection{Surjectivity}
Let
$\alpha\in H_c^{2-d}(R_m^{\mathrm{sm}}/\CY_4;\mathbb Q)$,
where $0\leq d\leq2$. After the graph trick Lemma~\ref{lemma:graph-trick}, $\alpha$ is represented by a class
supported on a compact subanalytic set $K$ which maps injectively to $B$. As in
\cite[Pf. of Theorem~5.15, p.46]{pardon2023}, embedding $B$ into a triangulated cube $\hat{B}$
and pulling $X\to B$ back along a local retraction $\hat{B}\to B$, we may assume that $B$
is an oriented manifold near the image of $K$. The pushforward along $B\subseteq\hat{B}$
does not change the class. Then a regularizing deformation
using Proposition~\ref{cor:maincor} gives us that
the Poincar\'e dual of $\alpha$ is a class in $H_d(S;\mathbb Q)$ for a smooth
rational semi-chart $S$ given by a smooth deformation space. Here, as in the construction, we form the semi-charts relative to $B$, regarding $\tilde X\to B\times B_0V$ as a family over $B$ with fibres $\tilde X|_{\{b\}\times B_0V}$, in which the cycles near $K$ are regular. They glue over the simplices of $B$ by \cite[Lemma~2.8]{pardon2023} to the smooth chart $S$. Since
$\Omega_q^{\mathrm{SO}}=0$ for $1\leq q\leq3$, the
Atiyah--Hirzebruch spectral sequence gives
\begin{equation}\label{eq:semi-chart-AHSS-isom}
    \Omega_d^{\mathrm{SO}}(S)\otimes\mathbb Q
 \xrightarrow{\ \sim\ } H_d(S;\mathbb Q),
 \qquad 0\leq d\leq2
\end{equation}
for any semi-chart $S$. Consequently, $\alpha$ is a finite rational sum of Thom representatives
of compact closed oriented $d$-manifolds $W_a$ which, up to homotopy and subdivision, using the graph trick, may be taken so that its triangulation makes it a subcomplex. Hence, it suffices to consider the Thom representative of a compact $d$-manifold $W$ of $V$-regular cycles $m[C_w]$ in a family $X\to B\times B_0V$, such as the regularizing
deformation above, with $W\subseteq B\times B_0V$ a subcomplex. Let $(C\to W,N,\theta)$ be
the induced rational normal data.

First, Proposition~\ref{cor:maincor} gives a deformation
$\widetilde X\to B\times B_0V\times B_0V_0$ of $X$ in which the cycles of $W$, lying over
$W\times\{0\}$, are $V_0$-regular. By Lemma~\ref{rat:lem:thom-restriction}, applied to the subcomplex
$B\times\{0\}\subseteq B\times B_0V_0$, the oriented manifold $W\subseteq B\times\{0\}$, and the additional $V$-directions, the Thom representative of $W$ is the same in $X$ and in $\widetilde X$. Applied instead to the subcomplex
$W\subseteq B\times B_0V$ with the additional $V_0$-directions, it is also the same as in
\[
    \widetilde X_0\coloneqq\widetilde X|_{W\times B_0V_0}\to W\times B_0V_0.
\]
Secondly, since $(C,N,\theta)$ is simplexwise holomorphic, $\Phi_{m,d}[W,C,N,\theta]$ is by definition the Thom representative of $W$ in a
deformation
\[
    \widetilde X_1\to W\times B_0V_1
\]
of $\Tot_C(N)\to W$ given by Proposition~\ref{cor:maincor}, in which $C$ is $V_1$-regular.

We compare $\widetilde X_0$ and $\widetilde X_1$ by deforming both to the normal cone of $C$.
For $i=0,1$, let $E_i$ be the normal bundle of $C$ in $\widetilde X_i$ relative to $W$. Since
$N$ is the normal bundle of $C$ in $\widetilde X_i|_{W\times\{0\}}$, this is an extension
\begin{equation}\label{rat:eq:surj-dnc-ses}
    0\longrightarrow N\longrightarrow E_i\longrightarrow\pi^*V_i\longrightarrow0,
\end{equation}
whose connecting map is the Kodaira--Spencer map $\kappa_i$, as in the proof of
Lemma~\ref{gt:lemma:v-reg-open}. Via $E_i\to\pi^*V_i$, the total space $\Tot(E_i)$ is a
CY$4$-family over $W\times V_i$, with volume form induced by $\theta$. The deformation to the
normal cone of $C$ in $\widetilde X_i$, together with that of $W\times\{0\}$ in
$W\times B_0V_i$, is a family $Y_i\to W\times \bR \times B_0V_i$ with
\[
    (Y_i)_{(w,t,v)}=(\widetilde X_i)_{(w,tv)}\quad(t\neq0),\qquad
    (Y_i)_{(w,0,v)}=\Tot(E_i)_{(w,v)},
\]
containing $C\times[0,1]$ over $W\times[0,1]\times\{0\}$. For each $(w,t)$, the Kodaira--Spencer map $V_i\to H^1(C_w,N_w)$ of $m[C_w]$ over $(w,t,0)$ is $t\kappa_i$ for $t\neq0$, since the slice of $Y_i$ over $(w,t)$ is that of $\widetilde X_i$ over
$w$ pulled back along $v\mapsto tv$. For $t=0$, it is the connecting map of
\eqref{rat:eq:surj-dnc-ses}, i.e.\ again $\kappa_i$. Hence $C\times[0,1]$ is $V_i$-regular
in $Y_i$.  Then Lemma~\ref{rat:lem:thom-restriction} with
$W\times\{0\},W\times\{1\}\subseteq W\times[-1,2]$ shows that the Thom representatives of
$W$ in $\widetilde X_i$ and in $\Tot(E_i)$ agree.

Finally, the pushout $E\coloneqq(E_0\oplus E_1)/N$ along the addition $N\oplus N\to N$, which
identifies the two copies of $N$, is an extension of $\pi^*(V_0\oplus V_1)$ by $N$. The family
$\Tot(E)\to W\times V_0\times V_1$ restricts to $\Tot(E_0)$ over $W\times V_0\times\{0\}$ and to
$\Tot(E_1)$ over $W\times\{0\}\times V_1$. By Lemma~\ref{rat:lem:thom-restriction} applied to
the subcomplexes $W\times\{0\}\subseteq W\times B_0V_1$ and $W\times\{0\}\subseteq W\times B_0V_0$,
with $V_0$ resp.\ $V_1$ as the regularizing directions, the Thom representatives of $W$ in
$\Tot(E_0)$ and in $\Tot(E_1)$ both agree with that in $\Tot(E)$.
This proves $\alpha=\Phi_{m,d}[W,C,N,\theta]$.
\end{proof}

\subsubsection{}\label{rnd:sec:rational-point-generator}
The discussion so far already lets us immediately determine a generator for the degree $2$ part of the cohomology.

\begin{definition}\label{rnd:def:rational-point-generator}
    For $m\geq 1$, we define the \textit{degree $2$ standard rational generator}
    \begin{equation*}
        p_m\in H_c^2(R_m^\sm/\CY_4;\bQ)
    \end{equation*}
    to be given by the Poincar\'e dual of a point in a semi-chart inside $R_m^\sm$.
\end{definition}

\begin{corollary}
    $H_c^2(R_m^\sm/\CY_4;\bQ)$ is generated by $p_m$.
\end{corollary}

\subsection{The Degree Two Normal Generator}
\label{rat:sec:degree-two-normal-generator}

We finish the argument by computing the part of
$\mathcal N_0$ relevant in degrees at most two.

\subsubsection{}

\begin{lemma}
\label{rat:lem:normal-data-low-homology}
There are natural isomorphisms
\begin{equation}
    H_1(\mathcal N_0;\mathbb Q)=0,
    \qquad
    H_2(\mathcal N_0;\mathbb Q)\cong\mathbb Q.
    \label{rat:eq:normal-data-low-homology}
\end{equation}
Moreover,
\begin{equation}
    \Omega_2^{SO}(\mathcal N_0)\otimes\mathbb Q
    \xrightarrow{\sim}
    H_2(\mathcal N_0;\mathbb Q).
    \label{rat:eq:normal-bordism-homology}
\end{equation}
\end{lemma}

\begin{proof}
First fix the standard oriented sphere $S^2$ and put
\(
    A_3=
    \operatorname{Map}_{-2}(S^2,BU(3)),
\)
where the subscript denotes the component on which the first Chern class
has degree $-2$. The determinant gives a map
\[
    \det:
    A_3\longrightarrow
    A_1=
    \operatorname{Map}_{-2}(S^2,BU(1)).
\]
Let $\mathcal F$ be its homotopy fiber over the map classifying
$K_{S^2}$. The choice of a point in $\mathcal F$ consists of a rank $3$
complex bundle $N$ on $S^2$, together with an identification
\(
    \det N\cong K_{S^2}.
\)

Let $x$ and $h$ denote the positive generators of the cohomology of the
curve and parameter factors in
\(
    S^2\times S^2.
\)
A $2$-sphere in $A_3$ is represented by a rank three bundle on this
product. Since the product has real dimension four, such a bundle is
classified in the relevant component by
\(
    c_1=-2x+ah,
    \
    c_2=bxh
\)
for integers $a$ and $b$. Thus
\(
    \pi_2(A_3)\cong\mathbb Z\oplus\mathbb Z,
\)
with coordinates $(a,b)$. Similarly,
\(
    \pi_2(A_1)\cong\mathbb Z.
\)
The determinant map records $a$ and forgets $b$:
\[
    \det_*:\pi_2(A_3)\longrightarrow\pi_2(A_1),
    \qquad
    (a,b)\longmapsto a.
\]
Since
\(
    \pi_3(BU(3))=\pi_3(BU(1))=0,
\)
we also get
\(
    \pi_1(A_3)=\pi_1(A_1)=0
\)
from the evaluation fibration.
The long exact sequence of the homotopy fiber therefore gives
\(
    \pi_1(\mathcal F)=0,
    \
    \pi_2(\mathcal F)
    =
    \mathbb Z\langle b\rangle.
\)
By the Hurewicz theorem,
\(
    H_1(\mathcal F;\mathbb Q)=0,
    \
    H_2(\mathcal F;\mathbb Q)\cong\mathbb Q.
\)

Allowing the oriented sphere itself to vary replaces $\mathcal F$ by its
homotopy quotient under
\(
    \operatorname{Diff}^+(S^2).
\)
By Smale's theorem,
\(
    \operatorname{Diff}^+(S^2)\simeq SO(3).
\)
Since
\(
    H_i(BSO(3);\mathbb Q)=0
\) for $0 < i < 4$, 
the Serre spectral sequence of
\[
    \mathcal F
    \longrightarrow
    \mathcal F_{hSO(3)}
    \longrightarrow
    BSO(3)
\]
shows that passing to the homotopy quotient does not change rational
homology in degrees one and two. This proves
\eqref{rat:eq:normal-data-low-homology}.

Finally,
\(
    \Omega_1^{SO}=\Omega_2^{SO}=0.
\)
The Atiyah--Hirzebruch spectral sequence therefore has
no additional coefficient terms in total degree two and gives
\eqref{rat:eq:normal-bordism-homology}.
\end{proof}

\subsubsection{}

We now choose an explicit generator. Let
\(
    B= C=\mathbb P^1,
\)
and let
\(
    L=\mathcal O_B(1).
\)
On
\(
    C\times B
\)
define the rank  $3$-bundle
\begin{equation}
    N_{\mathrm{Bor}}
    =
    \bigl(\mathcal O_C\boxtimes L\bigr)
    \oplus
    \bigl(\mathcal O_C(-1)\boxtimes L\bigr)
    \oplus
    \bigl(\mathcal O_C(-1)\boxtimes L^{-2}\bigr).
    \label{rat:eq:Borel-normal-bundle}
\end{equation}
Its determinant is
\[
    \det N_{\mathrm{Bor}}
    =
    \mathcal O_C(-2)\boxtimes\mathcal O_B
    =
    K_{C\times B/B}.
\]
Thus it defines an RND family over $B$.

Let
\(
    x=c_1(\mathcal O_C(1)),
    \ 
    h=c_1(L),
\)
oriented so that
\(
    \int_{C\times B}xh=1.
\)
The Chern roots of \eqref{rat:eq:Borel-normal-bundle} are
\(
    h,\qquad -x+h,\qquad -x-2h.
\)
Consequently,
\(
    c_1(N_{\mathrm{Bor}})=-2x
.\)
and
\[
\begin{aligned}
    c_2(N_{\mathrm{Bor}})
    &=
    h(-x+h)
    +
    h(-x-2h)
    +
    (-x+h)(-x-2h)  \\
    &=-xh.
\end{aligned}
\]
Therefore,
\begin{equation}
    \int_{C\times B}c_2(N_{\mathrm{Bor}})=-1.
    \label{rat:eq:Borel-c2}
\end{equation}
By the calculation in
Lemma~\ref{rat:lem:normal-data-low-homology}, this family represents a
generator of
\(
    H_2(\mathcal N_0;\mathbb Q).
\)

\begin{definition}\label{rat:def:degree-zero-rational-generator}
Define the \textit{degree $0$ standard rational generator} to be
\begin{equation}
    r_m
    =
    \Phi_{m,2}
    [B,C,N_{\mathrm{Bor}},\theta_{\mathrm{Bor}}]
    \in
    H_c^0(R_m^{\mathrm{sm}}/\CY_4;\mathbb Q).
    \label{rat:eq:degree-zero-rational-generator}
\end{equation}

\end{definition}

\subsubsection{}

\begin{proposition}
\label{rat:prop:rational-surface-reduction}
Let $W$ be a closed oriented surface and let
\(
    C_W\subset X
\)
be a family of smooth embedded rational curves with normal bundle $N$.
Then, 
\begin{equation}
    [W,C_W\subset X;m]
    =
    -
    \left(
        \int_{C_W}c_2(N)
    \right)r_m.
    \label{rat:eq:rational-surface-reduction}
\end{equation}
\end{proposition}

\begin{proof}
By Lemma~\ref{rat:lem:normal-data-low-homology}, the rational homology
class of the map
\(
    W\longrightarrow\mathcal N_0
\)
is determined by 
\(
    \int_{C_W}c_2(N).
\)
The generator \eqref{rat:eq:Borel-normal-bundle} has $b=-1$ by
\eqref{rat:eq:Borel-c2}. Hence
\[
    [W,C_W,N]
    =
    -
    \left(
        \int_{C_W}c_2(N)
    \right)
    [B,C,N_{\mathrm{Bor}}]
\]
in
\(
    H_2(\mathcal N_0;\mathbb Q)
    \cong
    \Omega_2^{SO}(\mathcal N_0)\otimes\mathbb Q.
\)
\end{proof}
\section{Connectedness of the Rational--Rational Meeting Locus}
\label{inc:sec:proof}

In this section, we study the locus of ordinary rational meeting cycles in the sense of Definition~\ref{gt:def:ordinary-rational-meeting-cycle}. More precisely, let \(X\to B\) be a family and let
\(u_i:C_i\to X_b\), \(i=1,2\), be two embedded rational curves.  Suppose
that they meet only at
\(q=u_1(p_1)=u_2(p_2)\). For this section, additionally suppose that their
tangent lines at \(q\) are distinct (we treat this constraint in Section \ref{rcl:subsubsec:connected-generation}).

We show here that this locus is connected, in the sense that every ordinary rational meeting cycle has a neighborhood, which may be deformed to a standard one. To do so, we deform the transition maps of a finite atlas around the cycle to become linear maps. The resulting model is called \textit{the linear plumbing}. The comparison of various linear plumbings then essentially uses that we can restrict to the regular case using \S\ref{rbr}, where the normal bundle is computable.

Though our exposition is self-contained, the ideas of the present section were inspired heavily by \cite{OzbagciPopescuPampu16}, \cite{Guadagni18}, \cite{Mumford61}. In particular, as far as we understand, the word ``plumbing'' goes back to Mumford, Milnor, and Hirzebruch.

\subsection{Branch-Regular Representatives}

\subsubsection{}\label{inc:subsec:common image}
Let \(S\) be Pardon's semi-chart at
\(m_1C_1+m_2C_2\). After replacing  \(S\) by a sufficiently small germ,
the deformations parametrized by \(S\) are given by families of smooth
curves \(\pi_i:\mathcal C_i\to S\), together with maps
\(F_i:\mathcal C_i\to X_S:=X\times_BS\).  Put
\(\mathcal C_{12}:=\mathcal C_1\times_S\mathcal C_2\).  

Consider
\begin{equation}
 \operatorname{ev}:\mathcal C_{12}\longrightarrow X_S\times_SX_S,
 \qquad
 (s,x_1,x_2)\longmapsto\bigl(F_1(x_1),F_2(x_2)\bigr),
 \qquad
 \widetilde D:=\operatorname{ev}^{-1}(\Delta_{X_S/S}).
 \label{inc:eq:common image-complex}
\end{equation}
The map \(\widetilde D\to S\) forgets the two chosen points.  Since
\(C_1\) and \(C_2\) have no other common point, after shrinking near
\((p_1,p_2)\) its image in $S$ is precisely the locus \(D\subset S\) on which
the two deformed curves meet (near $q$).

Over \(D\), let \(q:D\to X_S\) denote the common image of the two
distinguished points, and let \(T_i\subset q^*T_{X/B}\) be the tangent
line to the branch indexed by \(i\).  Set
\begin{equation*}
  Q\coloneqq q^*T_{X/B}/(T_1+T_2).
\end{equation*}  
Its fibre at the original meeting is
\(Q_q=T_qX_b/(T_qC_1+T_qC_2)\), which has complex dimension two.

\subsubsection{}
\begin{lemma}
\label{inc:lem:smooth-meeting-locus}
Suppose that the map
\(u_1\sqcup u_2:C_1\sqcup C_2\to X_b\) is regular in the family
\(X\to B\), in the sense of
\cite[Definition~5.1]{pardon2023}.

Then, after shrinking near
the meeting, \(D\) is smooth of complex codimension two in
\(S\), and its normal bundle is canonically isomorphic to \(Q\).

Suppose, in addition, that the same map is regular for the family
\(X_b\to\{b\}\), i.e. suppose that the deformation problem
remains unobstructed when the point \(b\) is not allowed to vary.  Writing
\(N_i=N_{C_i/X_b}\), we have
\begin{equation*}
 N_i\simeq\cO\oplus\cO(-1)^{\oplus2},
\end{equation*}
and an isomorphism
\begin{equation}
 H^0(C_1,N_1)\oplus H^0(C_2,N_2)
   \xrightarrow{\ \sim\ } Q_q,
 \qquad
 (s_1,s_2)\longmapsto s_1(p_1)-s_2(p_2).
 \label{inc:eq:effective-incidence-map}
\end{equation}

\end{lemma}

\begin{proof}
By Pardon's definition, regularity means that the deformation theory of
the triple
\((C_1\sqcup C_2,\{p_1,p_2\},u_1\sqcup u_2)\), subject to the condition
that the two special points have the same discrete data, is
unobstructed. For the map in \eqref{inc:eq:common image-complex}, the unobstructedness of the problem subject to the meeting condition in particular tells us that \(\operatorname{ev}\) is transverse to the relative diagonal
near \((p_1,p_2)\).  Since the relative diagonal has complex codimension
four, \(\widetilde D\) is smooth of complex codimension four in
\(\mathcal C_{12}\).

The projection \(\mathcal C_{12}\to S\) has relative complex dimension two.
We claim that its restriction to \(\widetilde D\) is an immersion.  A tangent vector in a fibre $C_1\times C_2$ of
\(\mathcal C_{12}\to S\) is a pair
\((v_1,v_2)\in T_{p_1}C_1\oplus T_{p_2}C_2\).  It is tangent to
\(\widetilde D\) precisely when
\(du_1(v_1)=du_2(v_2)\). By assumption, the $u_i$ are embeddings with distinct tangent lines $\rd u_i(T_{p_i} C_i)\subset T_qX_b$ at $q$, so that the equality $du_1(v_1)=du_2(v_2)$ forces \(v_1=v_2=0\). The curves have no other common point, hence,
after one further shrinking, \(\widetilde D\to S\) identifies
\(\widetilde D\) with \(D\).  It follows that \(D\) is smooth
and has complex codimension two in \(S\).

It remains to identify the normal bundle of $D$ in $S$.  Infinitesimal motion of the
two chosen points changes their difference in evaluations by
\(du_1(v_1)-du_2(v_2)\).  The image of this is exactly
\(T_qC_1+T_qC_2\).  Consequently, after the chosen points are forgotten by $\widetilde{D}\to S$,
the remaining normal part to the equality condition of $D$ takes values in \(Q_q\).
By transversality, one has therefore the exact sequence
\begin{equation}
 0\longrightarrow T_dD\longrightarrow T_dS
   \longrightarrow
   \frac{T_qX_b}{T_qC_1+T_qC_2}
   \longrightarrow0.
 \label{inc:eq:normal-evaluation-sequence}
\end{equation}
The same construction at every point of \(D\) identifies
\(N_{D/S}\) with \(Q\).

We now work in the single fourfold \(X_b\).  For an embedded curve
\(C_i\subset X_b\), first order deformations are
\(H^0(C_i,N_i)\), and their obstruction space is \(H^1(C_i,N_i)\).
Regularity of the deformation problem with the special points constraints in particular implies
regularity after the conditions at those points are omitted.  Hence
\(H^1(C_i,N_i)=0\) for both branches. 

Adjunction and the triviality of \(K_{X_b}\) give
\(\det N_i\simeq K_{C_i}\simeq\cO(-2)\), so \(\deg N_i=-2\).
By Riemann--Roch, \(\chi(N_i)=3-2=1\), and, since \(H^1(C_i,N_i)=0\), we
obtain \(h^0(C_i,N_i)=1\).  Thus, the surjection \eqref{inc:eq:effective-incidence-map} has a
source and a target both of complex dimension two, and is therefore an
isomorphism. Finally, using the Birkhoff--Grothendieck theorem, write
\(N_i\simeq\cO(a_1)\oplus\cO(a_2)\oplus\cO(a_3)\).  The vanishing of
\(H^1(C_i,N_i)\) implies \(a_j\geq-1\) for every \(j\), while
\(a_1+a_2+a_3=-2\).  Therefore, the three integers are
\(0,-1,-1\), which proves the asserted splitting.
\end{proof}

\subsubsection{}
\label{inc:subsec:divisors}

\begin{definition}
The ordinary rational meeting cycle $z$ is \emph{branch-regular} if its
associated stable map $f_z=u_1\sqcup u_2$, with special points $A=\{p_1,p_2\}$,
is regular in the sense of Definition~\ref{gt:def:regular} as a map to the single fourfold
$X_b$. That is, $\mathrm{Ob}_A(f_z)=0$.
\end{definition}

We now apply the branch regularization result of \S\ref{rbr} to
the degree zero classes supported on the locus of ordinary rational--rational
meetings.  It is clear that each such cycle can be
represented as a sum of a finite set of cocycles $\tau_{z_\alpha}$, where $\tau_{z_\alpha}$ denotes a compactly supported Thom cocycle of some point $z_\alpha$ in its meeting stratum. We can take each point
$z_\alpha$ to correspond to some meeting
$z_\alpha=m_1C_{1,\alpha}+m_2C_{2,\alpha}\subset X_\alpha$.

\begin{theorem}\label{th:rbr}
Every point
cocycle $\tau_{z_\alpha}$ in $H^0_{\rm c}(\cZ_{\rm con}^{\rm sr\circ}/\CY_4;\mathbb Q)$ of a point $z_\alpha$ in the ordinary rational--rational meeting locus 
of $\cZ_{\rm con}^{\rm sr\circ}$
represents a finite sum
\[
 [\tau_{z_\alpha}]
   =\sum_\beta\epsilon_{\alpha\beta}[\tau_{z_{\alpha\beta}}],
\]
where $z_{\alpha\beta}$ is an ordinary rational meeting cycle
\[
 m_1C_{1,\alpha\beta}+m_2C_{2,\alpha\beta}
 \subset X_{\alpha\beta}
\]
that is branch semi-regular in the fourfold $X_{\alpha\beta}$.
\end{theorem}

\begin{proof}
  We can immediately apply Proposition~\ref{gt:prop:rational-rational-gt} to this cycle, lying over $b\in B$. The fourfolds $X_{\alpha\beta}$ are then found as fibers of the resulting deformation over a common regular value $v$ in the interior of $S_0V$ of the three projections from the three deformation spaces in the slice over $b$ to $S_0V$. Pulling the holomorphic family $\tilde X\to B_\bC\times B_0W$ of Proposition~\ref{bg:prop:family-regluing} (which is used to obtain the deformation of Proposition~\ref{gt:prop:rational-rational-gt}) back along the simplexwise real analytic map $B\times[0,1]\to B_\bC\times B_0W$, $(b',t)\mapsto(b',tv)$, see \cite[\S2.4]{pardon2023}, gives a family over $B\times[0,1]$ whose restrictions to the subcomplexes $B\times\{0\}$ and $B\times\{1\}$ are $X$ and $\tilde X|_{B\times\{v\}}$.
\end{proof}

Consequently, for the remainder of this section it is enough to consider
one branch-regular ordinary rational meeting. The preceding Lemma~\ref{inc:lem:smooth-meeting-locus} describes its
normal bundles and the infinitesimal incidence map. The next construction
fixes such a meeting and begins its deformation to the toric model: we will use this model (essentially invented by \cite{KP}) in
the local calculation.

\subsection{Linear Plumbing}
\subsubsection{A Finite Meeting Atlas}\label{inc:subsec:finite-plumbing-atlas}
Fix one branch-regular ordinary rational meeting
\(m_1C_1+m_2C_2\subset X\), and denote its meeting point by \(q\).
The curves \(C_1\) and \(C_2\) are embedded, have no other common point,
and have distinct tangent lines at \(q\). Choose a polydisc \(P\subset X\) about \(q\) with coordinates
\((x_1,x_2,\eta_1,\eta_2)\) such that
\begin{equation}
 \begin{aligned}
 C_1\cap P&=\{x_2=\eta_1=\eta_2=0\},&
 C_2\cap P&=\{x_1=\eta_1=\eta_2=0\},\\
 \Omega&=dx_1\wedge dx_2\wedge d\eta_1\wedge d\eta_2.
 \end{aligned}
 \label{inc:eq:node-volume-coordinates}
\end{equation}
We regard \(P\) as the node chart. 

\begin{lemma}
\label{inc:lem:finite-plumbing-atlas}
A neighborhood of \(C_1\cup C_2\) is covered by \(P\) and finitely many
 charts
\begin{equation}
 \begin{gathered}
 T_{i,\alpha}\subset
 \mathbb C_{z_{i,\alpha}}\times\mathbb C^3_{\xi_{i,\alpha}},
 \qquad i=1,2,\qquad 0\leq\alpha\leq N_i,\\
 C_i\cap T_{i,\alpha}=\{\xi_{i,\alpha}=0\},
 \qquad
 \Omega=dz_{i,\alpha}\wedge d\xi_{i,\alpha,1}\wedge
 d\xi_{i,\alpha,2}\wedge d\xi_{i,\alpha,3}.
 \end{gathered}
 \label{inc:eq:outer-volume-charts}
\end{equation}
Put \(T_i=\bigcup_{\alpha=0}^{N_i}T_{i,\alpha}\).  The charts may be
indexed and trimmed so that the following properties hold.

\begin{enumerate}
 \item The neighborhoods \(T_1\) and \(T_2\) are disjoint.

 \item For each \(i\), only \(T_{i,0}\) meets \(P\), and
 \(R_i:=T_{i,0}\cap P\) is a neighborhood of an annulus in
 \(C_i\cap P\), away from \(q\).  The transition between \(P\) and
 \(T_{i,0}\) is the restriction of a biholomorphism defined near \(q\).

 \item The closures of \(R_1\) and \(R_2\) in \(P\) are disjoint.
This ensures that no triple overlap involves the node chart.

 \item All transition maps between charts belonging to the same branch
 preserve the zero section and the volume form.  (On \(R_1\), regard
 \(P\) as a branch chart with coordinates
 \((z,\xi)=(x_1,x_2,\eta_1,\eta_2)\), and on \(R_2\), use
 \((z,\xi)=(x_2,-x_1,\eta_1,\eta_2)\).)
\end{enumerate}
\end{lemma}

\begin{proof}
This is clear.
\end{proof}
\subsubsection{Contracting the Transition Maps}\label{inc:subsec:alexander-contraction}
We now deform the transition maps of the finite atlas to their first-order
normal parts.

\begin{lemma}
\label{inc:lem:alexander-deformation}
Let
\[
    g:D\longrightarrow D'
\]
be a holomorphic transition map between two charts of
Lemma~\ref{inc:lem:finite-plumbing-atlas} corresponding to the same branch (one of them may be a nodal chart).  Suppose that \(g\) preserves
the zero section and the standard volume forms.  In source and target
coordinates \((z,\xi)\in\mathbb C\times\mathbb C^3\), put
\[
    \delta_t(z,\xi)=(z,t\xi)
\]
and, for \(t\ne0\), define
\begin{equation}
    g_t
       =
    \delta_{t,\mathrm{target}}^{-1}
       \circ g\circ
    \delta_{t,\mathrm{source}}.
    \label{inc:eq:alexander-deformation}
\end{equation}
After shrinking, \(g_t\) and its
inverse extend holomorphically to \(t=0\).  The limiting map is of the form
\begin{equation}
    g_0(z,\xi)=\bigl(h(z),A(z)\xi\bigr),
    \label{inc:eq:normal-derivative}
\end{equation}
where $h(z)\in \bC$ and $A(z)\in\bC^{3\times 3}$ are holomorphic in $z$ and given by the normal
derivative. Every \(g_t\) preserves the volume form.  
\end{lemma}

\begin{proof}
Since \(g\) preserves the zero section, it has a convergent expansion
\[
    g(z,\xi)
      =
    \left(
       h(z)+\sum_{|I|\geq1}a_I(z)\xi^I,\,
       \sum_{|I|\geq1}b_I(z)\xi^I
    \right).
\]
The matrix
\[
    A(z)=\bigl(b_{e_j}(z)\bigr)_{j=1}^3
\]
is invertible since $g$ is a holomorphic transition map.  Substituting in
\eqref{inc:eq:alexander-deformation}, we get
\[
    g_t(z,\xi)
      =
    \left(
       h(z)+\sum_{|I|\geq1}t^{|I|}a_I(z)\xi^I,\,
       A(z)\xi+
       \sum_{|I|\geq2}t^{|I|-1}b_I(z)\xi^I
    \right).
\]
On a smaller relatively compact domain this is a convergent power series
in \((z,\xi,t)\), including at \(t=0\), and its value there is
\eqref{inc:eq:normal-derivative}.  Applying the same argument to \(g^{-1}\)
shows that the limiting maps are inverse.

Write
\[
    \Omega=dz\wedge d\xi_1\wedge d\xi_2\wedge d\xi_3.
\]
Normal scaling has volume character \(t^3\).  Hence, for \(t\ne0\),
\begin{align*}
    g_t^*\Omega_{\mathrm{target}}
      &=
    \delta_{t,\mathrm{source}}^*
       g^*
       (\delta_{t,\mathrm{target}}^{-1})^*
       \Omega_{\mathrm{target}}  \\
      &=
    t^{-3}\delta_{t,\mathrm{source}}^*
       \Omega_{\mathrm{source}}
     =\Omega_{\mathrm{source}}.
\end{align*}
The equality at \(t=0\) follows by holomorphic continuation, or by direct check using \eqref{inc:eq:normal-derivative}.  In
particular,
\[
    h'(z)\det A(z)=1.
\]
\end{proof}

\subsubsection{The Resulting Holomorphic Family}
\label{inc:subsec:holomorphic-plumbing-family}
The transition maps constructed above are initially defined on domains
which may depend on \(t\).  Since the atlas is finite, one uniform trimming
is sufficient to fix this. The following lemma follows.

\begin{lemma}
\label{inc:lem:uniform-trimming}
The charts of Lemma~\ref{inc:lem:finite-plumbing-atlas} may be shrunk,
while still covering a neighborhood of \(C_1\cup C_2\), so that there is
an open neighborhood \(\Delta\subset\mathbb C\) of \([0,1]\) with a separated holomorphic submersion
\begin{equation}
    \mathcal X\longrightarrow\Delta.
    \label{inc:eq:plumbing-family}
\end{equation}
Its fibre over \(t=1\) is the original germ along \(C_1\cup C_2\).
The volume forms on the charts glue to a relative volume form on
\(\mathcal X\to\Delta\).  The zero sections glue to submanifolds
\(\mathcal C_i\simeq C_i\times\Delta\).
\end{lemma}

\begin{definition}
  At \(t=0\), all transition maps along each
outer branch are linear in the normal directions, while the smooth node
chart and its two coordinate axes remain unchanged.  We call this endpoint
the \emph{linear plumbing} of $C_1\cup C_2$.
\end{definition} 

\subsection{Comparison to the Toric Linear Plumbing}
\label{inc:subsec:toric-plumbing-comparison}
\subsubsection{Toric Model}\label{inc:toric-model}
Let
\begin{equation*}
  \nu:S_0=\mathrm{Bl}_{(\infty,\infty)}\bP^1\times\bP^1\to\mathbb P^1_x\times\mathbb P^1_t
\end{equation*}
be the blow-up at
\((\infty,\infty)\). Denote the exceptional curve by \(E\), and
denote the strict transform of
\(\mathbb P^1_x\times\{\infty\}\) by \(A\). On \(S_0\), set $L_1\coloneqq \nu^*\mathcal O(-1,-1)$, $L_2\coloneqq \nu^*\mathcal O(-1,-1)(E)$, and consider the fourfold
\begin{equation}
 X_0=\operatorname{Tot}_{S_0}(L_1\oplus L_2). \label{eq:cap-local-fourfold}
\end{equation}
The canonical bundle formula gives \(K_{S_0}=L_1\otimes L_2\), and hence \(X_0\) is Calabi--Yau. Inside $X_0$, the embedded curves $A$ and $E$ meet transversely in one point $q$, and
\begin{equation*}
  N_{A/X_0}\cong \cO\oplus \cO(-1)^{\oplus 2},\qquad N_{E/X_0}\cong \cO\oplus \cO(-1)^{\oplus 2}.
\end{equation*}

\subsubsection{}
We use the toric fourfold \(X_0\) and the curves \(A,E\subset X_0\)
defined in \S\ref{inc:toric-model} above. For a pointed
rational curve with normal bundle
\(N=\cO\oplus\cO(-1)^{\oplus2}\), every bundle automorphism covering
the identity of the pointed curve has the form
\begin{equation}
 \begin{pmatrix}a&\ell\\0&B\end{pmatrix},
 \qquad
 a\in\mathbb C^*,\quad B\in GL_2(\mathbb C),\quad
 \ell\in H^0(\mathbb P^1,\cO(1))^{\oplus2},
 \label{inc:eq:normal-automorphisms}
\end{equation}
and it respects the determinant trivialisation precisely when \(a\det B=1\). By construction of the similar for the bundle automorphisms, taking fibers at the meeting point $q$ gives a surjection
\begin{equation}\label{inc:eq:surjection-fiber-bundle-aut}
  \Aut_{\bP^1}(N)\twoheadrightarrow \Aut(N|_q;\,H^0(N)|_q)
\end{equation}
onto the group of automorphisms of the
normal fibre $N|_q$ which preserve the image \(H^0(N)|_q\) of $H^0(N)$ in $N|_q$.  Once its value there is
fixed, the remaining choices form the affine space
\(H^0(\cO(1)(-q))^{\oplus2}\), which is the fiber of \eqref{inc:eq:surjection-fiber-bundle-aut}.

\subsubsection{}
\begin{lemma}
\label{inc:lem:toric-linear-plumbing}
The linear plumbing of a branch-regular ordinary rational meeting can be deformed, through holomorphic germs preserving the volume form, to the linear
plumbing of \(A\cup E\subset X_0\). 
\end{lemma}

\begin{proof}
Put \(W=T_qX\) and \(T_i=T_qC_i\).  Choose nonzero vectors in \(T_1\)
and \(T_2\), and choose lifts to \(W\) of generators of
\(H^0(N_i)|_q\subset W/T_i\).  By
\eqref{inc:eq:effective-incidence-map}, these four vectors form a basis
of \(W\).  Make analogous choices for the toric meeting.  After
rescaling one vector, there is a volume-preserving linear isomorphism
\(A_q:W\to T_{A\cap E}X_0\) carrying these vectors to their
toric counterparts.

Choose an isomorphism from each pointed curve to its toric counterpart
whose derivative at the meeting point $q$ is prescribed by the restriction of \(A_q\).
By the construction of $A_q$ above, the map $\bar{A}_q$ induced by \(A_q\) on the normal fibre preserves the line
obtained by evaluating global normal sections, so by \eqref{inc:eq:surjection-fiber-bundle-aut} we may extend $\bar{A}_q$ to a bundle morphism of the form \eqref{inc:eq:normal-automorphisms} over the complete curve.
The extension respects the determinant trivialisation: the ratio of the two
trivialisations is a nowhere vanishing function on \(\mathbb P^1\), hence is
constant, and it is one at the trivialised point.

On both the regular and the toric side, identify each branch neighborhood with a neighborhood of the zero section in its normal bundle. Using the global bundle isomorphism extending $\bar{A}_q$ above, we define comparison isomorphisms between the regular and the toric branch charts, restricting to the chosen isomorphisms of pointed curves and to $\bar{A}_q$ over $q$. We may choose these tubular identifications so that, at
the meeting point $q$, their differentials send the selected lifts to \(W\) of generators of
\(H^0(N_i)|_q\subset W/T_i\) to the
corresponding toric vectors.  The resulting identification of each
branch neighborhood then has full differential \(A_q\) at \(q\).

These choices identify the node charts and the normal bundles, but they
do not yet identify the maps gluing a branch neighborhood to the node
chart on its annulus.  Let \(g_i\) and \(g_i^0\) denote the two gluing
maps, after the node charts have been identified by \(A_q\).  Although
constructed only on the attachment annulus, both maps extend across the
trivialised point because they arise from tubular identifications of the
complete normal bundles.  Their ratio is therefore a holomorphic germ
at \(q\).  Their
differentials at \(q\) agree by construction.  Hence
\(r_i=(g_i^0)^{-1}\circ g_i\) has identity differential at \(q\).  Since
the two plumbings are linear in the normal directions, it has the form
\(r_i(z,\xi)=(h_i(z),A_i(z)\xi)\), where
\(h_i(0)=0\), \(h_i'(0)=1\), and \(A_i(0)=I\).  Preservation of the
volume form says \(h_i'(z)\det A_i(z)=1\).

The contraction
\begin{equation}
 r_{i,t}(z,\xi)=
   \bigl(t^{-1}h_i(tz),A_i(tz)\xi\bigr)\quad(t\ne0),
 \qquad r_{i,0}=\operatorname{id},
 \label{inc:eq:attachment-contraction}
\end{equation}
is holomorphic at \(t=0\), by Taylor expansion, automatically holomorphic at $t\neq 0$, and has Jacobian
\(h_i'(tz)\det A_i(tz)=1\).  Thus \(g_i^0\circ r_{i,t}\) deforms \(g_i\)
to \(g_i^0\) while preserving the volume form.  The two annuli have
disjoint closures, so the contractions can be made independently.
\end{proof}

Combining the preceding lemmas, we deform the original germ first to its
linear plumbing, then to the toric linear plumbing, and finally run the
rescaling deformation of the toric germ in the opposite direction to
reach \(A\cup E\subset X_0\).

\subsubsection{Semi-Charts along the Plumbing Path}
\label{inc:subsec:plumbing-semi-charts}
To finish the argument of the present section, it remains to state the following lemma (which is an easy exercise in the techniques used in the present paper, and is left to the reader.)

Let \(I=[0,1]\), subdivided by the vertices
\(0,\frac13,\frac23,1\).  The three resulting intervals parametrize,
respectively, the deformation from the original meeting to its linear
plumbing, the comparison with the toric linear plumbing, and the reverse
rescaling deformation from the toric linear plumbing to the toric
meeting. 

\begin{lemma}
\label{lem:plumbing-cylinder-semichart}
For every simplex \(\rho\subset I\) of the subdivision into three pieces above, the germ along \(z_\rho\) has a
smooth semi-chart
\(\mathcal S_\rho\to\rho\) of relative complex dimension two which
contains the distinguished section \(z_\rho\).  These semi-charts agree
under restriction to faces.

Let \(\mathcal D_\rho\subset\mathcal S_\rho\) be the locus on which the
two deformed curves have a common image point (near the deformation \(q_t\) of the original $q$ for some $t\in\rho$).  It is smooth
of complex codimension two in the holomorphic extension of
\(\mathcal S_\rho\).  After replacing the semi-charts by smaller
representatives around \(z_\rho\), one has
\(\mathcal D_\rho=z_\rho(\rho)\), i.e. the meeting locus is given by the section $z_\rho$.  If
\(q_\rho:\mathcal D_\rho\to X_\rho\) is the common image and
\(T_{i,\rho}\subset q_\rho^*T_{X_\rho/\rho}\) is the tangent line to the
curve indexed by \(i\), then
\begin{equation}
 N_{\mathcal D_\rho/\mathcal S_\rho}
 \simeq
 Q_\rho\coloneqq
 \frac{q_\rho^*T_{X_\rho/\rho}}
 {T_{1,\rho}+T_{2,\rho}}.
 \label{inc:eq:plumbing-incidence-normal}
\end{equation}
These identifications agree on faces. 
\end{lemma}
\section{Reduction to Generators}
\label{sec:connected-reduction}

We now explain how the results of the last sections help to prove the GW/PT-correspondence for Calabi-Yau $4$-olds.

We return to the notation of 
\S\ref{subsec:convolution-logarithms}.  Thus
\(\boldsymbol\Phi=(\Phi_n)_{n\in\mathbb Z}\) and
\(\boldsymbol\Psi=(\Psi_n)_{n\in\mathbb Z}\) denote the two normalized
\emph{total} theories, including possibly disconnected supports. We write
\(L_{\Phi,n}\coloneqq [q^n]\log^\star\boldsymbol\Phi\) and
\(L_{\Psi,n}\coloneqq [q^n]\log^\star\boldsymbol\Psi\); these are the coefficients
of the corresponding \emph{connected} theories.   Put \(\delta_n\coloneqq L_{\Phi,n}-L_{\Psi,n}\).  Recall that for every \(n\),
\begin{equation}
 \delta_n(1)=0,\qquad \delta_n(\mathfrak m^2)=0,
 \qquad
 \delta_n:\mathfrak m/\mathfrak m^2
 \longrightarrow\mathbb Q
 \quad(n\in\mathbb Z).
 \label{eq:coefficient-discrepancy-indecomposables}
\end{equation}
Our strategy is to reduce, first, the equality between generating functions $L_{\Phi,n}$ and $L_{\Psi,n}$ to a statement about the locus of connected cycles, and then prove the latter using an explicit description of lifts of the corresponding generators from Sections~\ref{rat:sec:rational-normal-data} and \ref{inc:sec:proof}.

\subsection{The Connected Locus}
\label{subsec:connected-locus}

\subsubsection{}

Recall from \S\ref{subsec:augmentation-ideal} that \(\mathcal Z^+=\mathcal Z\setminus\{0\}\) denotes the locus of positive
cycles.  For a
family \(X\to B\), denote by
\(\mathcal Z_{\mathrm{con}}(X/B)\coloneqq
\{z\in\mathcal Z^+(X/B):|z|_{\mathrm{red}}\text{ is connected}\}\) the locus of connected cycles, and
let \(\mathcal Z_{\mathrm{dis}}(X/B)\) be its complement in
\(\mathcal Z^+(X/B)\). The connected components of the support
of a disconnected cycle have pairwise disjoint isolating neighbourhoods,
and the same separation persists for every sufficiently nearby effective
cycle. Consequently, \(\mathcal Z_{\mathrm{dis}}(X/B)\) is open and \(\mathcal Z_{\mathrm{con}}(X/B)\) is closed in
\(\mathcal Z^+(X/B)\).  The zero section is open and closed, so the
connected locus is also closed in \(\mathcal Z(X/B)\). 

The resulting family of closed
loci will be denoted by
\(i:\mathcal Z_{\mathrm{con}}\hookrightarrow\mathcal Z^+\). More concretely, define $H_c^*(\mathcal Z_{\mathrm{con}}/\CY_4;\mathbb Q)$ as the colimit of \(H_c^*(\mathcal Z_{\mathrm{con}}(X/B)/B;\mathbb Q)\) over the
families used in Definition~\ref{coh:def:grothendieck-group-1-cycles}.
The proper inclusions \(i_{X/B}\) give the restriction map
\begin{equation*}
    i^*:\mathfrak m=H_c^*(\mathcal Z^+/\CY_4;\mathbb Q)
\to H_c^*(\mathcal Z_{\mathrm{con}}/\CY_4;\mathbb Q).
\end{equation*}
Notice that the connected locus need not be open.

\subsubsection{Disconnected Support and the Augmentation Ideal}
\label{subsec:disconnected-square}

Write \(\mathfrak m^{[q]}\) for the part of bivariant degree \(q\).  Let
\(j:\mathcal Z_{\mathrm{dis}}\hookrightarrow\mathcal Z^+\) and
\(i:\mathcal Z_{\mathrm{con}}\hookrightarrow\mathcal Z^+\) be the open
and closed inclusions.  The localization sequence gives a boundary map
\(\partial_{\mathrm{dis}}\) from the cohomology of the connected locus to
the cohomology of the disconnected locus, with degree raised by one.

\begin{proposition}
\label{prop:disconnected-square}
The image under \(j_!\) of the cohomology supported on
\(\mathcal Z_{\mathrm{dis}}\) is \(\mathfrak m^2\).  Consequently,
restriction induces, in every bivariant degree \(q\), an isomorphism
\begin{equation}
 (\mathfrak m/\mathfrak m^2)^{[q]}\xrightarrow{\ \sim\ }K^q,
 \qquad
 K^q\coloneqq \ker\!\left[
 H_c^q(\mathcal Z_{\mathrm{con}}/\CY_4;\mathbb Q)
 \xrightarrow{\ \partial_{\mathrm{dis}}\ }
 H_c^{q+1}(\mathcal Z_{\mathrm{dis}}/\CY_4;\mathbb Q)
 \right].
 \label{eq:indecomposables-connected-kernel}
\end{equation}
Thus a class on the connected locus belongs to \(K^q\) precisely when it
has an ambient lift in \(\mathfrak m^{[q]}\), and that lift is unique
modulo \(\mathfrak m^2\).
\end{proposition}

\begin{proof}
The inclusion \(\mathfrak m^2\subseteq\operatorname{im}(j_!)\) follows
from the results of Subsection~\ref{subsec:pardon-product}: a product of
two positive classes is represented on a disjoint union of target families,
and every cycle in its support has at least two nonempty connected
components.

For the converse, use convolution of endomorphisms of the bialgebra
\(\mathcal H\).  Put \(J=\mathrm{id}_{\mathcal H}-\epsilon\) and
define the Eulerian operator
\begin{equation}
 e\coloneqq \log^\star(\mathrm{id}_{\mathcal H})
 =\sum_{r\geq1}\frac{(-1)^{r-1}}rJ^{\star r}.
 \label{eq:eulerian-operator}
\end{equation}
The sum is finite on every compactly supported class.  Indeed, if
\(z=\sum_A m_A[A]\) is a nonzero effective cycle, set
\(\lambda(z)=\sum_A m_A\) and let
$c(z)$ denote the number of connected components of $\abs{z}_\red$.  The function \(\lambda\) is bounded
on compact subsets, and a decomposition of \(z\) into \(r\) positive
subcycles requires \(r\leq\lambda(z)\).

For \(N,k\geq1\), let \(V_N^k\) be the open locus on which
\(\lambda\leq N\) and \(c\geq k\), and let \(D_N^k\mathcal H\) be the
image in \(\mathfrak m\) of its compactly supported cohomology.  The first
condition is open by the multiplicity filtration, and the
second is open because we can find isolating neighborhoods of disjoint components.  The groups
\(D_N^k\mathcal H\) form a decreasing filtration, and
\(D_N^{N+1}\mathcal H=0\).

Let \(S(k,r)\) denote a Stirling number of the second kind. On the
quotient with exactly \(k\) connected components, one has
\begin{equation}
 J^{\star r}x\equiv r!S(k,r)x\pmod{D_N^{k+1}\mathcal H},
 \qquad x\in D_N^k\mathcal H.
 \label{eq:stirling-convolution}
\end{equation}
This is because the decompositions which do not split a connected component are indexed by
the surjections \(\pi_0(|z|_{\mathrm{red}})\twoheadrightarrow
\{1,\ldots,r\}\).  They hence form a finite cover of degree \(r!S(k,r)\).

Plugging \eqref{eq:stirling-convolution} into \eqref{eq:eulerian-operator}, the identity
\(\sum_{r=1}^k(-1)^{r-1}(r-1)!S(k,r)=\delta_{k1}\) now implies
\(e(D_N^k\mathcal H)\subseteq D_N^{k+1}\mathcal H\) for \(k\geq2\).
The usual argument for convolution powers shows that \(e\) is
infinitesimal, i.e. $e(xy)=e(x)\epsilon(y)+\epsilon(x)e(y)$, and hence\(e(\mathfrak m^2)=0\). Moreover,
\(e(x)-x\in\mathfrak m^2\) for \(x\in\mathfrak m\), and these identities
give \(e^2=e\).  If \(x\in D_N^2\mathcal H\), repeated application of the
preceding inclusion and the idempotence give \(e(x)=0\). It follows that
\(x\in\mathfrak m^2\).  Every compactly supported class on the locus of
disconnected support lies in some \(D_N^2\mathcal H\), proving the reverse
inclusion.  The remaining assertions follow from the localization
sequence.
\end{proof}

\subsection{Standard Classes and Lifts}
\label{subsec:standard-connected-classes}

We have the standard generators $e_m$ from \S\ref{sec:elliptic-locus}, $r_m$ from \S\ref{rat:sec:rational-normal-data} and the toric ordinary rational meeting model from \S\ref{inc:sec:proof}. We now want to define classes from such models which exhaust the cases we have to check later. 

This means we want to achieve two things: First, we want to extend these classes to $H_c^*(\cZ_\con/\CY_4)$. Secondly, we want to find a class in $\fm=H_c^*(\cZ^+/\CY_4)$ that extends the class in $H_c^*(\cZ_\con/\CY_4)$ under \eqref{eq:indecomposables-connected-kernel}. We call such a class an \textit{ambient lift}.

\subsubsection{The Rational Surface Class}\label{rcl:rational-surface-class}
Fix \(m\geq1\).  Recall that \(r_m\in H_c^0(R_m^\sm)\) is the standard class of
bivariant degree zero on the locus \(R_m^{\mathrm{sm}}\) of cycles
\(m[C]\), where \(C\) is a smooth embedded rational curve.  We first
construct an ambient lift of this class.

Put \(S=\mathbb P^1_C\times\mathbb P^1_W\) and
\(X_{\square}=\operatorname{Tot}_{S}
(\mathcal O_S(-1,-1)^{\oplus2})\), and let
\(\pi:S\to W\) be the second projection.  Since
\(\det(\mathcal O_S(-1,-1)^{\oplus2})=K_S\), the fourfold
\(X_{\square}\) is Calabi--Yau.  Let \(\beta\) be the class of a
fiber \(C_w=\mathbb P^1_C\times\{w\}\) of $\pi$.

\begin{lemma}
\label{lem:rm-compact-chow-representative}
The locus $Y_m^\square$ of compact cycles in \(X_{\square}\) of class \(m\beta\)
is an open and closed compact subspace of the cycle space, and
\begin{equation}
 Y_m^{\square}\cong\operatorname{Sym}^mW\cong\mathbb P^m,
 \qquad
 \Delta_m\coloneqq Y_m^{\square}\cap\mathcal \cZ_{\mathrm{con}}(X_\square)
 =\{m[C_w]:w\in W\}\cong W.
 \label{eq:rm-square-cycle-sector}
\end{equation}
Moreover, the restriction of the unit of \(Y_m^{\square}\) to the
connected locus is \(2r_m\).
\end{lemma}

\begin{proof}
Let \(p:X_{\square}\to S\) be the bundle projection.  If
\(D\subset X_{\square}\) is a compact irreducible curve whose class is
a positive multiple of \(\beta\), then \(p(D)\) is contained in a
fiber of \(\pi\).  Pulled-back to the normalization of \(D\), the two fiber
coordinates are sections of line bundles of negative degree.  These consequently
vanish, and hence \(D\) lies in the zero section. 

An effective curve
of class \(m\beta\) in \(S\) is a sum of fibers of \(\pi\).  Since the
homology class is locally constant on the cycle space, this proves
\eqref{eq:rm-square-cycle-sector}.

The relative normal bundle of the family \(C\times W=S\) over $B$ in the constant target
\(X_{\square}\) is
\begin{equation}
 N_{\square}=
 (\mathcal O_C\boxtimes\mathcal O_W(2))\oplus
 (\mathcal O_C(-1)\boxtimes\mathcal O_W(-1))^{\oplus2}.
 \label{eq:rm-square-normal-bundle}
\end{equation}
Write \(x=c_1(\mathcal O_C(1))\) and
\(h=c_1(\mathcal O_W(1))\), with
\(\int_{C\times W}xh=1\).  The Chern roots of
\(N_{\square}\) are \(2h,-x-h,-x-h\), so
\(c_2(N_{\square})=-2xh\).  The rational surface formula of
Proposition~\ref{rat:prop:rational-surface-reduction} therefore gives
\begin{equation}
 \mathbf1_{\Delta_m}
 =-\left(\int_{C\times W}c_2(N_{\square})\right)r_m
 =2r_m.
 \label{eq:rm-square-connected-restriction}
\end{equation}
\end{proof}

\subsubsection{}
\label{subsubsec:rm-is-a-lift}
Note now that $\Delta_m$ is compact, open, and closed in $\cZ_\con(X_\square)$, and is contained in the smooth rational locus $R_m^\sm(X_\square)$. Hence, $r_m=\frac{1}{2}\mathbf{1}_{\Delta_m}$ naturally extends to a class in 
\begin{equation*}
    r_m\in H_c^0(\cZ^{\sr\circ}_\con/\CY_4;\bQ).
\end{equation*}
We abuse notation slightly by denoting the classes in $H_c^0(R_m^\sm/\CY_4;\bQ)$, $H_c^0(\cZ^{\sr\circ}_\con/\CY_4;\bQ)$, and its image under extension by zero in $H_c^0(\cZ_\con/\CY_4;\bQ)$ all by $r_m$.

Because \(Y_m^{\square}\) is compact, open, and closed in the cycle
space, its unit defines a class in the augmentation ideal.  Set
\begin{equation}
 \widetilde r_m\coloneqq \frac12\mathbf1_{Y_m^{\square}}
 \in\mathfrak m^{[0]},
 \qquad
 [r_m]\coloneqq \widetilde r_m\bmod\mathfrak m^2
 \in(\mathfrak m/\mathfrak m^2)^{[0]}.
 \label{eq:rm-bracket-definition}
\end{equation}

\begin{proposition}\label{prop:rm-borel-lift}
One has \(i^*\widetilde r_m=r_m\).  In particular, \(r_m\in K^0\),
and \([r_m]\) is its unique ambient lift modulo \(\mathfrak m^2\).
\end{proposition}

\begin{proof}
The first assertion follows from
\eqref{eq:rm-square-connected-restriction}.  Since
\(\widetilde r_m\) is an ambient class, the localization sequence gives
\(\partial_{\mathrm{dis}}(r_m)=0\). 
\end{proof}

\subsubsection{The Rational Point Class}\label{subsec:rational-point-class}\label{subsubsec:definition-pm}
Fix \(m\geq1\), and use the same notation \(X_{\square}\),
\(Y_m^{\square}\cong\mathbb P^m\), and
\(\Delta_m\cong W\) as in
\S\ref{rcl:rational-surface-class} above.
Let \(H_m=c_1(\mathcal O_{\mathbb P^m}(1))\), and let
\(\nu_m:W\to\operatorname{Sym}^mW\), \(w\mapsto m[w]\), be the
small diagonal.  It is the Veronese embedding of degree \(m\), so
\(\nu_m^*H_m=m h_W\), where \(h_W=c_1(\mathcal O_W(1))\).

Recall from \S\ref{rnd:sec:rational-point-generator} that \(p_m\in H_c^2(R_m^{\mathrm{sm}}/\CY_4;\mathbb Q)\)
is the Poincar\'e dual of a point in a rational semi-chart.  On the
family \(\Delta_m\cong W\), it is represented by \(h_W\). As above, since $\Delta_m$ is open and closed in $\cZ_\con(X_\square)$, and contained in the smooth rational locus $R_m^\sm(X_\square)$, we can naturally extend to a class
\begin{equation*}
    p_m\in H_c^2(\cZ^{\sr\circ}_\con/\CY_4;\bQ).
\end{equation*}
By minor abuse of notation, we denote the classes in $H_c^2(R_m^\sm/\CY_4;\bQ)$, $H_c^2(\cZ^{\sr\circ}_\con/\CY_4;\bQ)$, and its image under extension by zero in $H_c^2(\cZ_\con/\CY_4;\bQ)$ all by $p_m$.

The description above in $X_\square$ additionally allows us to define
\begin{equation}
 \widetilde p_m\coloneqq \frac1mH_m\in\mathfrak m^{[2]},
 \qquad
 [p_m]\coloneqq \widetilde p_m\bmod\mathfrak m^2
 \in(\mathfrak m/\mathfrak m^2)^{[2]}.
 \label{eq:pm-bracket-definition}
\end{equation}
Here \(H_m\) is regarded as a class on the compact sector
\(Y_m^{\square}\) and extended by zero to the full cycle space.

\begin{proposition}
\label{prop:pm-borel-lift}
One has \(i^*\widetilde p_m=p_m\).  In particular, \(p_m\in K^2\),
and \([p_m]\) is its unique ambient lift modulo \(\mathfrak m^2\).
\end{proposition}

\begin{proof}
The restriction is
\(\nu_m^*(m^{-1}H_m)=h_B=p_m\). 
\end{proof}

\subsubsection{The Elliptic Point Classes}\label{subsec:elliptic-point-classes}
\begin{definition}
A \emph{regular CY elliptic normal datum} is a
triple
\begin{equation}
 (E,N,\theta),\qquad
 \theta:\det N\xrightarrow{\ \sim\ }K_E,\qquad H^1(E,N)=0,
 \label{eq:em-regular-elliptic-datum}
\end{equation}
where \(E\) is a smooth connected curve of genus one and \(N\) is a
holomorphic vector bundle of rank three.  Since \(\deg N=0\),
Riemann--Roch shows that regularity is equivalent to
\(H^0(E,N)=H^1(E,N)=0\). 

Let \(\mathscr N_1^{\mathrm{reg}}\) be the groupoid of these data,
with isomorphisms preserving \(\theta\), and put
\(\mathcal E=\pi_0(B\mathscr N_1^{\mathrm{reg}})\).
\end{definition}

First, we show that we may assume $N$ to be split for this data.

\begin{lemma}
\label{lem:em-regular-normal-form}
Every regular CY elliptic normal datum can be
deformed, through regular data, to one of the form
\begin{equation}
 N^\circ=L_1\oplus L_2\oplus L_3,\qquad
 L_i\in\operatorname{Pic}^0(E)\setminus\{\mathcal O_E\},\qquad
 L_1L_2L_3\cong K_E.
 \label{eq:em-split-normal-form}
\end{equation}
The line bundles may be chosen so that, after writing
\(L_3=(L_1L_2)^{-1}\), one has
\(H^*(E,L_1^aL_2^b)=0\) for every
\((a,b)\in\mathbb Z^2\setminus\{(0,0)\}\).  In fact,
\(B\mathscr N_1^{\mathrm{reg}}\) is connected.
\end{lemma}

\begin{proof}
The bundle \(N\) is semistable of degree zero.  Indeed, a subbundle of
positive degree has a nonzero section on an elliptic curve, contrary
to \(H^0(E,N)=0\).  By Atiyah's classification
\cite{Ati57}, every indecomposable summand of \(N\) is
\(F_r\otimes L\), where \(L\in\operatorname{Pic}^0(E)\) and \(F_r\)
has a filtration with quotients \(\mathcal O_E\).  The vanishing of
the cohomology of \(N\) implies \(L\neq\mathcal O_E\).

Scaling the extension classes in \(F_r\) deforms
\(F_r\otimes L\) to \(L^{\oplus r}\).  The cohomology remains zero,
and the determinant trivialisation can be transported along the
deformation: this gives \eqref{eq:em-split-normal-form}.  

For a fixed
elliptic curve, the resulting pairs \((L_1,L_2)\) form the connected complement
of the three loci \(L_1=\mathcal O_E\), \(L_2=\mathcal O_E\), and
\(L_1L_2=\mathcal O_E\) in \(\operatorname{Pic}^0(E)^2\). Remove the further loci
\(L_1^aL_2^b=\mathcal O_E\), for
\((a,b)\neq(0,0)\), and consider the resulting dense complement.

Any point of the
regular locus can be joined through other regular points to a point
in this dense complement.  Since the moduli space of elliptic curves and the
space of determinant trivialisations are connected as well, this proves the lemma.
\end{proof}

\subsubsection{}\label{subsubsec:definition-em}
Fix \(m\geq1\) and \((E,N)\in\mathcal E\).  Choose a general split
representative \(N^\circ=L_1\oplus L_2\oplus L_3\) as in
Lemma~\ref{lem:em-regular-normal-form}, and put
\(X_{\mathrm{ell}}=\operatorname{Tot}_E(N^\circ)\).

\begin{lemma}
\label{lem:em-cycle-sector-point}
The locus of compact effective cycles in \(X_{\mathrm{ell}}\) whose
pushforward to \(E\) has degree \(m\) consists of the single cycle
\(m[E]\).
\end{lemma}

\begin{proof}
Let \(D\subset X_{\mathrm{ell}}\) be a compact irreducible curve.
If its projection to \(E\) were constant, then \(D\) would lie in an
affine fiber of \(X_{\mathrm{ell}}\to E\), which is impossible.  After
normalization, its projection therefore gives a finite map
\(f:\widetilde D\to E\). 

Consider the three fiber coordinates are sections of
\(f^*L_i\).   If \(f^*L_i\) were trivial,
the norm map would imply that \(L_i^{\deg f}\) is trivial.  Hence
\(f^*L_i\) is a nontrivial line bundle of degree zero and has no
nonzero section.  In particular, all three coordinates vanish, so the reduced curve
is the zero section.  

The pushforward degree is locally constant on the
cycle space, so the asserted point is open and closed.
\end{proof}

\subsubsection{}\label{subsubsec:em-is-a-lift}
\begin{definition}
    The locus in Lemma~\ref{lem:em-cycle-sector-point} is an open and
    closed point of the cycle space, and contained in $E_m^\sm(X_\rmell)$. Hence, using the identity $\mathbf1_{\{m[E]\}}$ we can naturally extend $e_m$ of Definition~\ref{def:elliptic-standard-class-Em} to
    \begin{equation*}
        e_m\in H_c^0(\cZ_\con^{\sr\circ}/\CY_4).
    \end{equation*}  
    By slight abuse of notation, we again use the same symbol for the classes in $H_c^0(E_m^\sm/\CY_4)$, $H_c^0(\cZ_\con^{\sr\circ}/\CY_4)$, and its image in $H_c^0(\cZ_\con/\CY_4)$.
    The unit of the locus in Lemma~\ref{lem:em-cycle-sector-point} also defines an ambient class
    \(\widetilde e_m^\circ\coloneqq \mathbf1_{\{m[E]\}}\in\mathfrak m^{[0]}\).
    Transporting this class through a deformation supplied by
    Lemma~\ref{lem:em-regular-normal-form}, we obtain the ambient lift
    \begin{equation*}
        \tilde{e}_{m,(E,N)}\in \fm^{[0]}. 
    \end{equation*}  
    Finally, define
    \begin{equation}
     [e_{m,(E,N)}]
     \coloneqq \widetilde e_{m,(E,N)}\bmod\mathfrak m^2
     \in(\mathfrak m/\mathfrak m^2)^{[0]}.
     \label{eq:em-bracket-definition}
    \end{equation}
\end{definition}

\begin{proposition}
\label{prop:em-borel-lift}
One has
\(i^*\widetilde e_{m,(E,N)}=e_{m}\). In particular,
\(e_{m}\in K^0\), and \([e_{m,(E,N)}]\) is its unique ambient
lift modulo \(\mathfrak m^2\).
\end{proposition}

\begin{proof}
For the split representative, the whole sector of cycles and its connected
locus are both the point \(\{m[E]\}\), so the assertion is immediate.
It is preserved under the deformation to \((E,N)\).  
\end{proof}

\subsubsection{The Toric Ordinary Rational Meeting}
Recall the setup for the toric model $A\cup E\subset X_0=\Tot_{S_0}(L_1\oplus L_2)$ from \S\ref{inc:toric-model}. The map
\(\pi=p_2\circ\nu:S_0\to\mathbb P^1_t\) has general fiber \(F_t\) and
unique reducible fiber \(F_\infty=A+E\).  Thus,
\(A^2=E^2=-1\), \(A\cdot E=1\), and \([F_t]=[A]+[E]\). Let \(G=\mathbb C^*\) act trivially on \(S_0\), with weights \(1\) and
\(-1\) on \(L_1\) and \(L_2\). This action preserves the volume form,
and \(X_0^G=S_0\). 

\begin{definition}
    Fixing \(d_1,d_2\geq1\), and writing \(k=\min(d_1,d_2)\), we consider the class
    \[\beta_{d_1,d_2}=d_1[A]+d_2[E].\]
    Consider the projective completion
    \(\overline X_0=\mathbb P_{S_0}(\mathcal O\oplus L_1\oplus L_2)\),
    with divisor at infinity
    \(D_\infty=\mathbb P_{S_0}(L_1\oplus L_2)\), and write
    \(j:X_0\hookrightarrow\overline X_0\).  Choose a very ample divisor
    \(H\) on \(\overline X_0\).  In the projective Chow variety of
    \(\overline X_0\), take the union \(\overline Y_{d_1,d_2}\) of the
    components parametrising cycles of homology class
    \(j_*\beta_{d_1,d_2}\) and degree
    \(H\mathbin{\cdot}j_*\beta_{d_1,d_2}\), and put
    \[
     Y_{d_1,d_2}
     \coloneqq \overline Y_{d_1,d_2}
     \setminus\{\text{cycles meeting }D_\infty\}.
    \]
    Its analytification is the full locus of compact analytic cycles in
    \(X_0\) of class \(\beta_{d_1,d_2}\).
\end{definition}

\subsubsection{}
\label{subsubsec:toric-cap-is-a-lift}

\begin{proposition}
\label{prop:cap-connected-cycle-sector}
The connected locus in \(Y_{d_1,d_2}\), with its reduced structure, is
\begin{equation}
 C_{d_1,d_2}
 \coloneqq Y_{d_1,d_2}\cap\mathcal Z_{\mathrm{con}}(X_0)
 =\begin{cases}
   \{z_{d_1,d_2}\},&d_1\neq d_2,\\
   \Gamma_d,&d_1=d_2=d,
  \end{cases}
 \label{eq:cap-entire-connected-sector}
\end{equation}
where $z_{d_1,d_2}\coloneqq d_1A+d_2E$ is the ordinary rational meeting cycle, and \(\Gamma_d=\{dF_t:t\in\mathbb P^1\}\cong\mathbb P^1\).
\end{proposition}

\begin{proof}
Let \(p:X_0\to S_0\) be the bundle projection and put
\(\varpi=\pi\circ p\).  Every irreducible component of an effective
cycle has nonnegative degree against
\(\varpi^*\mathcal O_{\mathbb P^1}(1)\), while this degree is zero on
\(\beta_{d_1,d_2}\).  Thus every component is contained in a fiber of
\(\varpi\).

For \(t\neq\infty\), this fiber is
\(\operatorname{Tot}_{F_t}
(\mathcal O_{F_t}(-1)\oplus\mathcal O_{F_t}(-1))\).
The normalization argument used in
Lemma~\ref{lem:rm-compact-chow-representative} shows that its only
compact irreducible curve is the zero section.  A connected cycle over
such a fiber therefore has class \(d[A]+d[E]=d[F_t]\).

On the special fiber, the restrictions are
\begin{equation}
 (L_1,L_2)|_A=(\mathcal O_A(-1),\mathcal O_A),\qquad
 (L_1,L_2)|_E=(\mathcal O_E,\mathcal O_E(-1)).
 \label{eq:cap-special-fiber-restrictions}
\end{equation}
The compact irreducible curves are the constant translates \(A_a\)
in the trivial direction over \(A\), and \(E_b\) in the trivial
direction over \(E\).  They meet if and only if \(a=b=0\).  Since
\(d_1,d_2>0\), a connected cycle in the special fiber must contain
both kinds of component.  It is therefore \(d_1A+d_2E\).  This proves
\eqref{eq:cap-entire-connected-sector}.
\end{proof}

This allows us to define
\[
 \ell_{d_1,d_2}^{\mathrm{con}}
 \coloneqq \begin{cases}
   \mathbf1_{\{z_{d_1,d_2}\}},&d_1\neq d_2,\\
   \mathbf1_{\Gamma_d},&d_1=d_2=d
  \end{cases}
\]
in $H_c^0(\cZ^{\sr\circ}_\con/\CY_4;\bQ)$, which extends to a class in $H_c^0(\cZ_\con/\CY_4,\bQ)$ that we denote by the same notation. Note that, in the first case, this is exactly the point Thom representative of the ordinary rational meeting cycle in the toric model from \S\ref{inc:toric-model}.

\subsubsection{}
Since $Y_{d_1,d_2}$ is not compact, we instead work with the fixed locus and a Borel construction to find an ambient lift.
\begin{proposition}
\label{prop:cap-fixed-chow-family}
As a reduced space, the fixed locus is
\begin{equation}
 \begin{aligned}
 K_{d_1,d_2}&\coloneqq (Y_{d_1,d_2}^G)_{\mathrm{red}}
 =\left\{(d_1-k)A+(d_2-k)E+\sum_{i=1}^kF_{t_i}\right\},\\
 K_{d_1,d_2}&\cong\operatorname{Sym}^k\mathbb P^1\cong\mathbb P^k.
 \end{aligned}
 \label{eq:cap-fixed-chow-family}
\end{equation}
The point corresponding to \(k[\infty]\) is
\(z_{d_1,d_2}\coloneqq d_1A+d_2E\).
\end{proposition}

\begin{proof}
Every irreducible component of a fixed cycle is preserved by \(G\).
If \(G\) acted nontrivially on its normalization, its general orbit
would be dense.  Since the bundle projection \(X_0\to S_0\) is
\(G\)-invariant, the curve would lie in an affine fiber of this
projection, which is impossible for a compact curve.  Thus every fixed
cycle lies in \(S_0\).

Let \(F=[A]+[E]\).  Since \(F\) is nef and
\(F\cdot A=F\cdot E=0\), every component of a cycle in \(S_0\) of the
chosen class is contained in a fiber of \(\pi\).  Such a cycle has a
unique expression \(aA+bE+\sum_{t\neq\infty}n_tF_t\).  If
\(r=\sum n_t\), comparison of the coefficients of \([A]\) and \([E]\)
gives \(a=d_1-r\), \(b=d_2-r\), and \(0\leq r\leq k\).  It is therefore
determined by the divisor
\((k-r)[\infty]+\sum n_t[t]\) of degree \(k\).  This proves
\eqref{eq:cap-fixed-chow-family}.  
\end{proof}

\subsubsection{Borel Approximation}\label{rcl:sec:ell-borel-approx}
The unit of \(K_{d_1,d_2}\) is not a compactly supported class on the
generally noncompact space \(Y_{d_1,d_2}\).  Moreover,
\(K_{d_1,d_2}\) is closed rather than open, so its unit cannot be
extended by zero.  We use the following standard consequence of
localization.

\begin{lemma}
\label{lem:cap-finite-borel-clearing}
Let \(Y\) be a complex algebraic \(G=\mathbb C^*\)-space having the
homotopy type of a finite \(G\)-CW complex.  Suppose that
\(K=Y^G\) is compact and fixed pointwise.  For some \(r\geq0\) and
\(N\geq r\), there is a class
\begin{equation}
 \alpha_N\in H_c^{2N}(Y_N;\mathbb Q)
 =H_c^0(Y_N/\mathbb P^N;\mathbb Q),\qquad
 \alpha_N|_{K\times\mathbb P^N}=h^N\mathbf1_K,
 \label{eq:cap-cleared-borel-class}
\end{equation}
where \(Y_N=Y\times_G(\mathbb C^{N+1}\setminus\{0\})\) and
\(h=c_1(\mathcal O_{\mathbb P^N}(1))\).
\end{lemma}

\begin{proof}
Localization for equivariant cohomology with compact supports gives
\[
 H^*_{G,c}(Y;\mathbb Q)[u^{-1}]
 \xrightarrow{\ \sim\ }
 H^*_{G,c}(K;\mathbb Q)[u^{-1}],\qquad |u|=2.
\]
After clearing a denominator, we obtain
\(\gamma\in H^{2r}_{G,c}(Y;\mathbb Q)\) whose restriction to \(K\) is
\(u^r\mathbf1_K\).  Since
\(H^*_{G,c}(K;\mathbb Q)=H_c^*(K;\mathbb Q)\otimes\mathbb Q[u]\) has
no \(u\)-torsion, this equality holds before localization.

Since $\bC^{N+1}\setminus\{0\}\simeq S^{2N+1}$ is $2N$-connected and $Y$ has the homotopy type of a finite $G$-CW complex, $Y_N$ approximates the Borel construction, so that $H^{2r}_{G,c}(Y;\bQ)\cong H^{2r}_c(Y_N;\bQ)$ for sufficiently large $N$. Fix such an $N\geq r$, let $\gamma_N\in H_c^{2r}(Y_N;\bQ)$ correspond to $\gamma$, and let
\(p_N:Y_N\to\mathbb P^N\) be the projection.  Then
\(\alpha_N=p_N^*h^{N-r}\smile\gamma_N\) has ordinary degree \(2N\)
and hence bivariant degree zero.  Its restriction is
\(h^N\mathbf1_K\).
\end{proof}

\begin{definition}
    Apply the lemma to \(Y_{d_1,d_2}\) and \(K_{d_1,d_2}\).  Define
    \[
     Y_{d_1,d_2,N}
     \coloneqq 
     Y_{d_1,d_2}\times_G
     (\mathbb C^{N+1}\setminus\{0\}).
    \]
    This is the open and closed locus in
    \(\mathcal Z((X_0)_N/\mathbb P^N)\) consisting of cycles whose
    fiberwise homology class is \(\beta_{d_1,d_2}\).  Choose a class
    \(\alpha_{d_1,d_2,N}\) supplied by
    Lemma~\ref{lem:cap-finite-borel-clearing}.  Extension by zero from the open and closed
    locus \(Y_{d_1,d_2,N}\) to $\cZ((X_0)_N/\bP^N)$ gives a class whose
    image in the universal group defines
    \begin{equation}\label{rcl:eq:lifted-rational-meeting-class}
        \widetilde\ell_{d_1,d_2}\in\mathfrak m^{[0]},
    \end{equation}
    and
    \begin{equation}
        [\ell_{d_1,d_2}]\coloneqq \widetilde\ell_{d_1,d_2}\bmod\mathfrak m^2\in (\fm/\fm^2)^{[0]}.
    \end{equation}
\end{definition}

\subsubsection{}

\begin{proposition}
\label{prop:cap-is-a-lift}
One has
\(i^*\widetilde\ell_{d_1,d_2}
=\ell_{d_1,d_2}^{\mathrm{con}}\).  Consequently
\(\ell_{d_1,d_2}^{\mathrm{con}}\in K^0\), and
\([\ell_{d_1,d_2}]\) is its unique ambient lift modulo
\(\mathfrak m^2\).
\end{proposition}

\begin{proof}
The connected locus of
Proposition~\ref{prop:cap-connected-cycle-sector} is contained in
\(K_{d_1,d_2}\).  At the finite approximation its inverse image is
\(C_{d_1,d_2}\times\mathbb P^N\), and the class supplied by
Lemma~\ref{lem:cap-finite-borel-clearing} restricts to
\(h^N\mathbf1_{C_{d_1,d_2}}\).  Since
\(\int_{\mathbb P^N}h^N=1\), this is the relative unit.  The rest is clear.
\end{proof}

\subsection{Generators and the Logarithmic Reduction}
\label{rcl:subsec:generators-logarithmic-reduction}

We now prove that the standard connected classes constructed in
\S\ref{subsec:standard-connected-classes} generate the universal
cohomology of the connected locus in degrees zero, one, and two.
We use the multiplicity filtration to compare classes on the exact
multiplicity loci and then subtract the restrictions of their standard
ambient lifts. In this entire Section~\ref{rcl:subsec:generators-logarithmic-reduction} we work with rational coefficients, which we suppress from the notation for brevity in many instances. The existence of the universal bivariant groups we work with here, can be shown using Proposition~\ref{rcl:prop:cohomology-fiberwise-cond} for all loci considered in this section.

\subsubsection{Multiplicity and Localization}
\label{rcl:subsubsec:multiplicity-localization}

For a downward closed subset $I$ of the multiplicity poset of \S\ref{subsec:multiplicity-filtration},
write
\[
 C_I\coloneqq \mathcal Z_{\mathrm{con}}\cap\mathcal Z_I,
 \qquad
 C_{=\mathbf m}\coloneqq \mathcal Z_{\mathrm{con}}\cap
 \mathcal Z_{=\mathbf m},
\]
and $C^{\sr\circ}_I$ and $C_{=\bfm}^{\sr\circ}$ for the analogous intersections with $\cZ_\con^{\sr\circ}$.
The loci used below commute
with base change and restriction to open subfamilies. Hence their
localization sequences remain exact after passage to the filtered
colimit. If $I$ is finite and $\mathbf m$ is maximal in $I$, then
$I'\coloneqq I\setminus\{\mathbf m\}$ is downward closed,
$C_{I'}\subset C_I$ is open, since its closed complement is
$C_{=\mathbf m}$. In particular, localization \eqref{biv:localization} gives
\begin{equation}
 \cdots \to H_c^q(C_{I'}/\CY_4)\to H_c^q(C_I/\CY_4)
 \xrightarrow{i_{=\bfm}^*}H_c^q(C_{=\mathbf m}/\CY_4)\to\cdots.
 \label{rcl:eq:multiplicity-localization}
\end{equation}

We shall also use the following version of the dimension argument of
Proposition~\ref{inc:prop:cohomological-range} in degree $-1$.

\begin{proposition}\label{rcl:prop:cohom-range-1}
    Let $Y$ be a locally closed sublocus in $\cZ^{\sr\circ}_\con$, compatible with base change and local in the sense of \cite[Def. 3.13]{pardon2023}, and let $k$ be a fixed integer. If every $f$ in every $Y(X/B)$ has $\nu(f)\leq k$, then
\begin{equation}
 H_c^j(Y/\CY_4;\bQ)=0\qquad\text{for }j\geq 2k+1.
 \label{rcl:eq:negative-incidence-vanishing}
\end{equation}
\end{proposition}

\begin{proof}
    Represent a compactly supported class as a class $\alpha\in H_c^*(Y(X/B)/B)$ over a family \(X\to B\), and choose the finite Whitney stratification of Lemma~\ref{inc:lem:finite-incidence-stratification} near its support. Write $d_\sigma\coloneqq \dim_{\mathbb R}\sigma$. A stratum $T$ over the interior $\sigma^\circ$ of a real simplex $\sigma$ has real dimension $d_\sigma$ plus twice the complex index by semi-regularity, which is bounded above by $d_\sigma+2k$ by assumption, using that the strata have exact values of the complex index $\nu(f)$. The dualizing complex of the oriented real manifold \(\sigma^\circ\) gives, up to its orientation local system, \(H_c^j(T/\sigma^\circ)\cong H_c^{j+d_\sigma}(T)\). By the dimension bound above, the right-hand side vanishes for $j\geq 2k+1$. Successive localization using \eqref{biv:localization} gives $H_c^j(Y(X/B)/B)=0$ for $j\geq 2k+1$, and hence we obtain \eqref{rcl:eq:negative-incidence-vanishing} in the colimit.
\end{proof}

\subsubsection{}
\begin{lemma}\label{rcl:lemma:cohomology-identificaiton-localization}
    Let $j:U\hookrightarrow Y$ be an open subset with closed complement $i:K\hookrightarrow Y$. Assume both are compatible with base change and local in the sense of \cite[Def. 3.13]{pardon2023}.
    \begin{itemize}
        \item If $K$ satisfies $H_c^j(K/\CY_4)=0$ for $j\geq -1$, then extension by zero induces $H_c^q(U/\CY_4)\simeq H_c^q(Y/\CY_4)$ for $q\geq0$.
        \item If instead $U$ satisfies $H_c^j(U/\CY_4)=0$ for $j\geq -1$, then restriction induces $H_c^q(Y/\CY_4)\simeq H_c^q(K/\CY_4)$ for $q\geq -1$.
    \end{itemize}
\end{lemma}

\begin{proof}
    This immediately follows from the localization sequence 
    \begin{equation*}
        \cdots\to H_c^{j-1}(K/\CY_4)\to H_c^j(U/\CY_4)\xrightarrow{j_!} H_c^j(Y/\CY_4)\xrightarrow{i^*} H_c^j(K/\CY_4)\to H_c^{j+1}(U/\CY_4)\to\cdots
    \end{equation*}
    induced by \eqref{biv:localization}.
\end{proof}

\subsubsection{Generation}
\label{rcl:subsubsec:connected-generation}

\begin{theorem}
\label{rcl:thm:indecomposable-quotient-generators}
In bivariant degrees $q\in\{0,1,2\}$, restriction induces an
isomorphism
\begin{equation}
 (\mathfrak m/\mathfrak m^2)^{[q]}
 \xrightarrow{\ \sim\ }
 H_c^q(\mathcal Z_{\mathrm{con}}/\CY_4;\mathbb Q).
 \label{rcl:eq:indecomposables-connected-isomorphism}
\end{equation}
The part of degree
two is spanned by the $[p_m]$, the part of degree one vanishes, and
the part of degree zero is spanned by
\[
 [r_m],\qquad [e_{m,(E,N)}],\qquad [\ell_{d_1,d_2}],
 \qquad m,d_1,d_2\geq1,\quad (E,N)\in\mathcal E.
\]
Under \eqref{rcl:eq:indecomposables-connected-isomorphism}, these
classes correspond to $p_m,r_m,e_{m}$, and
$\ell_{d_1,d_2}^{\mathrm{con}}$, respectively.
\end{theorem}

\begin{proof}
 By Theorem~\ref{gt:thm:pardon-ring-sro-isom}, we may work entirely in the interior semi-regular locus $\cZ_\con^{\sr\circ}$, see Step~4 for details. 
 
 Below we first compute the cohomology of the loci where the multiplicity stays constant, and then apply \eqref{rcl:eq:multiplicity-localization} and \S\ref{uce:sec:total-support-multiplicity}. 
\medskip\noindent
\paragraph{Step 1. Profiles with one component.}
Fix the profile $[m]$. The Euler characteristic
$\widetilde\chi$ of normalization is upper semicontinuous on this locus by Lemma~\ref{inc:lem:semi-chart-stratification}, so
$Y_2\coloneqq C^{\sr\circ}_{=[m],\widetilde\chi=2}$ is closed in $C^{\sr\circ}_{=[m]}$, as $2$ is the maximal value of $\widetilde{\chi}$.
The smooth rational locus $R_m^{\mathrm{sm}}$ is open in $Y_2$:
in the semi-charts the normalization is a smooth rational curve,
and being an embedding is an open condition. Every $f$ in the complement over some family $X/B$ has
negative index by \eqref{inc:eq:connected-support-index}. By Proposition~\ref{rcl:prop:cohom-range-1} and Lemma~\ref{rcl:lemma:cohomology-identificaiton-localization} above, we hence get
\begin{equation*}
    H_c^{\geq 0}(R_m^{\mathrm{sm}}/\CY_4)\cong H_c^{\geq 0}(Y_2/\CY_4)
\end{equation*}
by extension by zero. Now consider the open complement $U_0\coloneqq C^{\sr\circ}_{=[m],\widetilde\chi\leq0}$ of $Y_2$.

The locus $Y_0\coloneqq C_{=[m],\widetilde\chi=0}^{\sr\circ}$ is closed in $U_0$ by upper semicontinuity.
Every $f$ in its open complement $C_{=[m],\widetilde\chi<0}^{\sr\circ}$ over some family $X/B$ has normalization genus at least two and hence
negative index by \eqref{inc:eq:connected-support-index}. By Proposition~\ref{rcl:prop:cohom-range-1} (for $k=-1$) and Lemma~\ref{rcl:lemma:cohomology-identificaiton-localization} above, we hence get
\begin{equation*}
    H_c^{\geq 0}(U_0/\CY_4)\cong H_c^{\geq 0}(Y_0/\CY_4)
\end{equation*}
via closed restriction. The smooth elliptic locus $E^\sm_m$ is open in $Y_0$, and its closed complement of singular cycles again has negative index by \eqref{inc:eq:connected-support-index}, which gives us
\begin{equation*}
    H_c^{\geq 0}(E_m^{\mathrm{sm}}/\CY_4)\cong H_c^{\geq 0}(Y_0/\CY_4)
\end{equation*}
via extension by zero, using Proposition~\ref{rcl:prop:cohom-range-1} (for $k=-1$) and Lemma~\ref{rcl:lemma:cohomology-identificaiton-localization}. We can collect these identifications in the localization sequence \eqref{biv:localization} to get
\begin{equation}\label{rcl:eq:m-localization-sequence}
    \begin{tikzcd}[column sep=small,row sep=small]
    \cdots\ar[r] & H_c^{j\geq 0}(U_0/\CY_4)\ar[r]\ar[d,"\cong"] & H_c^{j\geq 0}(C_{=[m]}^{\sr\circ}/\CY_4)\ar[r] & H_c^{j\geq 0}(Y_2/\CY_4)\ar[r]\ar[d,"\cong"] & \cdots\\
    & H_c^{j\geq 0}(E_m^{\mathrm{sm}}/\CY_4) & & H_c^{j\geq 0}(R_m^{\mathrm{sm}}/\CY_4). &
    \end{tikzcd}
\end{equation}
The rational normal data map of Proposition~\ref{prop:rational-normal-data-bordism} and the calculation in \S\ref{rat:sec:degree-two-normal-generator}
show that $H_c^{j\geq 0}(R_m^{\mathrm{sm}}/\CY_4)$ are spanned by $p_m$ in degree two, by
$r_m$ in degree zero, and vanish in degree one. 

Smooth elliptic curves have index $0$ by \eqref{inc:eq:connected-support-index}, so Proposition~\ref{rcl:prop:cohom-range-1} for $k=0$ shows that $H_c^{*}(E_m^\sm/\CY_4)$ vanishes in positive degrees. By semi-regularity, Poincar\'e duality in the smooth semi-charts, shows that $H_c^{0}(E_m^\sm/\CY_4)$ is generated by Thom representatives of points.  The elliptic comparison in \S\ref{sec:elliptic-locus} and \S\ref{subsec:elliptic-point-classes} identifies these with the standard elliptic point classes $e_{m,(E,N)}$.

It remains to check that, under the maps in the above localization sequence \eqref{rcl:eq:m-localization-sequence}, these generators coincide with the restrictions of $p_m$, $r_m$, and $e_{m,(E,N)}$ to $C_{=[m]}^{\sr\circ}$. More precisely, the standard classes $r_m$, $p_m$, and $e_{m,(E,N)}$ in $H_c^*(\cZ_\con^{\sr\circ}/\CY_4)$ are extended from $H_c^*(C_{\leq[m]}^{\sr\circ}/\CY_4)$ and we compare with their restrictions to $H_c^{*}(C_{=[m]}^{\sr\circ}/\CY_4)$. By definition of the groups as colimits, it suffices to show these identifications on specific geometries.

To compare the rational generators, consider the fourfold $X_\square$ of \S\ref{rcl:rational-surface-class}. The class $r_m$ is constructed in \S\ref{subsubsec:rm-is-a-lift} in $H_c^0(\cZ_\con^{\sr\circ}(X_\square))$ by extension from the open and closed sublocus $\Delta_m$ of $\cZ_\con^{\sr\circ}(X_\square)$. This open and closed locus already lies inside $R_m^\sm(X_\square)$, and hence, the restriction to $C_{=[m]}^{\sr\circ}$ of the corresponding class $r_m$ maps exactly to $r_m\in H_c^0(R_m^\sm/\CY_4)$ under the second map in \eqref{rcl:eq:m-localization-sequence}. The comparison of the class $p_m$ follows analogously.

To compare the elliptic generators, consider the fourfold $X_\rmell$ for a fixed split choice $N^\circ$ as in \S\ref{subsubsec:definition-em}. The class $e_m\in H_c^0(\cZ_\con^{\sr\circ}(X_\rmell))$ is extended from the open and closed sublocus of Lemma~\ref{lem:em-cycle-sector-point}, which already lies inside $E_m^\sm(X_\rmell)$, and, hence, the restriction to $C_{=[m]}^{\sr\circ}$ of the $e_m$ matches the image of $e_m\in H_c^0(R_m^\sm/\CY_4)$ under the first map in \eqref{rcl:eq:m-localization-sequence}. In sum, it follows that $H_c^2(C_{=[m]}^{\sr\circ})$ is spanned by the restriction of
$p_m$, that $H_c^1(C_{=[m]}^{\sr\circ})=0$, and that $H_c^0(C_{=[m]}^{\sr\circ})$ is
spanned by the restrictions of $r_m$ and the $e_{m}$.

\medskip\noindent
\paragraph{Step 2. Profiles with at least two components.}
Fix $\mathbf m=[d_1,d_2]$. Its maximal Euler
characteristic of normalization is $4$, and hence
$Y_4\coloneqq C_{=\mathbf m,\widetilde\chi=4}^{\sr\circ}$ is closed in
$C_{=\mathbf m}$ by upper semicontinuity. For every $f$ in the open complement $C_{=\mathbf m,\widetilde\chi<4}^{\sr\circ}$ over some family $X/B$, the sum of the
 genera of normalizations is positive, so Lemma~\ref{inc:lem:connected-support-index} gives negative index, so that
\begin{equation*}
    H_c^{\geq 0}(C_{=\bfm}^{\sr\circ}/\CY_4)\cong H_c^{\geq 0}(C_{=\mathbf m,\widetilde\chi=4}^{\sr\circ}/\CY_4)
\end{equation*}
via restriction, by Proposition~\ref{rcl:prop:cohom-range-1} and Lemma~\ref{rcl:lemma:cohomology-identificaiton-localization}.
Within $Y_4$, let $W$ be the open locus where both component maps
are embeddings. Since the meeting already cancels the positive $(1-g_i)$-terms of the rational components in \eqref{inc:eq:connected-support-index}, the complement of $W$, where the components may have singular points, also has negative index. Applying Proposition~\ref{rcl:prop:cohom-range-1} (for $k=-1$) and Lemma~\ref{rcl:lemma:cohomology-identificaiton-localization} twice, we therefore can identify
\begin{equation*}
    H_c^{\geq 0}(C_{=\mathbf m,\widetilde\chi=4}^{\sr\circ}/\CY_4)\cong H_c^{\geq 0}(Y_4/\CY_4)\cong H_c^{\geq 0}(W/\CY_4).
\end{equation*}
Let $O\subset W$ be the locus where there is a meeting at exactly one point, and that meeting point has distinct tangent lines, i.e. is transversal. The locus of transversal meetings is open in $W$, and since the components are proper, transversality of meetings implies that the locus $O$ where there is exactly one such meeting is open. Note that multiple meetings have relative complex index at most \(-2\), and non-transversal meetings have positive order, so that the index of \eqref{inc:eq:connected-support-index} is $\nu\leq -4$. Hence, we may identify
\begin{equation*}
    H_c^{\geq 0}(O/\CY_4)\cong H_c^{\geq 0}(W/\CY_4)
\end{equation*}
via extension by zero, by applying Proposition~\ref{rcl:prop:cohom-range-1} and Lemma~\ref{rcl:lemma:cohomology-identificaiton-localization}.

The ordinary meeting spaces $O$ have relative complex index zero.
By Proposition~\ref{rcl:prop:cohom-range-1} for $k=0$, their positive degree cohomology vanishes. By semi-regularity, Poincar\'e duality in the smooth semi-charts, shows that $H_c^{0}(O/\CY_4)$ is generated by point Thom representatives. The comparison with
the toric model identifies these with the standard toric ordinary meeting
class, which Proposition~\ref{prop:cap-is-a-lift} identifies with the restriction of $\ell_{d_1,d_2}^{\mathrm{con}}$ to the
exact multiplicity locus. In the geometry used to define $\ell_{d_1,d_2}^{\mathrm{con}}$ we can see that, if $d_1=d_2=d$, the other connected
members of its pencil have profile $[d]$ and do not occur in
$C_{=[d,d]}^{\sr\circ}$. Hence,
$H_c^0(C_{=[d_1,d_2]}^{\sr\circ})$ is spanned by the restriction of
$\ell_{d_1,d_2}^{\mathrm{con}}$, and its positive degree cohomology
vanishes.

For a multiplicity $\bfm$ with $\ell\geq3$ components, the estimate \eqref{inc:eq:connected-support-bound} is
$\nu\leq2-\ell-\sum_i g_i\leq-1$, so its nonnegative degree cohomology
therefore vanishes by Proposition~\ref{rcl:prop:cohom-range-1} for $k=-1$.

\medskip\noindent
\paragraph{Step 3. Induction.}
The generators in Steps~1--2 are restrictions of standard connected
classes represented on $C^{\sr\circ}_{\leq\mathbf m}$. For the rational and
elliptic classes, their models have connected support of the single
profile $\bfm=[m]$. For the meeting classes, by
Proposition~\ref{prop:cap-connected-cycle-sector}, their connected
support multiplicity is $\bfm=[d_1,d_2]$ when $d_1\neq d_2$, and $\bfm=[d]$ when
$d_1=d_2=d$. Since $[d]<[d,d]$, this representative also lies
in $C_{\leq[d,d]}^{\sr\circ}$. Via extension by zero, we can extend these to $C_I^{\sr\circ}$ of any downward closed subset $I$ containing $\bfm$.

Let $I$ be a finite downward closed set of nonempty multiplicities,
$\mathbf m\in I$ a maximal element, and
$I'=I\setminus\{\mathbf m\}$. For $u\in H_c^q(C_I^{\sr\circ})$, where
$q\in\{0,1,2\}$, Steps~1--2 give an expression
$i^*_{=\mathbf m}(u)=\sum_j b_ji^*_{=\mathbf m}(s_j)$, where
each $s_j$ is a standard connected representative on
$C_{\leq\mathbf m}^{\sr\circ}$. Using upper semicontinuity and that $I$ is downward closed, we have an open
inclusion $C_{\leq\mathbf m}^{\sr\circ}\subset C_I^{\sr\circ}$. Hence, we can extend the $s_j$ to $C_I^{\sr\circ}$.
Exactness of \eqref{rcl:eq:multiplicity-localization} implies that
$u-\sum_jb_js_j$ comes from $H_c^q(C_{I'}^{\sr\circ})$.
Induction on $|I|$ proves generation by the standard connected
classes in degrees zero and two, and vanishing in degree one.

By \eqref{uce:eq:H-total-mult-exhaustion}, every compactly supported class is represented on an
open locus where the total support multiplicity $\lambda$ is
bounded by some integer $N$. The nonempty profiles with
$\lambda\leq N$ form a finite downward closed set. Applying the
preceding induction proves the asserted generation on the whole
connected locus. 

Every standard connected class has an ambient lift by
\S\ref{subsec:standard-connected-classes}. Therefore, the proven generation
implies $H_c^q(\mathcal Z^{\sr\circ}_{\mathrm{con}}/\CY_4;\mathbb Q)\cong K^q$
for $q=0,1,2$. Proposition~\ref{prop:disconnected-square} and Theorem~\ref{gt:thm:pardon-ring-sro-isom} now give us the commutative diagram
\begin{equation*}
    \begin{tikzcd}
        & H_c^q(\cZ_\con^{\sr\circ}/\CY_4)\ar[d,"\cong"]\\
        (\fm/\fm^2)^{[q]}\cong K^q \ar[r,hookrightarrow]\ar[ur,"\cong"] & H_c^q(\cZ_\con/\CY_4),
    \end{tikzcd}
\end{equation*}
which implies \eqref{rcl:eq:indecomposables-connected-isomorphism}, and Proposition~\ref{prop:disconnected-square} identifies
the generators with the stated classes modulo $\mathfrak m^2$.
\end{proof}

\subsubsection{The Logarithmic Reduction}
\label{rcl:subsubsec:logarithmic-reduction}

\begin{corollary}
\label{rcl:cor:equivalence-num-checks}
The following are equivalent:
\begin{enumerate}[label=(\alph*)]
\item $\delta_1([p_m])=0$, and
$\delta_0([r_m])=\delta_0([e_{m,(E,N)}])
=\delta_0([\ell_{d_1,d_2}])=0$ for every index occurring in that range;
\item $L_{\Phi,1}=L_{\Psi,1}$ in bivariant degree two, and
$L_{\Phi,0}=L_{\Psi,0}$ in bivariant degree zero.
\end{enumerate}
There is no term of bivariant degree one.
\end{corollary}

\begin{proof}
Every $\delta_n$ vanishes on $\mathfrak m^2$ by
\eqref{eq:coefficient-discrepancy-indecomposables}, and hence
factors through $\mathfrak m/\mathfrak m^2$.
Theorem~\ref{rcl:thm:indecomposable-quotient-generators} shows that
$\delta_1$ vanishes in degree two precisely when it vanishes on
the $[p_m]$, and that $\delta_0$ vanishes in degree zero precisely
when it vanishes on the $[r_m]$, $[e_{m,(E,N)}]$, and
$[\ell_{d_1,d_2}]$. This proves the equivalence of (a) and (b).
The assertion in degree one follows from the same theorem.
\end{proof}

\section{Numerical Checks \& Applications}\label{liftschecks}
\subsection{Numerical Computations}
Now, we want to check the numerical calculation in the equivalence of Corollary~\ref{rcl:cor:equivalence-num-checks} above. Combining that corollary with Propositions~\ref{prop:gwpt-rm}, \ref{prop:gwpt-pm}, \ref{prop:gwpt-em}, and \ref{prop:gwpt-toric-cap} below then finally proves:
\begin{theorem}\label{nap:thm:L-phi-psi-correspondence}
     We have $L_{\Phi,1}=L_{\Psi,1}$ in bivariant degree two, and $L_{\Phi,0}=L_{\Psi,0}$ in bivariant degree zero, and no term of bivariant degree one.
\end{theorem}

\subsubsection{Orientations}
In the following paragraphs we reduce calculations of certain explicit family $GW$- and $PT$-invariants to those from \cite{MCT}\footnote{There is a discrepancy between the numbers of statements in the latest arXiv version of this paper, and those in the published version. We follow the latter.}.

We recall, however, that \cite{MCT} prove all their theorems modulo a sign: the issue they do not know which one comes from the canonical orientation choice. We now explain that the canonical orientations of \cite{karpov-cy4-orientations} yields the correct signs for the formulas in \cite{MCT}. Indeed, in all situations below when the relevant invariants are non-zero (and, hence, the choice of signs matters at all), the relevant geometries are \textit{local curves.} Hence, canonical orientations (in the sense of \cite{JU} and \cite{Bojko2021}) coincide with the explicit Shi's orientations for local CY $4$-folds, see \cite[Theorem 1.2]{karpov-cy4-orientations}. It is an easy exercise to compute the corresponding signs explicitly, and see that in all examples below they are, indeed, those which give the ``right'' answer in \cite{MCT}.

\subsubsection{}
Before starting with the actual calculations, we state here a useful Lemma, which we repeatedly use.
\begin{lemma}\label{num:lemma:PT-vanishing-global-gen-C-lemma}
Every proper $1$-dimensional subscheme $Z$ of a variety $X$ whose reduced support $C$ has globally generated conormal sheaf and $H^1(X,\cO_C)=0$ satisfies $\chi(\cO_Z)\geq 1$. 

Consequently, if $\beta\in H_2(X,\bZ)$ is a curve class, such that every proper $1$-dimensional subscheme $Z$ of $X$ with $[Z]=\beta$ has reduced support $C$ with globally generated conormal sheaf and $H^1(X,\cO_C)=0$, then $P_{\leq 0}(X,\beta)=\emptyset$.
\end{lemma}
\begin{proof}
     Let $C=Z_\red$ be the reduced support of $Z$ and $N^\vee_C=I_C/I_C^2$ its globally generated conormal sheaf. The filtration of $\cO_Z$ by powers of the ideal of $C$ has degree \(a\) part $I_C^a\cO_Z/I_C^{a+1}\cO_Z$. The ideal of $C$ in $Z$ is $I_C\cO_Z$, and hence by \cite[\href{https://stacks.math.columbia.edu/tag/062Z}{Tag 062Z}]{stacks-project}, we have surjections of $\cO_C$-modules
     \begin{equation*}
          \operatorname{Sym}^a(I_C/I_C^2)\twoheadrightarrow I_C^a/I_C^{a+1}\twoheadrightarrow I_C^a\cO_Z/I_C^{a+1}\cO_Z.
     \end{equation*}
     Global generation descends to symmetric powers and quotients, and hence \(H^1(C,I_C^a\cO_Z/I_C^{a+1}\cO_Z)\) vanishes, using that $C$ is $1$-dimensional. The degree zero part is \(\cO_Z/I_C\cO_Z\cong \mathcal O_C\), so additivity and $H^1(X,\cO_C)=0$ give \(\chi(\mathcal O_Z)\geq1\) on every nonempty connected component.
     
     Consider a stable pair $[\cO_X\xrightarrow{s}\cF]$ of class $\beta$ as specified above. Then, by \cite[(1.7)]{pt07}, there is an exact sequence
     \begin{equation*}
          0\to \cO_Z\to \cF\to \cQ \to 0,
     \end{equation*}
     where $Z$ is a Cohen--Macaulay subscheme of dimension $1$ and class $[Z]=\beta$, and $\cQ$ is $0$-dimensional. Hence, we get
     $\chi(\cF)=\chi(\mathcal O_Z)+\operatorname{length}(\cQ)>0$ by the first part of this lemma and the assumption on $\beta$. This shows \(P_{\leq 0}(X,\beta)=\emptyset\).
\end{proof}

\subsubsection{The GW/PT calculation for \texorpdfstring{\([r_m]\)}{[r(m)]}}
\label{subsubsec:gwpt-rm}

\begin{proposition}
\label{prop:gwpt-rm}
For every \(m\geq1\), one has
\(L_{\Phi,0}(\widetilde r_m)=L_{\Psi,0}(\widetilde r_m)=0\), and hence
\(\delta_0([r_m])=0\).
\end{proposition}

\begin{proof}
We first consider the PT theory of $X_\square$ in class $m\beta$,
where $\beta$ is the fiber class, as in \S\ref{rcl:rational-surface-class}. If $Z$ is a proper
one-dimensional subscheme with $[Z]=m\beta$, its reduced support $C$
is a disjoint union of distinct fibers of $\pi:S\to B$. On each component
$C_w\cong\mathbb P^1$, \eqref{eq:rm-square-normal-bundle} gives
\[
N_{C_w/X_\square}
 \cong \mathcal O_{C_w}\oplus
        \mathcal O_{C_w}(-1)^{\oplus2}.
\]
Thus, the conormal sheaf of $C$ is globally generated, and
$H^1(C,\mathcal O_C)=0$. Lemma~\ref{num:lemma:PT-vanishing-global-gen-C-lemma} above now implies
$P_{\leqslant0}(X_\square,m\beta)=\emptyset$. From $\Sigma^{-1}(Y_m)=\sqcup_{a+b=m}Y_a\times Y_b$, we can derive that $\sum_{m\geq 0} \mathbf{1}_{Y_m}Q^m$ is group-like, so that $L_{\Psi,0}(\sum_{m\geq 0} \mathbf{1}_{Y_m}Q^m)=0$ by Lemma~\ref{uce:lemma:conv-log-eval}. Linearity $L_{\Psi,0}$ now implies $L_{\Psi,0}(\widetilde r_m)=0$.

We next consider GW theory. For use throughout this computation, set $x=c_1(\cO_C(1))$, $h=c_1(\cO_W(1))$, and $P=xh$, which satisfy $x^2=h^2=0$ and $\int_S P =1$. Since the cohomology of $S$ and $X_\square$ are isomorphic, we denote the pullbacks of these to $X_\square$ by the same symbols throughout. Klemm--Pandharipande
\cite[\S3.3, Prop.~3]{KP}, with $q_1=Q$ and $q_2=0$ in their
notation, compute
\begin{equation}\label{num:eq:kp-X-square-comp}
\sum_{d\geqslant1}\mathrm{GW}_{1,d\beta}(X_\square)Q^d
 =\frac1{12}\log\frac1{1-Q}
 =\sum_{d\geqslant1}\frac{Q^d}{12d}.
\end{equation}
and
\begin{equation*}
     \mathrm{GW}_{0,d\beta}(P)=\frac1{d^2}.
\end{equation*}
Inserting the second equality into the genus $0$ formula \eqref{eq:gv-kp-definitions}, with one insertion,
therefore implies $n_{0,\beta}(P)=1$ and $n_{0,d\beta}(P)=0$ for $d>1$.

Restricted to $S$, the Chern roots of $TX_\square$ are
$2x,2h,-x-h,-x-h$. Hence, $c_2(TX_\square)|_S=-2P$, which, using the preceding calculation, gives
\[
n_{0,d\beta}\bigl(c_2(TX_\square)\bigr)
 =-2\delta_{d,1}.
\]

As $n_{0,d\beta}(P)=\delta_{d,1}$, the genus $0$ BPS invariants satisfy
$n^{X_\square}_{0,d\beta}(\gamma)=\delta_{d,1}\int_S\gamma|_S$
for every $\gamma\in H^4(X_\square;\bQ)=\bQ\cdot P$.
The self-intersection of $[S]$ is
\[
\int_S c_2(N_{S/X_\square})
 =\int_S(x+h)^2=2.
\]
Thus, the incidence term in the  rules \eqref{eq:universal-meeting-recursion} is $2$
for $(\beta,\beta)$ and is zero for all other pairs
$(a\beta,b\beta)$. One has
$m_{\beta,\beta}=2-2=0$. Using both rules in \eqref{app:eq:kp-meeting-rules} and inducting
on $a+b$, we obtain
$m_{a\beta,b\beta}=0$ for all $a,b\geqslant1$.

Consequently, the meeting contribution in \eqref{eq:gv-kp-definitions} vanishes, while its
$c_2$ contribution is
\[
\frac1{24}n_{0,\beta}\bigl(c_2(TX_\square)\bigr)
       \log(1-Q)
 =\frac1{12}\log\frac1{1-Q},
\]
which is the entire genus one series \eqref{num:eq:kp-X-square-comp} above. Uniqueness in the
genus one transformation \eqref{eq:gv-kp-definitions} now gives
$n_{1,d\beta}=0$ for every $d\geqslant1$. Hence
$L_{\Phi,0}(\widetilde r_m)=0$. As both logarithmic functionals
vanish on $\mathfrak m^2$, we conclude that
$\delta_0([r_m])=0$.
\end{proof}

\subsubsection{The GW/PT calculation for \texorpdfstring{\([p_m]\)}{[p(m)]}}
\label{subsubsec:gwpt-pm}

For the calculation below, we compare the representative $\tilde{p}_m$ of \eqref{eq:pm-bracket-definition} with
Pardon's equivariant local curve element.  Let
\(E_0=\mathcal O_{\mathbb P^1}\oplus
\mathcal O_{\mathbb P^1}(-1)^{\oplus2}\) and
\(X_{\mathrm{rat}}=\operatorname{Tot}_{\mathbb P^1}(E_0)\), and let
\(\mathbb C^*\) act on the three summands with weights \(1,1,-2\).
Pardon's calculation with the elementary symmetric functions on
\(\operatorname{Sym}^m\mathbb A^1\)
\cite[Def.~4.7, Prop.~4.8, and Lemma~4.14]{pardon2023} shows
that the connected restriction of \(x_{0,m}^{[1]}\)\footnote{Pardon calls this $\prescript{}{1}{x_{0,m}}$.} is \(p_m\).
Therefore,
\(x_{0,m}^{[1]}\equiv\widetilde p_m\pmod{\mathfrak m^2}\).

\begin{remark} We spell out the use of the Borel construction in this comparison. The elementary symmetric coordinates identify \(\operatorname{Sym}^m\mathbb A^1\) with the sum of the representations of weights \(1,\ldots,m\). Taking the finite approximation over \(B=\mathbb P^{m-1}\), with the convention \(t\mapsto h=c_1(\mathcal O_B(1))\), gives the bundle \(V_m=\bigoplus_{j=1}^m\mathcal O_B(j)\). If \(U_m\) denotes its compactly supported Thom class, then \(U_m/m!\) represents the class \(x_{0,m}^{[1]}\) of the bivariant degree \(2m-2(m-1)=2\).

The connected diagonal is the total space of \(\mathcal O_B(1)\), and its normal bundle is \(\bigoplus_{j=2}^m\mathcal O_B(j)\). The Thom self-intersection formula therefore identifies the restriction of \(U_m/m!\) with \(h^{m-1}\) times the Thom class of this line bundle. Since \(\int_Bh^{m-1}=1\), this gives the point class \(p_m\).

This explains how to compare equivariant invariants from \cite{MCT} with our non-equivariant invariants. Indeed, for either calculation below, write the raw equivariant invariant as \(I_m(t)=c_m/t\). The normalized Thom class restricts equivariantly to \(t^m\) on the zero section. Compatibility of the virtual classes with the Borel construction consequently gives
$$
\int_{\mathbb P^{m-1}}
\left.\bigl(t^m I_m(t)\bigr)\right|_{t=h}
=
\int_{\mathbb P^{m-1}}c_mh^{m-1}
=c_m.
$$
Thus the cancellation of \(t^{-1}\) computes an ordinary pairing in the family over \(B\).
\end{remark}

\begin{proposition}
\label{prop:gwpt-pm}
For every \(m\geq1\), one has \(\delta_1([p_m])=0\).  More precisely,
\begin{equation}
 L_{\Phi,1}(\widetilde p_m)=L_{\Psi,1}(\widetilde p_m)
 =\begin{cases}1,&m=1,\\0,&m>1.\end{cases}
 \label{eq:pm-common-value}
\end{equation}
\end{proposition}

\begin{proof}
Both logarithmic functionals vanish on \(\mathfrak m^2\), so we may
replace \(\widetilde p_m\) by \(x_{0,m}^{[1]}\).  We calculate on
\(X_{\mathrm{rat}}\).  Restrict the Calabi--Yau torus $T_0$ of dimension two
from \cite[\S4.2]{MCT} along the cocharacter with weights
\(1,1,-2\), and let \(t\) be its equivariant parameter.  The trivial
summand of \(E_0\), which is the deformation direction of the zero
section \(C_0\), has weight \(t\).

Lemma~\ref{num:lemma:PT-vanishing-global-gen-C-lemma} shows that every
nonempty stable pair on $X_{\mathrm{rat}}$ has positive charge, giving us $\Psi_{< 0}(x_{0,m}^{[1]})=0$ and $\Psi_0(x_{0,m}^{[1]})=1$. Cao--Maulik--Toda compute in \cite[Prop.~4.2]{MCT} that
\(\PT^{T_0}_{1,[C_0]}(X_{\mathrm{rat}})=t^{-1}\), and that
\(\PT^{T_0}_{1,m[C_0]}(X_{\mathrm{rat}})\) vanishes for \(m>1\). Moreover, the Thom class of the deformation space $H^0(C_0, \mathcal O_{C_0})$
cancels \(t^{-1}\); see \cite[Lemma~4.14]{pardon2023} and the remark above. This gives us $\Psi_1(x_{0,m}^{[1]})=\delta_{m,1}$ and hence
\begin{equation}
 \sum_{m\geq1}L_{\Psi,1}(\widetilde p_m)Q^m=Q
 \label{eq:pm-pt-log-series}
\end{equation}
can be read off from the definition of the convolution logarithm.

For $X_{\mathrm{rat}}$, the first case in the proof of
\cite[Cor.~4.3]{MCT}, using the Aspinwall--Morrison formula, yields \(\GW_{0,m}^{\mathrm{raw}} = m^{-3}t^{-1}\). We therefore get
$\sum_{m\geq1}\langle\mathbf{GW}_0^{\mathrm{raw}},\tilde{p}_m\rangle Q^m=\sum_{m\geq1}m^{-3}Q^m$, where we recognize the right-hand side as the third polylogarithm series $\operatorname{Li}_3(Q)$. We can then write out \eqref{eq:primitive-genus-zero-classes} to get
\begin{equation}\label{eq:pm-polylog-transform}
     \operatorname{Li}_3(Q)= \sum_{m\geq1}\langle\mathbf{GW}_0^{\mathrm{raw}},\tilde{p}_m\rangle Q^m = \sum_{m\geq 1} \langle\bfN_0,\tilde{p}_m\rangle \operatorname{Li}_3(Q^m),
\end{equation}
where we recognize the same polylogarithm inside the genus $0$ transformation. It follows that $\langle\bfN_0,\tilde{p}_m\rangle=\delta_{m,1}$. By Lemma~\ref{lem:normalized-theories-characters}, the same is true of
the \(L_{\Phi,1}(\widetilde p_m)=\delta_{m,1}\).  Together with
\eqref{eq:pm-pt-log-series}, this proves the proposition.
\end{proof}

\subsubsection{The GW/PT calculation for \texorpdfstring{\([e_{m,(E,N)}]\)}{[e(m,E,N)]}}\label{subsubsec:gwpt-em}
For \(r\geq1\), put
\(\sigma(r)=\sum_{d\mid r}d\) and
\(\sigma_2(r)=\sum_{d\mid r}d^2\), and let
\begin{equation}
 M(Q):=\prod_{k\geq1}(1-Q^k)^{-k}
 \label{eq:em-macmahon}
\end{equation}
be the MacMahon function.

\begin{proposition}
\label{prop:gwpt-em}
For every \(m\geq1\) and every regular elliptic normal datum
\((E,N)\), one has
\begin{equation}
 L_{\Phi,0}(\widetilde e_{m,(E,N)})
 =L_{\Psi,0}(\widetilde e_{m,(E,N)})
 =\frac{\sigma_2(m)}m,
 \qquad
 \delta_0([e_{m,(E,N)}])=0.
 \label{eq:em-common-logarithmic-value}
\end{equation}
Equivalently, 
\begin{equation}
 1+\sum_{m\geq1}\Phi(\widetilde e_{m,(E,N)})Q^m
 =1+\sum_{m\geq1}\Psi(\widetilde e_{m,(E,N)})Q^m
 =M(Q).
 \label{eq:em-common-total-series}
\end{equation}
\end{proposition}

\begin{proof}
Consider the setup in \S\ref{subsec:elliptic-point-classes}. By deformation invariance, it is enough to use the split representative
\(N^\circ\), so the cycle space under consideration is compact by Lemma~\ref{lem:em-cycle-sector-point}.  Cao--Maulik--Toda prove
\cite[Lemma~4.9]{MCT} that the fixed stable pairs of charge zero are
the quotients
\begin{equation}
 \mathcal O_{X_{\mathrm{ell}}}\twoheadrightarrow
 \bigoplus_{(i_1,i_2,i_3)\in\Delta}
 L_1^{-i_1}\otimes L_2^{-i_2}\otimes L_3^{-i_3},
 \label{eq:em-plane-partition-pair}
\end{equation}
where \(\Delta\subset\mathbb Z_{\geq0}^3\) is a plane partition of
cardinality \(m\).  The genericity of the \(L_i\) in Lemma~\ref{lem:em-regular-normal-form} makes every
nontrivial tensor product of them acyclic.  The remaining terms are
those in the local calculation of \cite[Theorem~4.10]{MCT}. Thus, every plane partition contributes \(1\), and
\[
 1+\sum_{m\geq1}\PT_{0,m[E]}(X_{\mathrm{ell}})Q^m=M(Q).
\]
On this cycle sector, since $\sum_{m\geq 0}\mathbf{1}_{m[E]}Q^m$ is group-like, Lemma~\ref{uce:lemma:conv-log-eval} gives us
\begin{equation}
 \sum_{m\geq1}L_{\Psi,0}(\widetilde e_m^\circ)Q^m=\log M(Q)
 =\sum_{k,r\geq1}\frac{k}{r}Q^{kr}
 =\sum_{m\geq1}\frac{\sigma_2(m)}mQ^m,
 \label{eq:em-log-macmahon}
\end{equation}
where the last two equations follows from expanding the logarithm.

For the GW theory, let \(f:E'\to E\) be a connected cover of degree
\(m\) from a curve of genus one.  It is \'etale, and each \(f^*L_i\)
is a nontrivial line bundle of degree zero.  Thus, by generality of the $L_i$, we get
\(H^0(E',f^*N^\circ)=H^1(E',f^*N^\circ)=0\).  There are
\(\sigma(m)\) index \(m\) sublattices of the lattice of \(E\), and
each corresponding cover has \(m\) deck transformations, which is how
Klemm--Pandharipande obtain \cite[\S2.1]{KP}
\begin{equation}
 \sum_{m\geq1}
 \left\langle\mathbf{GW}_{1}^{\mathrm{raw}},
 \widetilde e_m^\circ\right\rangle Q^m
 =\sum_{m\geq1}\GW_{1,m[E]} Q^m
 =\sum_{m\geq1}\frac{\sigma(m)}mQ^m.
 \label{eq:em-raw-gw-series}
\end{equation}
Since, we have no rational curves in this geometry, the genus $1$ part of \eqref{eq:gv-kp-definitions} identifies the left-hand side with $\sum_{m\geq 1}n_{1,m[E]}\sum_{d\geq 1}\frac{\sigma(d)}{d}Q^{md}$, so that \(n_{1,1}=1\) and
\(n_{1,k}=0\) for \(k>1\). By Proposition~\ref{uce:prop:transformed-theory-coefficients} and Lemma~\ref{lem:normalized-theories-characters}, we get $L_{\Psi,0}(\sum_{m\geq 0}\widetilde e_{m}^\circ Q^m)=\log M(Q)$. Comparison with \eqref{eq:em-log-macmahon} proves
\eqref{eq:em-common-logarithmic-value} for the split representative,
and deformation invariance proves it for \((E,N)\).  Exponentiating, we get
 \eqref{eq:em-common-total-series}.
\end{proof}

\subsubsection{The GW/PT calculation for
\texorpdfstring{\([\ell_{d_1,d_2}]\)}{[ell(d1,d2)]}}
\label{subsubsec:gwpt-toric-cap}

\begin{proposition}
\label{prop:gwpt-toric-cap}
For every \(d_1,d_2\geq1\), one has
\[
 L_{\Phi,0}(\widetilde\ell_{d_1,d_2})
 =L_{\Psi,0}(\widetilde\ell_{d_1,d_2})=0,
 \qquad
 \delta_0([\ell_{d_1,d_2}])=0.
\]
\end{proposition}

\begin{proof}
We first prove the PT vanishing.  The reduced support of a connected
cycle in \eqref{eq:cap-entire-connected-sector} is either a smooth
fiber \(F_t\cong\mathbb P^1\) or the nodal curve \(C=A\cup E\).
Their conormal bundles in \(X_0\) are
\begin{equation}
 I_{F_t}/I_{F_t}^2
 \cong\mathcal O_{F_t}\oplus\mathcal O_{F_t}(1)^{\oplus2},
 \qquad
 I_C/I_C^2
 \cong\mathcal O_C\oplus\mathcal O_C(1,0)\oplus\mathcal O_C(0,1).
 \label{eq:cap-conormal-bundles}
\end{equation}
Here \(\mathcal O_C(a,b)\) has degrees \(a\) and \(b\) on \(A\) and
\(E\), respectively.  Formula \eqref{eq:cap-conormal-bundles} follows
from \eqref{eq:cap-special-fiber-restrictions}. In particular, both conormal bundles
are globally generated with \(H^1(C,\mathcal O_C)=0\). Lemma~\ref{num:lemma:PT-vanishing-global-gen-C-lemma} above now implies
that there are no stable pairs of charge $\leq 0$ on the connected locus.
Since the convolution logarithm is infinitesimal, its pairing depends
only on the connected restriction $\ell^\con_{d_1,d_2}$ of the insertion. From the definition of the convolution logarithm, we can then read off \(L_{\Psi,0}(\widetilde\ell_{d_1,d_2})=0\).

For the GW theory, put \(\beta_1=[A]\) and \(\beta_2=[E]\).
Klemm--Pandharipande use precisely the pencil
\(\pi:S_0\to\mathbb P^1\) to determine the meeting correction. They compute
\cite[\S2.3.1, Prop.~1]{KP} the following contribution of $\pi$ to the GW theory
\begin{equation}
 \left\langle\mathbf{GW}_{1}^{\mathrm{raw}},
 \widetilde\ell_{d_1,d_2}\right\rangle
 =\operatorname{Cont}_{\pi}
 \bigl(\GW_{1,d_1\beta_1+d_2\beta_2}\bigr)
 =\frac{\delta_{d_1,d_2}}{12d_1},
 \label{eq:cap-kp-raw-number}
\end{equation}
where the equality with the pairing follows by equivariant localization: on
the fixed locus the insertion is \(h^N\), and
\(\int_{\mathbb P^N}h^N=1\), see the remark below.

Since \(X_0\) retracts onto \(S_0\) and \(H^4(S_0;\mathbb Q)=\mathbb Q\), we can compute
\[
 \int_{S_0}c_2(TX_0)
 =\int_{S_0}\bigl(c_2(TS_0)-c_1(TS_0)^2
   +c_1(L_1)c_1(L_2)\bigr)=5-7+2=0
\]
via Chern roots. Thus, \(c_2(TX_0)=0\), so that there is no contribution from the term involving \(c_2(TX_0)\) in \eqref{eq:gv-kp-definitions} and \eqref{eq:universal-gw-transformation}. Since \(A\) and \(E\) meet transversely once, the only meeting numbers in this geometry are \(m_{\beta_1,\beta_2}=m_{\beta_2,\beta_1}=1\), and we can write out the pairing of \eqref{uce:eq:n1-def-rewriting} with the $\widetilde{\ell}_{d_1,d_2}$ to get
\begin{equation*}
     \begin{aligned}
          &\sum_{(d_1,d_2)>0}\langle \bfN_1,\widetilde{\ell}_{d_1,d_2}\rangle\sum_{k>0}\frac{\sigma(k)}{k}Q_1^{kd_1}Q_2^{kd_2}-\frac{1}{12}\log(1-Q_1Q_2)\\
          =&\sum_{(d_1,d_2)>0}\langle \mathbf{GW}^{\mathrm{raw}}_1,\widetilde{\ell}_{d_1,d_2}\rangle Q_1^{d_1}Q_2^{d_2}
          =\sum_{d>0}\frac{1}{12d}Q_1^d Q_2^d = -\frac{1}{12}\log(1-Q_1Q_2),
     \end{aligned}
\end{equation*}
where the second equality uses \eqref{eq:cap-kp-raw-number} and the last equality can be found by expanding the logarithm. Hence, $\langle \bfN_1,\widetilde{\ell}_{d_1,d_2}\rangle=0$ for all $d_1$, $d_2$, and, using Lemma~\ref{lem:normalized-theories-characters}, it follows that
\(L_{\Phi,0}(\widetilde\ell_{d_1,d_2})=0\).  The two vanishings prove the
proposition.
\end{proof}

\subsection{GW/PT-correspondence and Gopakumar--Vafa integrality}
\label{sec:gv-integrality}

We will now specialize our  correspondence between connected theories to a fixed smooth projective
Calabi--Yau fourfold.  

\subsubsection{The GW/PT correspondence}
We now explain how the correspondence in Theorem~\ref{nap:thm:L-phi-psi-correspondence} yields the conjectures in \S\ref{intro:sec:conjectures}.

\begin{proposition}
\label{prop:gv-cmt-identities}
Using Proposition~\ref{prop:family-cy4-virtual-package}, the
identities \eqref{intro:eq:cmt-conj} hold for \(X\).
\end{proposition}

\begin{proof}
First note, that Lemma~\ref{uce:lemma:divisor-insertions} immediately reduces \eqref{intro:eq:cmt-conj} to the case where $m_i\geq 3$ for all $i$, since it implies $\PT_{1,\beta}(\gamma_\bullet)=\PT_{1,\beta}(\tilde{\gamma}_\bullet)\cdot\prod_j\int_\beta\delta_j$. So we may now assume that all insertions are in degree $m_i\geq 3$.

Secondly, every \(n<0\), the oriented
stable pair virtual class belongs to
\(H_{2n}(P_n(X,\beta);\mathbb Z)=0\)
\cite[Theorem~2.4]{MCT}. Hence, the specialization of PT theory to
\(X\) has no terms of negative charge. The transformed connected GW
class has only charges zero and one, so its exponential likewise has no
terms of negative charge. Hence, the discussion in \S\ref{uce:sec:log-insertion-evaluation} applies and \eqref{uce:eq:log-coefficient-identities} holds for both theories and we can prove the following.

So, we can evaluate at $(X,\beta)$ and $(X,\beta,\bullet_\gamma)$ following \S\ref{uce:sec:log-insertion-evaluation}. Then, combining \eqref{uce:eq:log-coefficient-identities}, Theorem~\ref{nap:thm:L-phi-psi-correspondence}, Lemma~\ref{lem:normalized-theories-characters}, and Proposition~\ref{uce:prop:transformed-theory-coefficients}, we obtain 
\begin{equation}\label{eq:gv-connected-identities}
     \begin{aligned}
          \frac{\sum_{\beta>0}\PT_{1,\beta}(\gamma_\bullet)Q^\beta}{\sum_{\beta\geq 0}\PT_{0,\beta}Q^\beta}&=\sum_{\beta>0}n_{0,\beta}(\gamma_\bullet)Q^\beta,\\ 
          \log\bigg(\sum_{\beta\geq 0}\PT_{0,\beta}Q^\beta\bigg)&=\sum_{\beta>0}n_{1,\beta}\log\big(M(Q^\beta)\big).
     \end{aligned}
\end{equation}
Multiplying the first identity by \(\sum_{\beta\geq 0}\PT_{0,\beta}Q^\beta\), and
exponentiating the second, gives \eqref{intro:eq:cmt-conj}.
\end{proof}

\subsubsection{}

\begin{example}\label{ex:mct-divisor-counterexample}
     We show here in an example that \cite[Conj. 0.1]{MCT} is inconsistent without our mild correction in \eqref{intro:eq:cmt-conj}. Consider the example in \cite[\S3.4]{MCT}: Let $Y$ be a general quintic $3$-fold, $E$ an elliptic curve and $X\coloneqq Y\times E$. Denote by $\alpha$ the class of a line on $Y$ and by $e$ the class of $E$ in $X$. Set
     \begin{equation*}
          \gamma\coloneqq\pr_Y^*H\cup \pr_E^*[\pt],\qquad D\coloneqq\pr_E^*[\pt],
     \end{equation*}
     where $H$ is the hyperplane class in $Y$. We have $2875$ lines on $Y$ so that $n_{0,\alpha}(\gamma)=2875$, since every line intersects $H$ once. By \cite[Lemma 3.14]{MCT}, we have $\PT_{0,e}=\chi(Y)=-200$. We have $H_2(X)=\bZ\alpha\oplus\bZ e$ using the Lefschetz hyperplane theorem for $Y$, and every genus $0$ map is constant in the $E$ direction. Hence, for $\beta=\alpha+e$, the only decomposition contributing on the right-hand side of \eqref{intro:eq:cmt-conj} is $\alpha+e$. 

     Now, using the divisor equation for stable pairs and $\int_{\alpha+e}D=1$, we obtain
     \begin{equation*}
          \PT_{1,\alpha+e}(\gamma,D)=1\cdot \PT_{1,\alpha+e}(\gamma)= n_{0,\alpha}(\gamma)\cdot \PT_{0,e}=-575000,
     \end{equation*}
     where the second equality is \cite[Conj. 0.1]{MCT} applied to $\gamma$, and the last uses the computations above. On the other hand, \cite[Conj. 0.1]{MCT} applied to $(\gamma, D)$ gives us 
     \begin{equation*}
          \PT_{1,\alpha+e}(\gamma,D)= n_{0,\alpha}(\gamma,D)\cdot \PT_{0,e}= \int_\alpha D\cdot n_{0,\alpha}(\gamma)\cdot \PT_{0,e} = 0,
     \end{equation*}
     since $\int_\alpha D$ vanishes. This is a contradiction, but is fixed by our small correction in the first formula in \eqref{intro:eq:cmt-conj}.
\end{example}

\subsubsection{Integrality}

\begin{corollary}
\label{cor:gv-integrality}
For every integral insertion datum \(\gamma_\bullet\) as above and every
\(\beta>0\),
\(n_{0,\beta}(\gamma_\bullet)\) and \(n_{1,\beta}\) are integers.  The
meeting numbers \(m_{\beta_1,\beta_2}\) are integers as well.
\end{corollary}

\begin{proof}
By Proposition~\ref{prop:family-cy4-virtual-package} (from \cite{karpov-bj}), the oriented stable pair virtual classes have integral coefficients (as in \cite{Borisov_2017}). Moreover,
\(\operatorname{ch}_3(\mathbb F)\) is the fundamental cycle of the
universal sheaf, so the classes \(\tau(\gamma_i)\) are integral
\cite[\S2.3]{MCT}. Hence,
\(\PT_{0,\beta}\) and \(\PT_{1,\beta}(\gamma_\bullet)\) are integers.
The series \(Z_{\mathrm{PT},0}\) has constant term one, and consequently
its inverse also has integral coefficients.  The first identity in
\eqref{eq:gv-connected-identities} hence implies that every
\(n_{0,\beta}(\gamma_\bullet)\) is an integer.

We next use a standard triangular argument.  We claim that every series \(F\) in $Q^\beta$ with $\beta$ ranging over effective curve classes, with
integral coefficients and constant term one has a unique expression
\(F=\prod_{\beta>0}M(Q^\beta)^{a_\beta}\) with
\(a_\beta\in\mathbb Z\).  Indeed, order the effective classes by
\(H\cdot\beta\).  Since \(M(t)=1+t+O(t^2)\), the coefficient of
\(Q^\beta\) is \(a_\beta\) plus an integral expression in exponents
already determined in smaller degree.  Thus induction gives
\(a_\beta\in\mathbb Z\).  Apply this to
\(F=Z_{\mathrm{PT},0}\) and compare with the second identity in
\eqref{intro:eq:cmt-conj}. Uniqueness gives
\(a_\beta=n_{1,\beta}\), which implies integrality of $n_{1,\beta}$.

Finally, choose an integral basis $S_i$ of
\(H^4(X,\mathbb Z)/\mathrm{tors}\).  Its intersection matrix $g_{ij}=\int_X S_iS_j$ is
unimodular by writing the Poincar\'e duality isomorphism of \cite[Prop. 3.3.8]{hatcher2002algebraic} in the chosen basis. Hence, its inverse $g^{ij}$ also has integral entries. The rules
\eqref{app:eq:kp-meeting-rules}, together with the integrality just proved
for insertions from \(H^4(X,\mathbb Z)\) and for \(c_2(TX)\), determine
the meeting numbers inductively as integers.
\end{proof}

\phantomsection
\addcontentsline{toc}{section}{References}

\begin{small}
\bibliographystyle{alphaurl}
\bibliography{paper}

@article{pardon2023,
      author = {{Pardon}, John},
      title = "{Universally counting curves in Calabi--Yau threefolds}",
      journal = {arXiv e-prints},
      year = 2023,
      month = aug,
      eid = {arXiv:2308.02948},
      pages = {arXiv:2308.02948},
      doi = {10.48550/arXiv.2308.02948},
      archivePrefix = {arXiv},
      eprint = {2308.02948},
      primaryClass = {math.AG},
      adsurl = {https://ui.adsabs.harvard.edu/abs/2023arXiv230802948P}
}

@article{Ionel2018,
  title = {The Gopakumar-Vafa formula for symplectic manifolds},
  volume = {187},
  ISSN = {0003-486X},
  url = {http://dx.doi.org/10.4007/annals.2018.187.1.1},
  DOI = {10.4007/annals.2018.187.1.1},
  number = {1},
  journal = {Annals of Mathematics},
  publisher = {Annals of Mathematics},
  author = {Ionel,  Eleny-Nicoleta and Parker,  Thomas},
  year = {2018},
  month = Jan 
}

@unpublished{Hatcher-VBKT,
  author = {Hatcher, Allen},
  title  = {Vector Bundles and {K}-Theory},
  note   = {Unpublished book, version 2.2},
  year   = {2017},
  url    = {https://pi.math.cornell.edu/~hatcher/VBKT/VBpage.html},
}

@book{Lee2012,
  title = {Introduction to Smooth Manifolds},
  ISBN = {9781441999825},
  ISSN = {0072-5285},
  url = {http://dx.doi.org/10.1007/978-1-4419-9982-5},
  DOI = {10.1007/978-1-4419-9982-5},
  journal = {Graduate Texts in Mathematics},
  publisher = {Springer New York},
  author = {Lee,  John M.},
  year = {2012}
}

@book{Goresky1988,
  title = {Stratified Morse Theory},
  ISBN = {9783642717147},
  url = {http://dx.doi.org/10.1007/978-3-642-71714-7},
  DOI = {10.1007/978-3-642-71714-7},
  publisher = {Springer Berlin Heidelberg},
  author = {Goresky,  Mark and MacPherson,  Robert},
  year = {1988}
}

@book{cartan1995elementary,
  title={Elementary Theory of Analytic Functions of One or Several Complex Variables},
  author={Cartan, Henri},
  year={1995},
  publisher={Dover Publications},
  address={New York},
  note={Reprint of the 1963 edition published by Addison-Wesley}
}

@article{Demailly1990,
  title = {Cohomology of q-convex spaces in top degrees},
  volume = {204},
  ISSN = {1432-1823},
  url = {http://dx.doi.org/10.1007/BF02570874},
  DOI = {10.1007/bf02570874},
  number = {1},
  journal = {Mathematische Zeitschrift},
  publisher = {Springer Science and Business Media LLC},
  author = {Demailly,  Jean-Pierre},
  year = {1990},
  month = Dec,
  pages = {283–295}
}

@book{Grauert1984,
  title = {Coherent Analytic Sheaves},
  ISBN = {9783642695827},
  ISSN = {0072-7830},
  url = {http://dx.doi.org/10.1007/978-3-642-69582-7},
  DOI = {10.1007/978-3-642-69582-7},
  journal = {Grundlehren der mathematischen Wissenschaften},
  publisher = {Springer Berlin Heidelberg},
  author = {Grauert,  Hans and Remmert,  Reinhold},
  year = {1984}
}

@article{Siu1976,
  title = {Every Stein subvariety admits a Stein neighborhood},
  volume = {38},
  ISSN = {1432-1297},
  url = {http://dx.doi.org/10.1007/BF01390170},
  DOI = {10.1007/bf01390170},
  number = {1},
  journal = {Inventiones Mathematicae},
  publisher = {Springer Science and Business Media LLC},
  author = {Siu,  Yum-Tong},
  year = {1976},
  month = feb,
  pages = {89–100}
}

@article{OzbagciPopescuPampu16,
  author  = {{\"O}zba{\u g}c{\i}, Burak and Popescu-Pampu, Patrick},
  title   = {Generalized plumbings and {Murasugi} sums},
  journal = {Arnold Mathematical Journal},
  volume  = {2},
  number  = {1},
  pages   = {69--119},
  year    = {2016},
  doi     = {10.1007/s40598-015-0033-3}
}

@misc{Guadagni18,
  author        = {Guadagni, Roberta},
  title         = {Symplectic neighborhood of crossing divisors},
  year          = {2018},
  eprint        = {1611.02363},
  archivePrefix = {arXiv},
  primaryClass  = {math.SG}
}

@article{Mumford61,
  author  = {Mumford, David},
  title   = {The topology of normal singularities of an algebraic
             surface and a criterion for simplicity},
  journal = {Publications Math{\'e}matiques de l'IH{\'E}S},
  volume  = {9},
  pages   = {5--22},
  year    = {1961},
  doi     = {10.1007/BF02698717}
}

@article{JU,
      author = {{Joyce}, Dominic and {Upmeier}, Markus},
      title = "{Bordism categories and orientations of moduli spaces}",
      journal = {arXiv e-prints},
      year = {2025},
      month = mar,
      eid = {arXiv:2503.20456},
      pages = {arXiv:2503.20456},
      doi = {10.48550/arXiv.2503.20456},
      archivePrefix = {arXiv},
      eprint = {2503.20456},
      primaryClass = {math.AT},
}

@article{cao-local-p2,
       author = {{Cao}, Yalong},
        title = "{A proof of GW/DT4 conjecture on local projective plane}",
      journal = {arXiv e-prints},
         year = 2026,
        month = oct,
          eid = {arXiv:2610.06024},
        pages = {arXiv:2610.06024},
archivePrefix = {arXiv},
       eprint = {2610.06024},
 primaryClass = {math.AG},
       adsurl = {https://ui.adsabs.harvard.edu/abs/2026arXiv261006024C}
}

@article{wendl-symplectic-lectures,
       author = {{Wendl}, Chris},
        title = "{Lectures on Holomorphic Curves in Symplectic and Contact Geometry}",
      journal = {arXiv e-prints},
         year = 2010,
        month = nov,
          eid = {arXiv:1011.1690},
        pages = {arXiv:1011.1690},
          doi = {10.48550/arXiv.1011.1690},
archivePrefix = {arXiv},
       eprint = {1011.1690},
 primaryClass = {math.SG},
       adsurl = {https://ui.adsabs.harvard.edu/abs/2010arXiv1011.1690W}
}

@article{Smale1965,
  title = {An Infinite Dimensional Version of Sard’s Theorem},
  volume = {87},
  ISSN = {0002-9327},
  url = {http://dx.doi.org/10.2307/2373250},
  DOI = {10.2307/2373250},
  number = {4},
  journal = {American Journal of Mathematics},
  publisher = {JSTOR},
  author = {Smale,  S.},
  year = {1965},
  month = Oct,
  pages = {861}
}

@inbook{Kontsevich1995,
  title = {Enumeration of Rational Curves Via Torus Actions},
  ISBN = {9781461242642},
  url = {http://dx.doi.org/10.1007/978-1-4612-4264-2_12},
  DOI = {10.1007/978-1-4612-4264-2_12},
  booktitle = {The Moduli Space of Curves},
  publisher = {Birkh\"{a}user Boston},
  author = {Kontsevich,  Maxim},
  year = {1995},
  pages = {335–368}
}

@article{Behrend-Manin-Mgn,
  title = {Stacks of stable maps and Gromov-Witten invariants},
  volume = {85},
  ISSN = {0012-7094},
  url = {http://dx.doi.org/10.1215/S0012-7094-96-08501-4},
  DOI = {10.1215/s0012-7094-96-08501-4},
  number = {1},
  journal = {Duke Mathematical Journal},
  publisher = {Duke University Press},
  author = {Behrend,  K. and Manin,  Yu.},
  year = {1996},
  month = Oct 
}

@article{Bojko2021,
  author  = {Bojko, Arkadij},
  title   = {Orientations for {DT} invariants on quasi-projective {Calabi--Yau} 4-folds},
  journal = {Advances in Mathematics},
  volume  = {388},
  year    = {2021},
  pages   = {107859},
  doi     = {10.1016/j.aim.2021.107859},
  eprint  = {2008.08441},
  archivePrefix = {arXiv},
  primaryClass  = {math.AG}
}

@book{hatcher2002algebraic,
  author    = {Hatcher, Allen},
  title     = {Algebraic Topology},
  publisher = {Cambridge University Press},
  address   = {Cambridge},
  year      = {2002},
  pages     = {xii+544},
  isbn      = {0-521-79540-0},
  mrnumber  = {1867354}
}

@misc{kr, 
    title={Proof of a Magnificent Conjecture},
    author={Kool, Martijn and Rennemo, Jørgen Vold}, 
    note={In Preparation}
}

@article{bf97_intrinsic_normal_cone,
	doi = {10.1007/s002220050136},
	url = {https://doi.org/10.1007%2Fs002220050136},
	year = 1997,
	month = {3},
	publisher = {Springer Science and Business Media {LLC}},
	volume = {128},
	number = {1},
	pages = {45--88},
	author = {K. Behrend and B. Fantechi},
	title = {The intrinsic normal cone},
	journal = {Inventiones Mathematicae}
}

@article{Behrend1997GW,
  author        = {Behrend, Kai},
  title         = {Gromov--{W}itten invariants in algebraic geometry},
  journal       = {Inventiones Mathematicae},
  volume        = {127},
  number        = {3},
  pages         = {601--617},
  year          = {1997},
  doi           = {10.1007/s002220050132},
  eprint        = {alg-geom/9601011},
  archivePrefix = {arXiv}
}

@misc{karpov-bj,
      author = {{Karpov}, Ivan},
      title = "{GW/PT-correspondence for CY4-folds II: Virtual Classes}",
      year = {2026}
}

@article{KP,
  author  = {Klemm, Albrecht and Pandharipande, Rahul},
  title   = {Enumerative geometry of {Calabi--Yau} 4-folds},
  journal = {Communications in Mathematical Physics},
  volume  = {281},
  year    = {2008},
  pages   = {621--653},
  doi     = {10.1007/s00220-008-0490-9},
  eprint  = {math/0702189},
  archivePrefix = {arXiv}
}

@article{Ati57,
  author  = {Atiyah, Michael F.},
  title   = {Vector bundles over an elliptic curve},
  journal = {Proceedings of the London Mathematical Society},
  series  = {3},
  volume  = {7},
  year    = {1957},
  pages   = {414--452},
  doi     = {10.1112/plms/s3-7.1.414}
}

@article{MCT,
  author  = {Cao, Yalong and Maulik, Davesh and Toda, Yukinobu},
  title   = {Stable pairs and {Gopakumar--Vafa} type invariants for {Calabi--Yau} 4-folds},
  journal = {Journal of the European Mathematical Society},
  volume  = {24},
  number  = {2},
  year    = {2022},
  pages   = {527--581},
  doi     = {10.4171/JEMS/1110},
  eprint  = {1902.00003},
  archivePrefix = {arXiv}
}

@book{Steenrod1951,
 ISBN = {9780691005485},
 URL = {http://www.jstor.org/stable/j.ctt1bpm9t5},
 author = {Steenrod, Norman E.},
 publisher = {Princeton University Press},
 title = {The Topology of Fibre Bundles. (PMS-14)},
 urldate = {2026-09-05},
 year = {1951}
}

@misc{bryan-thimm,
      author={Jim Bryan and Felix Thimm},
      title={--\,},
      year={},
      note ={In preparation},
      url={},
}

@misc{karpov-cy4-orientations,
      title={{GW/PT-correspondence for CY 4-folds-I: Canonical Orientations}}, 
      author={Ivan Karpov},
      year={2026},
      eprint={2610.02061},
      archivePrefix={arXiv},
      primaryClass={math.AG},
      url={https://arxiv.org/abs/2610.02061}, 
}

@incollection{Tian1987,
  author    = {Tian, Gang},
  title     = {Smoothness of the universal deformation space of compact {C}alabi-{Y}au manifolds and its {P}etersson-{W}eil metric},
  booktitle = {Mathematical Aspects of String Theory},
  editor    = {Yau, Shing-Tung},
  series    = {Advanced Series in Mathematical Physics},
  volume    = {1},
  pages     = {629--646},
  publisher = {World Scientific Publishing},
  address   = {Singapore},
  year      = {1987}
}

@article{Todorov1989,
  author  = {Todorov, Andrey N.},
  title   = {The {W}eil-{P}etersson geometry of the moduli space of {SU}($n \geq 3$) ({C}alabi-{Y}au) manifolds. {I}},
  journal = {Communications in Mathematical Physics},
  volume  = {126},
  number  = {2},
  pages   = {325--346},
  year    = {1989},
  doi     = {10.1007/BF02125128}
}

@article{Bogomolov1974,
  author  = {Bogomolov, Fedor A.},
  title   = {K\"ahler manifolds with trivial canonical class},
  journal = {Mathematics of the {USSR}-{I}zvestija},
  volume  = {8},
  number  = {1},
  pages   = {9--20},
  year    = {1974}
}

@article{siu-analytic-sheaf-cohom,
 ISSN = {00029947},
 URL = {http://www.jstor.org/stable/1995234},
 author = {Yum-Tong Siu},
 journal = {Transactions of the American Mathematical Society},
 pages = {77--94},
 publisher = {American Mathematical Society},
 title = {Analytic Sheaf Cohomology Groups of Dimension n of n-Dimensional Complex Spaces},
 urldate = {2026-08-21},
 volume = {143},
 year = {1969}
}

@article{Borisov_2017,
	doi = {10.2140/gt.2017.21.3231},
	url = {https://doi.org/10.2140%2Fgt.2017.21.3231},
	year = 2017,
	month = {8},
	publisher = {Mathematical Sciences Publishers},
	volume = {21},
	number = {6},
	pages = {3231--3311},
	author = {Dennis Borisov and Dominic Joyce},
	title = {Virtual fundamental classes for moduli spaces of sheaves on Calabi{\textendash}Yau four-folds},
	journal = {Geometry {\&} Topology}
}

@article{oh2022counting,
      author = {Jeongseok Oh and Richard P. Thomas},
      title = {{Counting sheaves on Calabi–Yau 4-folds, I}},
      volume = {172},
      journal = {Duke Mathematical Journal},
      number = {7},
      publisher = {Duke University Press},
      pages = {1333 -- 1409},
      year = {2023},
      doi = {10.1215/00127094-2022-0059},
      URL = {https://doi.org/10.1215/00127094-2022-0059}
}

@misc{stacks-project,
  author       = {The {Stacks project authors}},
  title        = {The Stacks project},
  howpublished = {\url{https://stacks.math.columbia.edu}},
  year         = {2026},
}

@article{BryanPandharipande2008,
  author        = {Bryan, Jim and Pandharipande, Rahul},
  title         = {The local {Gromov--Witten} theory of curves},
  journal       = {J. Amer. Math. Soc.},
  volume        = {21},
  number        = {1},
  pages         = {101--136},
  year          = {2008},
  doi           = {10.1090/S0894-0347-06-00545-5},
  mrnumber      = {2350052},
  eprint        = {math/0411037},
  archivePrefix = {arXiv},
  primaryClass  = {math.AG}
}

@article{Smale1959,
  author  = {Smale, Stephen},
  title   = {Diffeomorphisms of the {$2$}-sphere},
  journal = {Proc. Amer. Math. Soc.},
  volume  = {10},
  number  = {4},
  year    = {1959},
  pages   = {621--626},
  doi     = {10.1090/S0002-9939-1959-0112149-8},
  issn    = {0002-9939},
  mrnumber = {0112149}
}

@article{pt07,
       author = {{Pandharipande}, R. and {Thomas}, R.~P.},
        title = "{Curve counting via stable pairs in the derived category}",
      journal = {Inventiones Mathematicae},
         year = 2009,
        month = may,
       volume = {178},
       number = {2},
        pages = {407-447},
          doi = {10.1007/s00222-009-0203-9},
archivePrefix = {arXiv},
       eprint = {0707.2348},
 primaryClass = {math.AG},
       adsurl = {https://ui.adsabs.harvard.edu/abs/2009InMat.178..407P}
}

@article{PandharipandeThomas2010,
  author        = {Pandharipande, Rahul and Thomas, Richard P.},
  title         = {Stable pairs and {BPS} invariants},
  journal       = {J. Amer. Math. Soc.},
  volume        = {23},
  number        = {1},
  pages         = {267--297},
  year          = {2010},
  doi           = {10.1090/S0894-0347-09-00646-8},
  mrnumber      = {2552254},
  eprint        = {0711.3899},
  archivePrefix = {arXiv},
  primaryClass  = {math.AG}
}

@article{OkounkovPandharipande2010,
  author        = {Okounkov, Andrei and Pandharipande, Rahul},
  title         = {The local {Donaldson--Thomas} theory of curves},
  journal       = {Geom. Topol.},
  volume        = {14},
  number        = {3},
  pages         = {1503--1567},
  year          = {2010},
  doi           = {10.2140/gt.2010.14.1503},
  mrnumber      = {2679579},
  eprint        = {math/0512573},
  archivePrefix = {arXiv},
  primaryClass  = {math.AG}
}
\end{small}

\end{document}